\documentclass{article}

\usepackage{float}
\PassOptionsToPackage{numbers,sort&compress}{natbib}
\usepackage[main,final]{neurips_2026}

\usepackage[utf8]{inputenc} %
\usepackage[T1]{fontenc}    %
\usepackage{hyperref}       %
\usepackage{url}            %
\usepackage{booktabs}       %
\usepackage{amsfonts}       %
\usepackage{nicefrac}       %
\usepackage{microtype}      %
\usepackage{xcolor}         %
\usepackage{amsmath,amssymb,amsfonts,amsthm}
\usepackage{mathtools}
\usepackage{bm}
\usepackage{graphicx}
\graphicspath{{./}}
\usepackage{booktabs}
\usepackage{enumitem}
\usepackage[dvipsnames]{xcolor}
\usepackage{hyperref}   %
\hypersetup{hidelinks}
\usepackage[nameinlink, capitalize, noabbrev]{cleveref}
\usepackage{tikz}
\usepackage{mathrsfs}  
\usepackage{verbatim}
\usepackage{multirow}

\newtheorem{theorem}{Theorem}
\numberwithin{theorem}{section}
\newtheorem{lemma}[theorem]{Lemma}

\newtheorem{remark}[theorem]{Remark}

\newtheorem{assumption}[theorem]{Assumption}
\newtheorem{corollary}[theorem]{Corollary}

\newtheorem{definition}[theorem]{Definition}

\usepackage{thmtools,thm-restate}
\crefname{assumption}{Assumption}{Assumptions}
\Crefname{assumption}{Assumption}{Assumptions}

\makeatletter
\AtBeginDocument{%
    \@for\revisiontheoremtype:=theorem,lemma,claim,remark,proposition,assumption,corollary,fact,observation,definition,conjecture\do{%
        \expandafter\def\csname theH\revisiontheoremtype\endcsname{\thesection.\arabic{theorem}}%
    }%
}
\makeatother

\newcommand{\Ocal}{{\cal O}}

\newcommand{\E}{\mathbb{E}}
\newcommand{\dd}{\,\mathrm{d}}
\newcommand{\Renyi}{R\'{e}nyi}
\renewcommand{\leq}{\leqslant}
\renewcommand{\geq}{\geqslant}
\newcommand{\norm}[1]{\lVert #1 \rVert}

\newcount\Comments
\newcommand{\kibitz}[2]{\ifnum\Comments=1\textcolor{#1}{#2}\fi}

\title{Learning Rate Decay Can Exponentially Accelerate SGD for Global Optimization of Nonconvex Functions}

\author{%
    Guneykan Ozgul\thanks{\texttt{guneykan.ozgul@jpmchase.com}.} \quad
    Dylan Herman \quad
    Junhyung Lyle Kim \\
    \textbf{Shree Hari Sureshbabu} \quad
    \textbf{Shouvanik Chakrabarti}\thanks{\texttt{shouvanik.chakrabarti@jpmchase.com}.} \\
    Global Technology Applied Research, JPMorganChase \\
    New York, NY 10001, USA
}

\begin{document}

\maketitle

{\renewcommand{\thefootnote}{}%
\footnotetext{\copyright~2026 JPMorgan Chase \& Co. All rights reserved.}}

\begin{abstract}
Stochastic Gradient Descent (SGD) is a workhorse algorithm for continuous optimization. It is a folk belief among deep learning practitioners that learning rate decay in SGD substantially improves performance. However, it has remained unclear whether learning rate decay offers an asymptotic advantage for \emph{global} nonconvex optimization. We answer this question in the affirmative by demonstrating a natural class of nonconvex functions for which, we prove that SGD with learning rate decay requires \emph{exponentially fewer} gradient queries for global optimization than SGD with any fixed learning rate in the small-step size regime. 
To the best of our knowledge, this is the first exponential advantage that has been demonstrated for learning rate decay. The class of functions we consider includes many popular benchmark nonconvex functions. Our results provide a robust explanatory theory for the empirically observed benefits of learning rate decay for global optimization. We show that annealing performs the non-logconcave sampling task associated to the optimization problem in polynomial time, via a novel application of weak Poincar\'{e} inequalities. 
\end{abstract}

\section{Introduction}

Minimizing nonconvex objective functions to global optimality is a central challenge in modern machine learning and scientific computing. From training deep neural networks \citep{bengio2017deep} to solving inverse problems in the physical sciences \citep{doi:10.1137/1.9780898717921}, practitioners routinely face optimization landscapes with spurious local minima, saddle points, and flat regions that impede first-order methods such as gradient descent and its variants. 
Given a continuous nonconvex objective function $f:\mathbb{R}^d \to \mathbb{R}$, we consider the problem of finding an approximate global optimum
$\widetilde{x}\in \mathbb{R}^d$ satisfying
\begin{equation*}
    f(\widetilde{x}) - f^{\star}\leq \epsilon,
\end{equation*}
where $f^{\star} =\min_{x\in \mathbb{R}^d} f(x)$. This problem is fundamentally harder than finding approximate stationary points \citep{nemirovskij1983problem, vavasis1995complexity}, for which gradient-based methods perform well. Despite the pessimistic worst-case complexity of first-order methods for nonconvex optimization \citep{Ma_2019, azizian2025the}, stochastic algorithms equipped with appropriate sources of randomness have demonstrated empirical success, motivating a rich line of theoretical inquiry into when and why noise-driven methods can efficiently navigate complex landscapes \citep{belloni2015escaping, jin2017escape, jin2021nonconvex}. The workhorse algorithm in this setting is stochastic gradient descent (SGD) \citep{robbins1951stochastic, bottou1991stochastic, amari1993backpropagation}, which updates iterates as 
\begin{equation*}
    x_{t+1} = x_t - \eta_t \nabla f(x_t) + \eta_t \xi_t,
\end{equation*}
 where $\eta_t > 0$ is the learning rate and $\xi_t$ is gradient noise arising either from minibatch approximation or intentional injection. The noise is precisely what gives SGD the ability to escape local minima and explore the landscape, making it a natural candidate for global optimization despite being a first-order method.

Unfortunately, for functions with non-trivial energy barriers between local and global minima, constant-step-size SGD requires exponential time. \citet{azizian2025the} recently established tight bounds showing that the expected hitting time of fixed-learning-rate SGD to the global minimum scales as $\exp\!\bigl(\Theta(E/\eta)\bigr)$, where $E$ is the depth of the deepest spurious local minimum along the worst-case path and $\eta$ is the step size. For learning rates $\eta$ of order $\mathcal{O}(\epsilon/d)$ chosen so that the iterates concentrate near the global minimum \citep{shi2023learning, herman2025mechanismsquantumadvantageglobal}, this yields a runtime exponential in both dimension and precision. Similar conclusions follow from the analysis of \citet{shi2023learning}. In this step-size regime, sub-exponential time is possible when all spurious minima are extremely shallow, scaling as $\mathcal{O}(\epsilon/d)$, or under other restrictive conditions such as the PL-inequality that are common proxies for convexity in optimization theory~\citep{belloni2015escaping, ozgul2025quantumspeedupsmarkovchain}. For interesting nonconvex landscapes with barriers of constant depth, constant-step-size SGD is provably intractable.

It is a folklore belief among practitioners that decaying the learning rate over the course of training substantially improves optimization: large learning rates explore broadly while small learning rates concentrate near optima. This schedule is universally deployed in deep learning, yet it has remained unclear whether learning rate decay offers a provable \emph{asymptotic} advantage for global nonconvex optimization. Crucially, the lower bound of \citet{azizian2025the} applies only to fixed-step-size SGD. An annealing schedule that varies the learning rate over time could, in principle, circumvent this barrier entirely. This raises the central question we address:
\begin{center}
    \emph{Can a simple, structure-agnostic learning rate schedule provably reduce the asymptotic runtime of SGD compared to a fixed learning rate?}
\end{center}

{
We answer this question in the affirmative. We prove that for a natural class of block-separable nonconvex functions, including multi-well potentials with locally strongly convex basins separated by constant-height barriers, and satisfying \cref{ass:sgd-spectral-transfer}, SGD with an arithmetically decaying learning rate finds an $\epsilon$-optimal solution with high probability in $\mathrm{poly}(d, 1/\epsilon)$ time, whereas any fixed learning rate $\eta=\mathcal O(\epsilon/d)$ provably requires $\exp(\Omega(d))$ time on the same instances. To our knowledge, this is the first polynomial-versus-exponential separation demonstrated for learning rate decay in nonconvex optimization for SGD.
}

\subsection{Main Result}
Our theoretical results are about the following class of functions.
\begin{definition}[Block-Separable Function]
\label{def:separable-func}
A function $f : \mathcal{X} \to \mathbb{R}$, $\mathcal{X} \subseteq \mathbb{R}^d$, is \emph{block-separable} if there exist an orthogonal matrix $R \in \mathcal{O}(d)$, a function $g: \mathcal{X} \to \mathbb{R}$, an integer $k \leq d$, and a partition of coordinates into blocks $\hat{x}_1, \ldots, \hat{x}_k$ of sizes $d_1,\ldots,d_k$ such that $f(x) = g(Rx)$, $g(x) = \sum_{i=1}^k g_i(\hat{x}_i)$, where $g_i : \mathcal{X}_i \to \mathbb{R}$ and $\sum_{i=1}^k d_i = d$, and the number of coordinates in each block is bounded above by a universal constant. If $k = d$, we say $g$ is \emph{completely separable}.
\end{definition}

The class of (block-)separable functions encompasses many standard benchmarks used in nonconvex optimization. In the widely used BBOB test suite \citep{hansen2009real}, 5 of 24 functions are explicitly separable and an additional 8 become separable after removing the orthogonal rotations that BBOB applies to test coordinate-dependence, yielding roughly 54\% of the suite that is intrinsically (block-)separable. Beyond BBOB, canonical multimodal benchmarks including Rastrigin, Schwefel, L\'{e}vy, Styblinski--Tang, and Michalewicz are all separable \citep{jamil2013literature}. 

We compare the following variants of SGD.
\begin{definition}[Fixed-step-size SGD]
\label{def:fixed-sgd-intro}
Given a fixed learning rate $\eta > 0$, \emph{fixed-step-size SGD} iterates
\begin{align*}
    x_{t+1} = x_t - \eta \nabla f(x_t) + \eta \xi_t, \qquad \xi_t \sim \mathcal{N}(0, I_d).
\end{align*}
For $\epsilon$-optimality, we consider the step size regime $\eta = \Ocal(\epsilon/d)$ of \cref{thm:lower-bound-main}. The algorithm runs until the iterate first enters an $\epsilon$-optimal region.
\end{definition}

\begin{definition}[Annealed SGD]
\label{def:annealed-sgd-intro}
\emph{Annealed SGD} uses a  learning rate schedule $\theta:[0,T]\to \mathbb{R}_{>0}$ and iterates as
\begin{align*}
    x_{t+1} = x_t - \eta(\theta_t) \nabla f(x_t) + \eta(\theta_t) \xi_t, \qquad \xi_t \sim \mathcal{N}(0, I_d),
\end{align*}
where $t \mapsto \eta(\theta_t)$ is  positive, non-increasing, and non-constant.
\end{definition}

\begin{figure}
    \centering
    \includegraphics[width=0.9\linewidth]{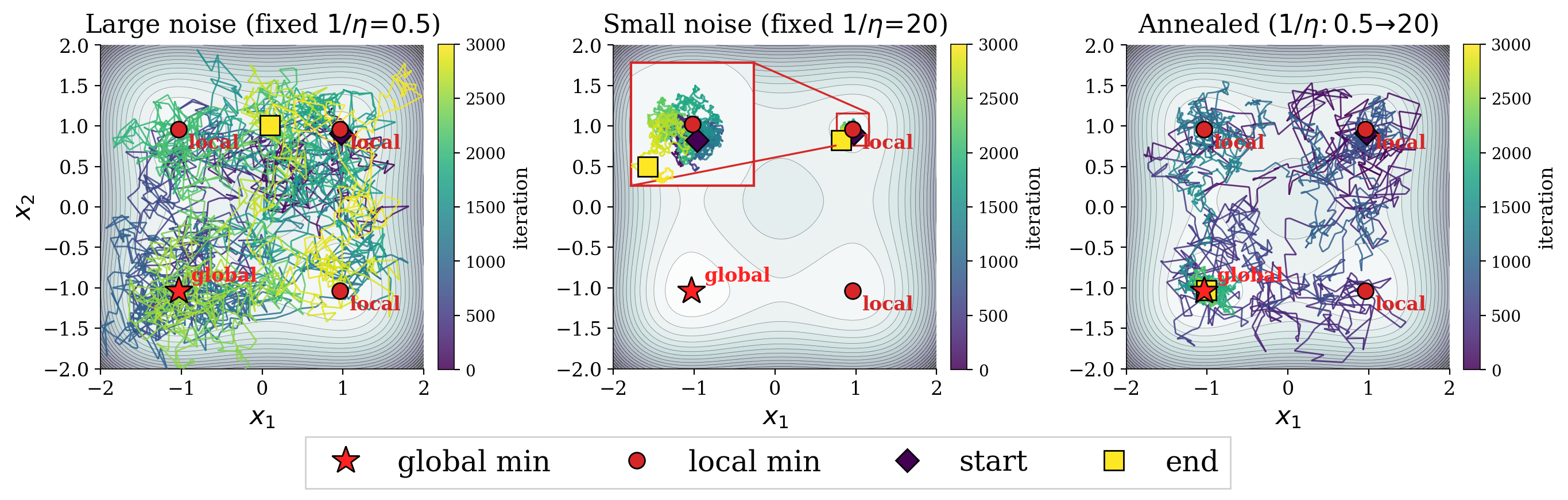}
    \includegraphics[width=0.85\linewidth]{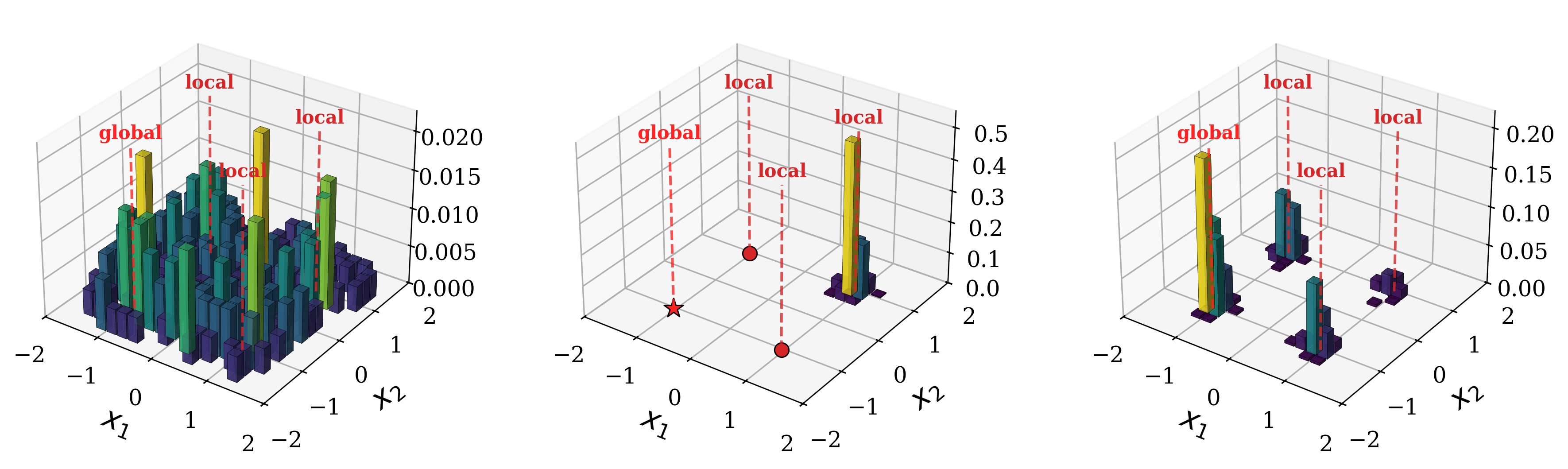}
    \caption{The trajectories of the SGD algorithm with fixed large and small step sizes (Columns 1 and 2, respectively) and a decaying learning rate on a toy problem. While a large step size allows the algorithm to explore the landscape efficiently, the small step size causes the algorithm to be trapped in a local minimum, leading to exponential runtime. The decaying step size allows the algorithm to visit the global minimum with high probability in polynomial time (Column 3).}
    \label{fig:trajectory-2d}
\end{figure}

\begin{theorem}[Informal]
\label{thm:informal-main}
Let $f$ belong to the block-separable class with bounded barriers (\cref{def:separable-func}) satisfying the conditions of \cref{thm:global_opt_sep}. Then:
\begin{enumerate}[leftmargin=*]
    \item \textbf{Lower bound for fixed step size (\cref{thm:lower-bound-main}).} For any fixed learning rate $\eta = \Ocal(\epsilon/d)$, the expected number of iterations for fixed-step-size SGD (\cref{def:fixed-sgd-intro}) to reach an $\epsilon$-optimal point is $\exp(\Omega(d))$.
    \item \textbf{Upper bound for decaying step size (\cref{thm:global_opt_sep}).} There exists an annealing schedule such that annealed SGD (\cref{def:annealed-sgd-intro}) outputs an $\epsilon$-optimal point with probability at least $2/3$ using at most $\mathrm{poly}(d, 1/\epsilon)$ gradient queries.
\end{enumerate}
In particular, under the hypotheses of \cref{thm:lower-bound-main,thm:global_opt_sep}, learning rate decay provides an exponential acceleration for global optimization of this function class when both algroithms are started from the same point from the same initial point.
\end{theorem}
An intuitive illustration of the beneficial effects of learning rate annealing is given in \cref{fig:trajectory-2d}.
Prior theoretical work on learning rate decay has focused on implicit bias~\citep{wu2018sgd, wu2020direction, andriushchenko2023sgd, even2023s} and generalization~\citep{li2019towards, you2019does, ren2024understanding}, showing that decaying schedules select flatter minima or improve test error, but has not established any asymptotic runtime improvement for the underlying optimization problem. Our result demonstrates that the benefit of learning rate decay is not merely statistical but computational: it provably transforms an exponentially hard optimization task into a polynomially easy one. The exponential lower bound for fixed step size is not an artifact of a pathological construction; it follows from a direct bound on the exit time of SGD with Gaussian noise from a spurious well (\cref{app:lower-bound}) and applies to the same concrete, well-studied functions on which practitioners routinely observe the advantages of annealing.

{
A striking feature of this result is that annealed SGD does not directly exploit the block-separable structure of the objective. The algorithm requires no knowledge of the coordinate decomposition, the barrier heights, the basin widths, or the number of local minima. A standard hyperparameter search over decay schedules suffices to achieve polynomial runtime, in contrast to methods such as coordinate descent that require knowing the decomposition. That a completely structure-agnostic modification to SGD achieves an exponential speedup suggests the benefit of learning rate decay is intrinsic to the dynamics of the algorithm rather than to any structural exploitation. Finally, while SGD is rotationally invariant, the rotation included in the function class makes them formally non-separable and challenging for algorithms like coordinate descent that can trivially optimize separable functions.
}

In general, no learning rate schedule can succeed in polynomial time for all nonconvex functions due to the NP-hardness of global optimization. Nevertheless, the magnitude of the separation we establish, combined with the structure-agnostic nature of the algorithm, suggests that qualitatively similar improvements should persist for broader function classes with weak inter-variable coupling or approximate separability. Our result provides a rigorous explanation for the empirically observed benefits of learning rate decay from a purely optimization-theoretic perspective: the early large-learning-rate phase performs implicit global exploration that no fixed learning rate can replicate, and this exploration compounds through the subsequent convergence phase.

\subsection{Theoretical Contributions}
We now discuss the main theoretical ingredients that lead to \cref{thm:informal-main}.
While the intuition for learning rate decay is natural, translating it into a rigorous guarantee is far from straightforward. SGD is a local algorithm: a large learning rate that enables broad exploration can just as easily cause the iterates to overshoot and leave a favorable basin as it can help them arrive there. Unlike a hypothetical structure-aware version of SGD that explicitly detects when a good region has been found and then switches strategy, SGD with a decaying schedule commits to reducing the learning rate on a predetermined timetable, regardless of the current state of the iterates. There is no a priori reason to expect that the iterates will be in the correct basin at the moment the learning rate becomes small enough for the iterates to stabilize. Establishing that this in fact occurs with high probability, and does so in polynomial time, requires a new analysis that goes beyond existing fixed-learning-rate convergence analyses. Our proof has three main ingredients which may each be of independent interest.
\begin{itemize}[leftmargin=6pt]
    \item \textbf{Analysis of Annealed Langevin Diffusion.} To demonstrate the effectiveness of annealing for SGD, we proceed via a correspondence with the annealed version of a continuous-time sampling algorithm based on Langevin diffusion, which forms the basis of discretized algorithms such as Langevin Monte Carlo (LMC)~\citep{Ma_2019,kinoshita2022improvedconvergenceratestochastic}. We analyze the convergence of Langevin diffusion under a temperature annealing schedule, where the learning rate plays the role of temperature. The central difficulty is that at any single low temperature, the standard Poincar\'{e} inequality yields an exponentially small spectral gap, reflecting the exponential time required for the dynamics to cross energy barriers from a worst-case initialization. Our key insight is that annealing avoids this worst case: when the temperature decreases slowly, the distribution inherited from the previous high-temperature phase serves as a favorable starting point for the next phase. To formalize this, we employ the \emph{weak Poincar\'{e} inequality}~\citep{huang2025weakpoincare, chen2025optimizationisoperimetricinequalitiessampling}, a relaxation of the standard Poincar\'{e} inequality that allows a defect term depending on the quality of the initialization. We show that for block-separable functions, the weak Poincar\'{e} constant remains polynomially bounded across all phases of the annealing schedule. Combining these per-phase bounds via an annealing theorem adapted from~\citep{huang2025weakpoincare}, we establish that the continuous annealed Langevin diffusion reaches the Gibbs measure concentrated on the global minimum in $\mathrm{poly}(d, 1/\epsilon)$ time.
    
    \item \textbf{Discretization Analysis for Annealed ULA.} We combine our diffusion analysis with a R\'{e}nyi discretization bound to prove an end-to-end polynomial-time optimization guarantee for the annealed unadjusted Langevin algorithm (ULA). The temperatures and step sizes are chosen independently, controlling discretization error and the number of gradient queries across all phases. For block-separable functions satisfying our smoothness and quadratic growth conditions, the resulting schedule finds an $\epsilon$-optimal point with probability at least $2/3$ using $\mathrm{poly}(d,1/\epsilon)$ gradient queries (\cref{thm:ula-optimization-formal}).

    \item \textbf{Spectral Transfer to Annealed SGD.} The next problem is to translate our result on annealing the temperature of Langevin diffusion into one about annealing the learning rate of SGD. Standard pathwise discretization estimates can be too restrictive when the early learning rates are large. We derive a discrete R\'{e}nyi contraction bound, combine it across phases, and prove vanishing stationary bias at the final learning rate under spectral transfer assumption~(\cref{ass:sgd-spectral-transfer}). Each phase supplies the starting distribution for the next; this distribution is a conclusion of the induction over phases using Poincar\'{e} inequality.
\end{itemize}

We finally comment on the lower bound in \cref{thm:informal-main}. This result shows concretely that existing fixed-rate analyses of SGD still yield exponential runtimes  for block-separable functions as for general nonconvex functions. The main idea is the following: we first explain why a step size $\eta = \mathcal{O}(\epsilon/d)$ is the regime we consider for global optimization. We then bound the exit time from a spurious well directly to prove that there exist completely separable functions for which SGD requires $\exp(\Omega(d))$ iterations at any sufficiently small step size. Details of this calculation are given in \cref{app:lower-bound}.

\textbf{Organization. }
The rest of the paper is organized as follows. Our complete proofs are given in the appendix. \cref{sec:main-annealing} gives more detail and intuition on the mixing analysis of annealed Langevin diffusions. \cref{sec:discretization-main} gives more detail on the translation of this analysis to results on the convergence of SGD with decaying learning rate, and how this results in a proof of the upper bound in \cref{thm:informal-main}. \cref{sec:experiments} describes numerical simulations that confirm our theoretical results and indicate the applicability of the identified phenomena beyond the scope of the analysis.

\subsection{Related Work}

SGD has been extensively analyzed in the convex setting \citep{nemirovski2009robust, moulines2011non,  rakhlin2011making,  schmidt2013fast, needell2014stochastic, gower2019sgd}.
For general nonconvex functions, SGD has been typically analyzed for convergence to a stationary point \citep{ghadimi2013stochastic, khaled2020better,  arjevani2023lower, patel2021stochastic,   sebbouh2021almost,
liu2022almost}.
When the objective has favorable structure, such as the Polyak--{\L}ojasiewicz condition~\citep{karimi2016linear, khaled2020better}, quasar convexity~\citep{gower2021sgd}, or a star-convex path~\citep{zhou2019sgd}, fixed-step-size SGD can converge to the global minimum efficiently. SGD has also been analyzed for stylized nonconvex problems like matrix completion or phase retrieval~\citep{de2015global, ma2023global}, where stationary points are globally optimal \citep{ge2016matrix, ge2015escaping}.

The Langevin Monte Carlo (LMC) algorithm, a discretization of overdamped Langevin diffusion, has well-established convergence guarantees in the log-concave setting~\citep{dalalyan2016theoreticalguaranteesapproximatesampling, JMLR:v20:18-173, vempala2019rapid}. For nonconvex objectives, \citet{raginsky2017non} and \citet{zhang2018hittingtimeanalysisstochastic} gave non-asymptotic bounds that depend exponentially on barrier height, and \citet{azizian2025the} provided tight matching bounds confirming this dependence.

Simulated annealing is a classical strategy for escaping local minima~\citep{belloni2015escaping}. Recently, \citet{huang2025weakpoincare} introduced weak Poincar\'{e} inequalities for analyzing simulated annealing on discrete spin-glass models, and \citet{chen2025optimizationisoperimetricinequalitiessampling} developed connections between isoperimetric inequalities and optimization. Our work adapts these ideas to continuous SGD with learning rate decay and adds a discretization analysis bridging the Langevin diffusion and the discrete iterates. The annealing of noise scales also appears in score-based generative models~\citep{song2019generative, song2021scorebasedgenerativemodelingstochastic, chen2023improved, benton2024nearly}, which share technical ingredients with our analysis despite targeting a different goal.

Learning rate decay of SGD itself has been studied through the lens of implicit bias~\citep{wu2018sgd, wu2020direction, andriushchenko2023sgd, even2023s} and generalization~\citep{li2019towards, you2019does, ren2024understanding}. Large learning rates have been connected to a catapult mechanism~\citep{lewkowycz2020large}, while \citet{sclocchi2024different} studied how step-size-to-batch-size ratios determine different SGD regimes. \citet{mohtashami2023large} also considered a large learning rate for nonconvex optimization, but is restricted to gradient descent with a large but constant step size.

To summarize, to the best of our knowledge, this is the first work that studies the provable performance of SGD with \textit{decaying step size} for global optimization, establishing a provable separation in global optimization runtime between annealed and non-annealed SGD.

\section{Analysis of Annealed Langevin Diffusion}
\label{sec:main-annealing}
We now present our proof strategy and main results. Our analysis proceeds by relating annealed SGD to a learning-rate-dependent Langevin SDE,
\begin{align}
    \label{eq:learning-rate-sde}
    \mathrm{d}X_t = -\nabla f(X_t)\mathrm{d}t+\sqrt{\eta(\theta(t))}\mathrm{d}B_t
\end{align}
where $B_t$ is a standard Brownian motion and $\theta:[0,T]\to \mathbb{R}_{>0}$ is the learning rate schedule. Let $\mu_{t}$ denote the distribution of $X_t$. For a constant learning rate $\eta$, as $t \rightarrow \infty$, we have
$
    \mu_{t} \rightarrow \pi_{f,\beta}(x) \propto e^{-\beta f(x)},
$ 
known as the Gibbs distribution, where $\beta = 2/\eta$. The Euler--Maruyama discretization of \eqref{eq:learning-rate-sde} yields the SGD update with Gaussian noise:
\begin{equation*}
    x_{t+1} = x_t - \eta(\theta_t) \nabla f(x_t)+ \eta(\theta_t)\xi_t
\end{equation*}
where $\xi_t \sim \mathcal{N}(0, I_d)$ are independent standard Gaussian random vectors. While the general form of SGD does not have Gaussian normal noise in practice, we make this simplification for the theoretical analysis as done in~\citep{shi2023learning} and expect the results to hold as long as the noise distributions have decaying tails. We then implement the annealing through a decreasing sequence of learning rates $\eta_1 > \eta_2 > \cdots > \eta_M > 0$. In phase $k = 1, \ldots, M$, the algorithm runs with constant learning rate $\eta_k$ for a duration $\tau_k$ chosen so that the constant-temperature Langevin SDE, $dX_t = -\nabla f(X_t)\dd t + \sqrt{\eta_k}\,dB_t$,  approximately reaches its stationary distribution. In this sense, our analysis uses the sampling problem as a proxy for the optimization, since when the step size $\eta_k$ is small, the stationary distribution $\pi_{\eta_k}$ has large probability around the global minimum. These techniques are used in the context of the Langevin Monte Carlo (LMC) algorithm for nonconvex optimization~\citep{Ma_2019,kinoshita2022improvedconvergenceratestochastic}, although the runtime is usually given by an exponentially small geometric quantity such as the log Sobolev or Cheeger constant characterizing the conductance of this SDE.

We now define the key functional inequality underlying our analysis. A reversible Markov semigroup $P$ with stationary distribution $\pi$ and Dirichlet form $\mathcal{E}_{\pi}$ satisfies a $\delta_{\textup{p}}$-Poincar\'{e} inequality if $\mathcal{E}_{\pi}(h,h) \geq \delta_{\textup{p}}\, \mathrm{Var}_\pi[h]$ for all functions $h$. For the Langevin diffusion \eqref{eq:learning-rate-sde} with a constant learning rate $\eta$, the Dirichlet form is $\mathcal{E}_{\pi}(h,h) = \frac{\eta}{2}\int \|\nabla h\|^2 d\pi$. For multimodal distributions, the Poincar\'{e} constant $\delta_p$ is typically exponentially small in the barrier height and inverse learning rate and characterizes the optimization hardness from the worst initialization. The weak Poincar\'{e} inequality relaxes the standard Poincar\'{e} inequality by allowing a defect term and characterizes the hardness from most initializations~\citep{chen2025optimizationisoperimetricinequalitiessampling}. While we provide the self-contained definitions in our proofs, for a more detailed background on Markov semigroups and functional inequalities, we refer the reader to~\citep{Bakry2013AnalysisAG}.

\begin{definition}[Weak Poincar\'{e} inequality]
\label{def:weak-poincare}
We say $P$ with state space $\Omega$ satisfies a $(\delta_{\textup{wp}}(s), s)$-weak Poincar\'{e} inequality with respect to its stationary distribution $\pi$ if for all $h : \Omega \to \mathbb{R}$,
{
\begin{equation*}
\mathcal{E}_{\pi}(h,h)\geq \delta_{\textup{wp}}(s)\, \operatorname{Var}_\pi[h] - s\, \|h\|^2_{\infty}.
\end{equation*}
}
\end{definition}
When $s=0$, we obtain the standard Poincar\'{e} inequality. The defect term $s\|h\|_\infty^2$ allows the constant $\delta_{\textup{wp}}(s)$ to be much larger than the standard Poincar\'{e} constant, which is critical for obtaining polynomial runtime bounds. Crucially, under a $(\delta_{\textup{wp}}(s), s)$-weak Poincar\'{e} inequality, the convergence of Langevin dynamics starting from an initial distribution $\nu_0$ satisfies
{
\begin{equation*}
    \chi^2(\nu_t \| \pi) \leq e^{-\delta_{\textup{wp}}(s) t}\, \chi^2(\nu_0 \| \pi) + \frac{s}{\delta_{\textup{wp}}(s)}\left\|\frac{d\nu_0}{d\pi} - 1\right\|_\infty^2,
\end{equation*}
}
where 
$\chi^2(\nu || \pi) = \int \left(\frac{d\nu}{d\pi}\right)^2 d\pi - 1$ is the chi-square divergence. The key consequence is that if the initial distribution $\nu_0$ is a warm start with bounded density ratio (i.e., $\|d\nu_0/d\pi - 1\|_\infty = \Ocal(1)$), then the convergence rate is governed by $\delta_{\textup{wp}}$ rather than by the standard Poincar\'{e} constant. In the annealing setting, the output of each high learning rate phase provides precisely such a warm start for the low learning rate phase, enabling the algorithm to bypass the exponential mixing bottleneck that slows down the fixed learning rate dynamics.

\textbf{Proof Sketch.} The proof has three main steps. First, we show that the continuous annealed dynamics~\eqref{eq:learning-rate-sde} converge to a distribution $\mu_T$ that assigns probability $\Omega(1)$ to a small neighborhood of $x^{\star}$, where $x^{\star}$ is the global minimum, with total runtime governed by the minimum weak Poincar\'{e} constant across phases. Next, we derive a lower bound on the weak Poincar\'{e} constant via perturbation analysis, establishing that this constant remains large even at small learning rates. Third, we show that the discretized algorithm (annealed SGD) converges to the same distribution with only polynomial overhead, using a new discretization argument that decouples the early exploratory phase from the late convergent phase.

\subsection{Bounding the Weak Poincar\'{e} Constant}
\begin{figure}
    \centering
    \includegraphics[width=\linewidth]{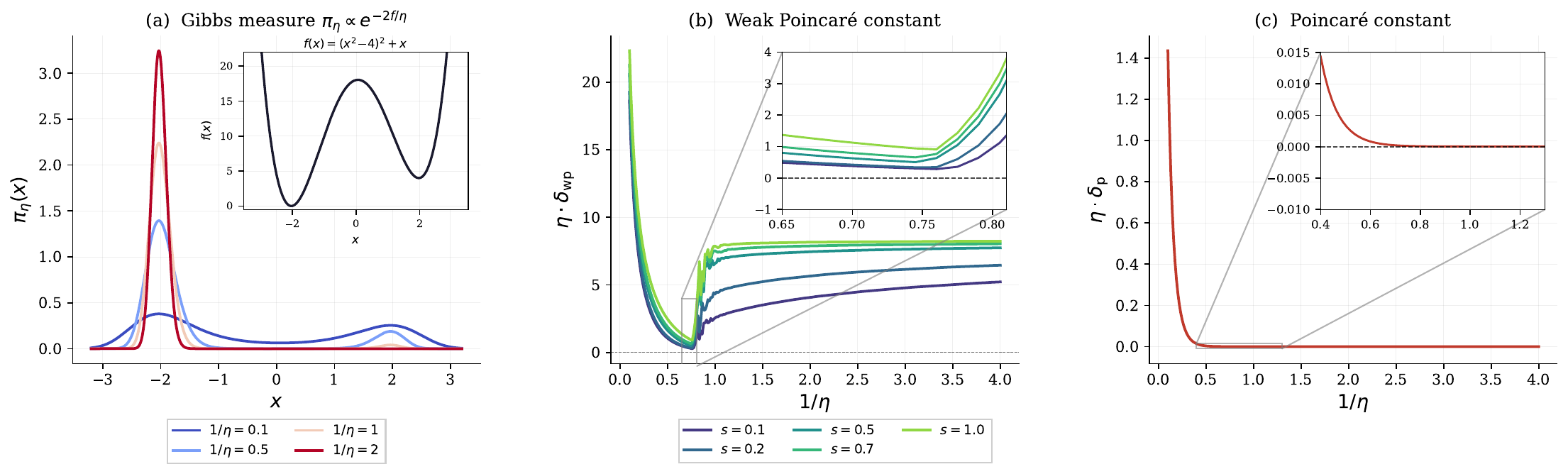}
    \caption[Comparison of weak Poincar\'e to Poincar\'e constant]{Comparison of the weak Poincar\'e constant (normalized by $1/\eta$) to the standard Poincar\'e constant (normalized by $1/\eta$) for a 1-D toy function\footnotemark. While the standard Poincar\'e constant decays exponentially, the weak Poincar\'e constant decays initially but then starts increasing, indicating that the effect of local minima can be benign with a proper annealing scheme.}
    \label{fig:poincarecomparison}
\end{figure}
\footnotetext{The plot for the standard Poincar\'e constant is obtained by diagonalization of the actual generator of the dynamics. The plot for the weak Poincar\'e constant is obtained via \cref{thm:weak-poincare-perturbation} by choosing various Gaussian distributions as base measure to get useful lower bound estimates causing some discontinuities in the plot. We expect the actual weak Poincar\'e constant to be smooth and even larger than our estimates. Since these are extensive quantities in $1/\eta$ for small $\eta$, we normalize the constants.}

While the weak Poincar\'{e} inequality relaxes the standard Poincar\'{e} inequality significantly, its constant is still a global geometric quantity that depends heavily on the function geometry. To quantify this constant, we derive a novel comparison method in~\cref{thm:weak-poincare-perturbation}. Informally, given a measure $m$ that satisfies a Poincar\'{e} inequality, \cref{thm:weak-poincare-perturbation} allows us to quantify the weak Poincar\'{e} constant of a measure $\pi$. This theorem generalizes the standard Holley--Stroock type of perturbation results used previously in~\citep{Ma_2019} to the weak Poincar\'{e} setting. However, as opposed to the Holley--Stroock perturbation theorem, our result is also applicable to distributions with an unbounded density ratio. As an application of this theorem, we compare the weak Poincar\'{e} constant of the SGD dynamics to the standard Poincar\'{e} constant.
As opposed to the standard Poincar\'{e} constant, \cref{fig:poincarecomparison} shows that for a 1-D toy problem, the weak Poincar\'{e} constant initially decreases and then starts increasing again around $\eta \sim 1$, where the stationary distribution concentrates around the basin of the global minimum. This indicates that the dynamics stabilize early on, at a learning rate much larger than what standard analyses require. As a further corollary of this theorem, by choosing $m$ to be a Gaussian around $x^{\star}$ and $\pi = \pi_{f,\beta}$, we show the following theorem:

{
\begin{theorem}[Informal]
\label{cor:weak-poincare-linit}
Let $f:\mathbb{R}^d \to \mathbb{R}$ be a $C^3$ function with a unique global minimizer $x^{\star}$, positive definite Hessian $H \coloneqq \nabla^2 f(x^{\star})$ with smallest eigenvalue $\omega$, and quadratic growth away from $x^{\star}$. Then for all $\beta$ sufficiently large, the Gibbs measure $\pi_\beta(x) \propto \exp(-\beta f(x))$ satisfies the weak Poincar\'e inequality
\begin{equation*}
\textup{Var}_{\pi_\beta}[h] \leq \frac{1 + \mathcal{O}(\beta^{-1/4})}{\omega\beta} \int \|\nabla h\|^2 \, d\pi_\beta + C\exp\!\big({-c\,\beta^{1/6}}\big) \, \|h\|_\infty^2,
\end{equation*}
where $c, C > 0$ are constants independent of $\beta$. The formal statement is \cref{thm:weakpoincare-formal}.
\end{theorem}
}
The importance of \cref{cor:weak-poincare-linit} is that for $\beta = 2/\eta$, the weak Poincar\'{e} constant does not exponentially decay, but converges to $\beta\omega = \frac{2\omega}{\eta}$ asymptotically (equivalently, to $\omega$ when the inequality is stated with the Dirichlet form $\frac{\eta}{2}\int \|\nabla h\|^2 \, d\pi_\beta$ of the generator of Langevin dynamics). The asymptotic result alone does not imply polynomial runtimes for general nonconvex functions since it is only valid in the small $\eta$ limit and $\delta_{\textup{wp}}$ can be very small for larger learning rates. However, the result is also important from a mathematical point of view, since for functions with unique global minima, similar asymptotic behavior is known for the spectral gap of Schr\"odinger operators~\citep{herman2025mechanismsquantumadvantageglobal} that are connected to Langevin diffusion. While the spectrum of the generator of Langevin dynamics contains an eigenstate that could make the spectral gap very small, this eigenstate is not significant in the limit where $\beta\to \infty$.

 The fact that the limit of the weak Poincar\'{e} constant is $\frac{2\omega}{\eta}$ is not coincidental. This matches the Poincar\'{e} constant $\beta\omega$ of the Gaussian $\mathcal{N}(x^{\star}, (\beta H)^{-1})$, the Gibbs measure of a strongly convex quadratic function with its minimum at $x^{\star}$. From an algorithmic point of view, this result implies that existence of a learning rate schedule that can escape from local minima early on will localize around the global minima and effectively the algorithm will only be affected by the strongly convex well around the global minima.

For the setting of block separable functions with constant size blocks, we prove that each block satisfies a weak Poincar\'e inequality with
{
\begin{equation}
(\delta_{\textup{wp}}(s),s)=\left(\frac\xi4,\frac\xi4Ce^{-2c/\eta}\right),
\qquad \eta<\eta_*,
\label{eq:block-wpi-main}
\end{equation}
}
where $\eta_* = \Omega_d(1)$ and the constants $\xi,c,C>0$ do not depend on $d$ (\cref{cor:wpi-exp-defect}). The exponentially small defect in \eqref{eq:block-wpi-main} provides the stronger inequality used at low temperature. The uniform weak Poincar\'e family of \cref{lem:ula-block-estimates} handles the earlier temperatures.

\section{Discretization and Learning Rate Schedule}
\label{sec:discretization-main}
When the goal is to simulate Langevin dynamics, it is common in the literature to select a time-discretization step size that is different from the corresponding learning rate \citep{zou2021faster}. However, this algorithm is slightly different from the form of SGD in (\cref{def:annealed-sgd-intro}).

{
We first state the optimization result for this intermediate algorithm. The unadjusted Langevin algorithm (ULA) uses the update
\begin{equation}
\label{eq:ula-main}
    x_{t+1} = x_t - h_t \nabla f(x_t) + \sqrt{h_t \eta_t}\, \xi_t, \qquad \xi_t \sim \mathcal{N}(0, I_d).
\end{equation}
Here $\eta_t$ is the noise level and $h_t$ is the discretization step. The choice $h_t = \eta_t$ gives annealed SGD (\cref{def:annealed-sgd-intro}). We begin with the optimization guarantee obtained when the steps are free.
}

{
\begin{theorem}[Informal]
\label{thm:ula-optimization-main}
Under the hypotheses of \cref{thm:ula-optimization-formal}, the annealed ULA algorithm \eqref{eq:ula-main} uses independently chosen noise levels $\eta_r$ and discretization steps $h_r$. For every $\epsilon\in(0,1]$, it outputs an $\epsilon$-optimal point with probability at least $2/3$ using $\mathrm{poly}(d/\epsilon)$ gradient queries.
\end{theorem}
}

Unlike ULA, the discrete iterates of SGD can deviate from the continuous Langevin diffusion despite its convergence. Similar phenomena have been observed for the gradient descent algorithm, that is, despite the success of gradient descent, it can follow a significantly different path from continuous gradient flow~\citep{de2015global}.

To analyze SGD, we instead transfer the weak Poincar\'e inequalities of Langevin diffusion to the discrete SGD kernel. This gives an annealing proof analogous to the continuous-time argument, with phase lengths determined by the discrete contraction bound.
Additionally, the fact that the learning rate is changing over time also introduces challenges into the discretization analysis. Particularly, convergence results on continuous-time Langevin dynamics are all in terms of quantities that fail to satisfy the triangle inequality. Thus we have to carefully divide the learning rate sequence in such a way that we can compose all sources of error properly.

The next theorem states the resulting SGD optimization guarantee.

{
\begin{theorem}[Informal]
\label{thm:sgd-optimization-main}
Under the hypotheses of \cref{thm:global_opt_sep} and \cref{ass:sgd-spectral-transfer}, the annealed SGD algorithm uses $h_r=\eta_r$ in every phase. For every $\epsilon\in(0,1]$, it outputs an $\epsilon$-optimal point with probability at least $2/3$ using $\mathrm{poly}(d/\epsilon)$ gradient queries.
\end{theorem}
}
The learning rate schedule is an arithmetically decaying sequence from $\eta_1$ to $\eta_M$ with sufficiently small constant decrements. 
We select a constant $\eta_1$. To obtain an $\epsilon$ approximate global minimizer with constant probability, it suffices to take the final learning rate to be $\eta_M=\Theta(\epsilon^2/d)$ (proven in \cref{sec:optimization_block_sep}).

The spectral-transfer condition in \cref{ass:sgd-spectral-transfer} compares one-step contraction of each fixed-rate SGD kernel with Langevin diffusion, up to a constant or polylogarithmic factor $\kappa_d$. Constant comparisons between invariant and Gibbs laws control initialization and transfers between phases. These are mild requirements on each constant-size block, one-dimensional for completely separable functions, with bounds independent of the ambient dimension and target accuracy.  

At low temperature, the the SGD update yield an order-$\eta$ stationary bias (\cref{lem:sgd-stationary-bias}), so the sum of the block biases at $\eta_M=\Theta(\epsilon^2/d)$ is small enough for optimization. On the other hand, for general nonconvex functions, an initial learning rate $\eta_1=\mathcal O(1/d)$ can lead to exponential runtime whereas \cref{thm:sgd-optimization-main} allows a constant initial learning rate and proves convergence using the inherited weak Poincar\'e inequalities when the learning rate is decayed slowly instead of initialized at such small values.

\raggedbottom
\section{Numerical Results}
\label{sec:experiments}

We now empirically examine the
polynomial-versus-exponential separation predicted by our theory through controlled experiments on
the following three functions: 
\begin{enumerate}[leftmargin=7mm,topsep=4pt,itemsep=2pt,parsep=0pt,partopsep=0pt]
    \item $f(x) = \sum_{i=1}^d [(x_i^2 - 1)^2 + 0.3 x_i]$ \quad (Separable Double-Well)
    \item $f(x) = Ad + \sum_{i=1}^d [ x_i^2 - A \cos (2 \pi x_i) ]$ \quad (Rastrigin with $A=1$)
    \item $f(x) = \sum_{i=1}^d [(x_i^2 - 1)^2 + 0.3 x_i] + 0.5 \sum_{i=1}^{d-1} (x_i - x_{i+1})^2  $ \quad (Coupled Double-Well)
\end{enumerate}
The first two functions belong to the function class satisfying \cref{def:separable-func}. We test the third function, which is non-separable, to assess how annealed SGD performs in the setting beyond \cref{thm:informal-main}.

\begin{figure}[h]
    \centering
    \includegraphics[width=1.0\linewidth]{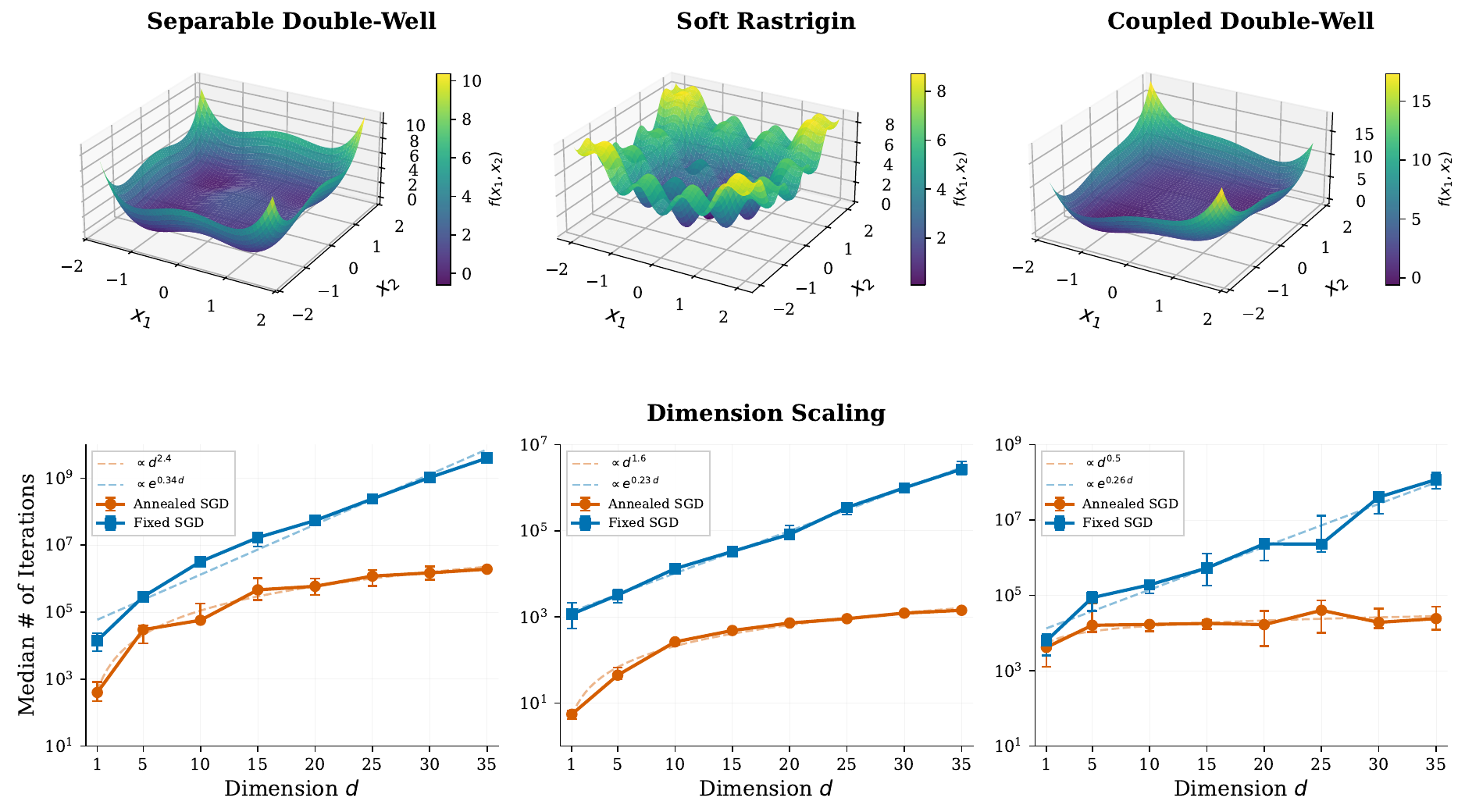}
    \caption{Median number of iterations to convergence as a function of dimension $d$ for fixed-rate SGD (blue) and annealed SGD (orange), both with $h_k=\eta_k$. Error bars show the interquartile range over 10 runs. Dashed curves are descriptive fits; the coupled annealed power-law fit is weak.
    }
    \label{fig:dimension}
\end{figure}

\textbf{Setup.}
We implement fixed-step-size and annealed SGD (\cref{def:fixed-sgd-intro,def:annealed-sgd-intro}) with Gaussian gradient noise, using the same learning rate $\eta_t$ to scale the gradient and noise terms.
Both algorithms are initialized near the spurious local minimum for all three functions.
Fixed-step-size SGD keeps its learning rate constant throughout each run.
Annealed SGD uses function-specific arithmetic learning-rate decay across constant-rate phases.
Convergence is declared when $\|x - x^\star\|_\infty$ falls below a tolerance.
Specific parameters for each function are detailed in \cref{app:numeric-details}.

\textbf{Results.}
Figure~\ref{fig:dimension} shows a clear separation between the fixed-step-size SGD algorithm and the annealed SGD schedules over the tested dimensions. Fixed-step-size SGD becomes rapidly more expensive as dimension increases, with approximately exponential growth across all three functions.
In contrast, Annealed SGD exhibits polynomial growth in median iteration counts, $\tau\propto d^{c'}$, for the separable double-well and Rastrigin functions, while the coupled double-well median grows only mildly with dimension. This polynomial-versus-exponential pattern is consistent across the three functions tested. Indeed, Annealed SGD performs well in higher dimensions using the same function-specific schedules, as shown in \cref{tab:high-dim-sgd-arithmetic}. In particular, even the largest median recorded run time (separable double well, $d=400$) is about 42 seconds, which is less than one minute.

\noindent\begin{minipage}{\linewidth}
\begin{table}[H]
\centering
\small
\setlength{\tabcolsep}{4pt}

\caption{Annealed SGD with arithmetic learning-rate decay and $h_t=\eta_t$ in high dimensions. Median recorded update counts ($\times 10^5$) and median run times (seconds), $\pm$ one empirical standard deviation over 10 trials per function and dimension. The stopping criterion is $\|x-x^\star\|_\infty<\delta$, with $\delta=0.5$ for both double-well functions and $\delta=0.3$ for Rastrigin.}
\label{tab:high-dim-sgd-arithmetic}
\resizebox{\linewidth}{!}{%
\begin{tabular}{@{}ll cccc@{}}
\toprule
& & $d=100$ & $d=200$ & $d=300$ & $d=400$ \\
\midrule
\multirow{3}{*}{\shortstack[l]{Iterations\\($\times 10^5$)}}
  & Separable Double-Well & $32.31\pm67.26$ & $36.95\pm40.53$ & $76.67\pm20.20$ & $79.68\pm42.17$ \\
  & Rastrigin & $0.0375\pm0.0036$ & $0.1844\pm0.2210$ & $3.0936\pm0.7004$ & $4.1200\pm0.0013$ \\
  & Coupled Double-Well & $0.1875\pm0.2942$ & $0.2822\pm0.2131$ & $0.2265\pm0.1189$ & $0.5260\pm0.7577$ \\
\midrule
\multirow{3}{*}{\shortstack[l]{Run time\\(seconds)}}
  & Separable Double-Well & $4.20\pm8.75$ & $9.42\pm10.34$ & $29.60\pm8.03$ & $41.57\pm21.87$ \\
  & Rastrigin & $0.0076\pm0.0006$ & $0.0671\pm0.0851$ & $1.67\pm0.37$ & $2.96\pm0.05$ \\
  & Coupled Double-Well & $0.0251\pm0.0401$ & $0.0742\pm0.0575$ & $0.0905\pm0.0468$ & $0.284\pm0.415$ \\
\bottomrule
\end{tabular}%
}

\end{table}

\section{Conclusion}
We have shown that for a special but broad class of nonconvex functions satisfying \cref{ass:sgd-spectral-transfer}, learning rate decay converts the exponential runtime of fixed-step-size SGD in the regime $\eta=\mathcal O(\epsilon/d)$ into a polynomial one, establishing the first polynomial-versus-exponential separation for this ubiquitous heuristic. The primary limitation is the current restriction to block-separable functions, and the primary open question is whether our techniques can be used to prove (potentially smaller) separations for more general functions, such as those with weak inter-variable couplings.

\end{minipage}
\clearpage
\flushbottom
\section*{Disclaimer}
This paper was prepared for informational purposes by the Global Technology Applied Research
center of JPMorgan Chase \& Co. This paper is not a product of the Research Department of
JPMorgan Chase \& Co. or its affiliates. Neither JPMorgan Chase \& Co. nor any of its affiliates
makes any explicit or implied representation or warranty and none of them accept any liability
in connection with this paper, including, without limitation, with respect to the completeness,
accuracy, or reliability of the information contained herein and the potential legal, compliance, tax,
or accounting effects thereof. This document is not intended as investment research or investment
advice, or as a recommendation, offer, or solicitation for the purchase or sale of any security,
financial instrument, financial product or service, or to be used in any way for evaluating the merits
of participating in any transaction.

\section*{AI Usage}
The initial version of this paper was written without the use of AI tools. The authors subsequently used GPT 6.1 Sol and Claude 5.5 Opus to correct the errors identified in the proofs during the review process, expand the details of the proofs and polish the writing of the appendix. The authors verified the details of the these proofs and take responsibility for the correctness and presentation of the paper.

\bibliographystyle{plainnat}
\bibliography{biblio}

\newpage
\appendix
\crefalias{section}{appendix}

\noindent\makebox[\linewidth]{\rule{\linewidth}{1.5pt}}
\section*{\centering \large Technical Appendices and Supplementary Material for\\ ``Learning Rate Decay Can Exponentially Accelerate SGD for Global Optimization of Nonconvex Functions''}
\noindent\makebox[\linewidth]{\rule{\linewidth}{1.5pt}}
\tableofcontents

\section{Weak Poincar\'e Inequalities}
\label{sec:weak-poincare-proofs}

In this section, we transfer a Poincar\'e inequality for a reference measure to a weak Poincar\'e inequality for a measure whose density ratio with respect to the reference measure is bounded above and below in a ball, and we apply this comparison to the Gibbs measures used in the ULA and SGD analyses.

For Langevin diffusion with potential $f$ at a constant noise level $\eta>0$, the stationary law is $\pi \propto \exp(-2f/\eta)$ and its Dirichlet form is $\mathcal{E}_{\pi}(h,h) = \frac{\eta}{2}\int \|\nabla h\|^2 d\pi$.

We say that $\pi$ satisfies a $(\delta, \epsilon)$-WPI, where $\delta > 0$ and $\epsilon \geq 0$, if every bounded smooth function $h$ satisfies
{
\begin{align}
\label{eq:wpi-convention}
    \mathcal{E}_{\pi}(h,h) \geq \delta\, \textup{Var}_{\pi}(h) - \epsilon \|h\|_{\infty}^2.
\end{align}
}
By the formula for $\mathcal{E}_{\pi}$, if every bounded smooth function $h$ satisfies
{
\begin{align}
\label{eq:wpi-variance-form}
    \textup{Var}_{\pi}(h) \leq C_P \int \|\nabla h\|^2 d\pi + s\|h\|_{\infty}^2
\end{align}
}
for constants $C_P > 0$ and $s \geq 0$, then $\pi$ satisfies a $(\delta, \epsilon)$-WPI with $\delta = \eta/(2C_P)$ and $\epsilon = \delta s$.

\begin{restatable}[Weak Poincar\'e Inequality via Density Ratio Perturbation]{theorem}{weakPoincareComparison}
\label{thm:weak-poincare-perturbation}
Let $m$ be a probability measure satisfying the Poincar\'e inequality
\[
\gamma \, \textup{Var}_{m}[h] \leq \int \|\nabla h\|^2 \, dm
\]
for all smooth $h$, where $\gamma > 0$ is the spectral gap. Let $\pi$ be a probability measure absolutely continuous with respect to $m$, and let $\rho = d\pi/dm$ denote the Radon--Nikodym derivative. For a given $\Delta > 0$ and a point $x^{\star} \in \mathbb{R}^d$, define the local density ratio bounds
\[
K_{\textup{max}} \coloneqq \sup_{\|x - x^{\star}\| \leq \Delta} \rho(x), \qquad K_{\textup{min}} \coloneqq \inf_{\|x - x^{\star}\| \leq \Delta} \rho(x),
\]
and let $R_\Delta = \mathcal{B}(x^{\star}, \Delta)$. Then for every $\epsilon > 0$, $\pi$ satisfies the weak Poincar\'e inequality
{
\begin{align*}
\textup{Var}_{\pi}[h] &\leq \frac{(1+\epsilon)^2 K_{\textup{max}}}{\gamma K_{\textup{min}}} \int \|\nabla h\|^2 \, d\pi \\
&\quad + \left( 4\Pr_{\pi}(x \notin R_\Delta) + \Big(1 + \frac{1}{\epsilon}\Big)\Big(1 + \frac{16(1+\epsilon)}{\gamma \Delta^2}\Big) K_{\textup{max}} \Pr_{m}(x \notin R_{\Delta/2}) \right) \|h\|_{\infty}^2.
\end{align*}
}
\end{restatable}
\begin{proof}
{
For every bounded real-valued $h$ and constant $c$, expanding the square gives
\begin{align*}
\mathbb{E}_{\pi}[(h-c)^2]
&=\textup{Var}_{\pi}[h]+(\mathbb{E}_{\pi}h-c)^2.\\
\intertext{Taking $c=\mathbb{E}_{m}h$ yields}
\mathbb{E}_{\pi}[(h-\mathbb{E}_{m}h)^2]
&=\textup{Var}_{\pi}[h]+(\mathbb{E}_{\pi}h-\mathbb{E}_{m}h)^2
\geq\textup{Var}_{\pi}[h].
\end{align*}
We split this upper bound over the ball $R_\Delta$ and its complement:
}
\begin{align*}
\textup{Var}_{\pi}[h] &\leq \mathbb{E}_{\pi}\|h - \mathbb{E}_{m} h\|^2 \\
&= \int_{\|x - x^{\star}\| \leq \Delta} \|h(x) - \mathbb{E}_{m} h\|^2 \, d\pi(x) + \int_{\|x - x^{\star}\| > \Delta} \|h(x) - \mathbb{E}_{m} h\|^2 \, d\pi(x) \\
&\leq K_{\textup{max}} \, \mathbb{E}_{m}\|h - \mathbb{E}_{m} h\|^2 + 4\Pr_{\pi}(\|x - x^{\star}\| > \Delta) \, \|h\|_{\infty}^2.
\end{align*}

It remains to bound $\textup{Var}_m[h]$ in terms of the Dirichlet form under~$\pi$. Let $\phi: [0,\infty) \to [0,1]$ be a smooth non-increasing function satisfying $\phi(t) = 1$ for $t \leq 1$, $\phi(t) = 0$ for $t \geq 2$, and $|\phi'(t)| \leq 2$ for all~$t$. Define the cutoff $\chi(x) = \phi(2\|x - x^{\star}\| / \Delta)$, so that $\chi = 1$ on $R_{\Delta/2}$, $\textup{supp}(\chi) \subseteq R_{\Delta}$, and $\|\nabla \chi(x)\| \leq 4/\Delta$.
Fix $\epsilon > 0$, and recall that $\|u + v\|^2 \leq (1+\epsilon)\|u\|^2 + (1+1/\epsilon)\|v\|^2$ for all vectors $u, v$.
The same expansion of the square, applied under $m$ with $c=\mathbb E_m(h\chi)$, justifies the first inequality next.

{
\begin{align*}
\mathbb{E}_{m}\|h - \mathbb{E}_{m} h\|^2 &\leq \mathbb{E}_{m}\|h - \mathbb{E}_{m}(h\chi)\|^2 \\
&\leq \Big(1 + \frac{1}{\epsilon}\Big) \mathbb{E}_{m}\|h(1 - \chi)\|^2 + (1+\epsilon)\, \mathbb{E}_{m}\|h\chi - \mathbb{E}_{m}(h\chi)\|^2 \\
&\leq \Big(1 + \frac{1}{\epsilon}\Big) \Pr_{m}\!\big(\|x - x^{\star}\| > \Delta/2\big) \|h\|_\infty^2 + (1+\epsilon)\, \mathbb{E}_{m}\|h\chi - \mathbb{E}_{m}(h\chi)\|^2.
\end{align*}
}

For the second term, we apply the Poincar\'e inequality for $m$ to $h\chi$ and expand:
{
\begin{align*}
\mathbb{E}_{m}\|h\chi - \mathbb{E}_{m}(h\chi)\|^2 &\leq \frac{1}{\gamma} \int \|\nabla(h\chi)\|^2 \, dm = \frac{1}{\gamma} \int \|h \nabla\chi + \chi \nabla h\|^2 \, dm \\
&\leq \frac{1+\epsilon}{\gamma} \int \|\nabla h\|^2 \chi^2 \, dm + \frac{1+1/\epsilon}{\gamma} \int \|\nabla\chi\|^2 h^2 \, dm \\
&\leq \frac{1+\epsilon}{\gamma} \int_{\|x - x^{\star}\| \leq \Delta} \|\nabla h\|^2 \, dm + \frac{(1+1/\epsilon) \|h\|_\infty^2}{\gamma} \int \|\nabla\chi\|^2 \, dm.
\end{align*}
}
Since $\chi$ is supported on $\mathcal{B}(x^{\star}, \Delta)$ and $\rho \geq K_{\textup{min}}$ therein, the first integral satisfies
{
\[
\int_{\|x - x^{\star}\| \leq \Delta} \|\nabla h\|^2 \, dm \leq \frac{1}{K_{\textup{min}}} \int \|\nabla h\|^2 \, d\pi.
\]
}
For the gradient of $\chi$, we have $\|\nabla\chi(x)\|^2 \leq 16/\Delta^2$ whenever $\Delta/2 < \|x - x^{\star}\| < \Delta$, and $\nabla\chi = 0$ otherwise. Therefore,
{
\begin{align*}
\mathbb{E}_{m}\|h\chi - \mathbb{E}_{m}(h\chi)\|^2 &\leq \frac{1+\epsilon}{\gamma K_{\textup{min}}} \int \|\nabla h\|^2 \, d\pi + \frac{16(1+1/\epsilon)\|h\|_\infty^2}{\gamma \Delta^2} \Pr_{m}\!\big(\|x - x^{\star}\| > \Delta/2\big).
\end{align*}
}

Combining, we arrive at
{
\begin{align*}
\textup{Var}_{\pi}[h] &\leq \frac{(1+\epsilon)^2 K_{\textup{max}}}{\gamma K_{\textup{min}}} \int \|\nabla h\|^2 \, d\pi + 4\Pr_{\pi}(\|x - x^{\star}\| > \Delta) \, \|h\|_{\infty}^2 \nonumber \\
&\quad + \Big(1 + \frac{1}{\epsilon}\Big)\Big(1 + \frac{16(1+\epsilon)}{\gamma\Delta^2}\Big) K_{\textup{max}} \nonumber \\
&\qquad \times \Pr_{m}(\|x - x^{\star}\| > \Delta/2) \, \|h\|_\infty^2, 
\end{align*}
}
which concludes the proof.
\end{proof}

{
The Gaussian tail estimate below is used in the density-ratio and weak Poincar\'e bounds that follow.

\begin{lemma}[Gaussian tail integrals]
\label{lem:wpi-gaussian-tail}
Let $\ell\geq1$, $a\in\mathbb R^\ell$ and $\lambda>0$, and let $\gamma_\lambda=\mathcal N(a,\lambda^{-1}I_\ell)$. For $r>0$ with $\lambda r^2\geq4\ell$,
\begin{align}
\gamma_\lambda(\|x-a\|>r)
&=\left(\frac{\lambda}{2\pi}\right)^{\ell/2}
\int_{\|x-a\|>r}e^{-\lambda\|x-a\|^2/2}\,dx,\nonumber\\
&\leq e^{-\lambda r^2/8}.
\label{eq:wpi-gaussian-tail}
\end{align}
\end{lemma}
\begin{proof}
For $Z\sim\mathcal N(0,I_\ell)$, Gaussian integration gives $\mathbb E e^{\|Z\|^2/4}=2^{\ell/2}$. Markov's inequality with $t=\lambda r^2\geq4\ell$ yields
\[
\mathbb P(\|Z\|^2>t)\leq2^{\ell/2}e^{-t/4}\leq e^{-t/8}.
\]
The substitution $x=a+\lambda^{-1/2}Z$ gives \eqref{eq:wpi-gaussian-tail}, including the Gaussian normalizing factor.
\end{proof}
}

\begin{lemma}[Density Ratio Bounds Near the Global Minimum]
\label{lemma:bounded-ratio}
Let $f:\mathbb{R}^d \to \mathbb{R}$ be a $C^3$ function with a unique global minimizer $x^{\star}$ satisfying $f(x^{\star})=0$ and $H \coloneqq \nabla^2 f(x^{\star}) \succ 0$, and let $\omega \coloneqq \lambda_{\min}(H) > 0$. Suppose that $f(x) \geq \frac{a}{2}\|x - x^{\star}\|^2$ for some $a > 0$ and all $x \in \mathbb{R}^d$. Define the Gibbs measure $\pi(x) \propto \exp(-\beta f(x))$ and the comparison Gaussian $m(x) \propto \exp(-\beta q(x))$, where $q(x) \coloneqq \frac{1}{2}\langle x - x^{\star}, H(x - x^{\star})\rangle$, so that $m = \mathcal{N}(x^{\star}, (\beta H)^{-1})$. Let $\Delta = \beta^{-5/12}$ and $R_\Delta = \mathcal{B}(x^{\star}, \Delta)$. Then the density ratio $\rho_\beta = d\pi/dm$ satisfies
\[
1 - \mathcal{O}(\beta^{-1/4}) \;\leq\; \rho_\beta(x) \;\leq\; 1 + \mathcal{O}(\beta^{-1/4}), \qquad \forall\, x \in R_\Delta,
\]
for all $\beta$ sufficiently large.
\end{lemma}

\begin{proof}
The density ratio factorizes as
{
\[
\rho_\beta(x) = \frac{d\pi}{dm}(x) = \frac{Z_m}{Z_\pi}\exp\!\big(-\beta(f(x) - q(x))\big).
\]
}

Since $m$ is a Gaussian with covariance $(\beta H)^{-1}$, its normalizing constant is $Z_m = (2\pi/\beta)^{d/2}(\det H)^{-1/2}$.
For $Z_\pi = \int \exp(-\beta f(x))\,dx$, we split the integral over $R_\Delta$ and its complement. On $R_\Delta$, the estimate~\eqref{eq:C2-local-taylor} below for $\exp(-\beta(f(x) - q(x)))$ gives
{
\[
\int_{R_\Delta} \exp(-\beta f(x))\,dx = \big(1 + \mathcal{O}(\beta^{-1/4})\big)\, Z_m \big(1 - \Pr_m(x \notin R_\Delta)\big).
\]
}
Writing $x = x^{\star} + (\beta H)^{-1/2} Z$ with $Z \sim \mathcal{N}(0, I_d)$ under $m$, we have $\|x - x^{\star}\|^2 \leq \|Z\|^2/(\beta\omega)$, so \cref{lem:wpi-gaussian-tail} gives $\Pr_m(x \notin R_\Delta) \leq \Pr(\|Z\|^2 > \omega\beta^{1/6}) \leq \exp(-\omega\beta^{1/6}/8)$ once $\omega\beta^{1/6} \geq 4d$. On the complement of $R_\Delta$, the bound $f(x) \geq \frac{a}{2}\|x - x^{\star}\|^2$ and \cref{lem:wpi-gaussian-tail} applied to $\mathcal{N}(x^{\star}, (\beta a)^{-1} I_d)$ give
{
\begin{equation}
\int_{\|x - x^{\star}\| > \Delta} \exp(-\beta f(x))\,dx \leq Z_m \Big(\frac{\det H}{a^d}\Big)^{1/2} \exp\!\big(-a\beta^{1/6}/8\big) \label{eq:C2-outside-ball}
\end{equation}
}
once $a\beta^{1/6} \geq 4d$. Since $\exp(-\omega\beta^{1/6}/8)$ and $\exp(-a\beta^{1/6}/8)$ are $\mathcal{O}(\beta^{-1/4})$, adding the two parts of $Z_\pi$ gives $Z_\pi = Z_m(1 + \mathcal{O}(\beta^{-1/4}))$. Therefore,
{
\[
\frac{Z_m}{Z_\pi} = 1 + \mathcal{O}(\beta^{-1/4}).
\]
}

For the exponential term, by Taylor's theorem with remainder, since $f(x^{\star})=0$, $\nabla f(x^{\star})=0$, and $\nabla^2 f(x^{\star}) = H$, we have
{
\[
f(x) - q(x) = \mathcal{O}\!\big(\|x - x^{\star}\|^3\big),
\]
}
with $|f(x) - q(x)| \leq \frac{M_3}{6}\|x-x^{\star}\|^3$ for all $x \in \mathcal{B}(x^{\star}, 1)$, where $M_3 = \sup_{\mathcal{B}(x^{\star}, 1)}\|\nabla^3 f\|_{\textup{op}}$ is finite since $f$ is $C^3$. For $x \in R_\Delta$ with $\Delta = \beta^{-5/12} \leq 1$:
{
\[
\beta\,\big|f(x) - q(x)\big| \leq \frac{\beta M_3}{6}\,\Delta^3 = \frac{M_3}{6}\,\beta^{-1/4}.
\]
}
Since $|e^{-u} - 1| \leq 2|u|$ for $|u| \leq 1$, we conclude that for $\beta$ large enough,
{
\begin{equation}
\exp\!\big(-\beta(f(x) - q(x))\big) = 1 + \mathcal{O}(\beta^{-1/4}). \label{eq:C2-local-taylor}
\end{equation}
}

Multiplying the two estimates gives $\rho_\beta(x) = (1+\mathcal{O}(\beta^{-1/4}))(1+\mathcal{O}(\beta^{-1/4})) = 1+\mathcal{O}(\beta^{-1/4})$ for all $x\in R_\Delta$, completing the proof.
\end{proof}
\begin{restatable}[Weak Poincar\'e Bound for Low Temperature Gibbs Distributions]{theorem}{weakpoincarelimit}
\label{thm:weakpoincare-formal}
    Let $f:\mathbb{R}^d \to \mathbb{R}$ be a $C^3$ function with a unique global minimizer $x^{\star}$ and $H \coloneqq \nabla^2 f(x^{\star}) \succ 0$, and let $\omega \coloneqq \lambda_{\min}(H) > 0$. Suppose additionally that $f$ has quadratic growth away from $x^{\star}$, i.e., $f(x) - f(x^{\star}) \geq \frac{a}{2}\|x - x^{\star}\|^2$ for some $a > 0$ and all $x \in \mathbb{R}^d$. Then for all $\beta$ sufficiently large, the Gibbs measure $\pi_\beta(x) \propto \exp(-\beta f(x))$ satisfies the weak Poincar\'e inequality
{
\[
\textup{Var}_{\pi_\beta}[h] \leq \frac{1 + \mathcal{O}(\beta^{-1/4})}{\omega\beta} \int \|\nabla h\|^2 \, d\pi_\beta + C\exp\!\big({-c\,\beta^{1/6}}\big) \, \|h\|_\infty^2,
\]
}
where $c, C > 0$ are constants depending only on $H, a$ and $d$.
\end{restatable}

\begin{proof}
We apply \cref{thm:weak-poincare-perturbation} with the comparison Gaussian $m(x) \propto \exp\!\big(-\frac{\beta}{2}\langle x - x^{\star}, H(x - x^{\star})\rangle\big)$, $\Delta = \beta^{-5/12}$ and $\epsilon = \beta^{-1/4}$. Since adding a constant to $f$ does not change $\pi_\beta$, we may assume $f(x^{\star}) = 0$, so that $f$ satisfies the assumptions of \cref{lemma:bounded-ratio}.

The measure $m$ is $\mathcal{N}(x^{\star}, (\beta H)^{-1})$, which satisfies the Poincar\'e inequality with spectral gap $\gamma = \beta\omega$.

By \cref{lemma:bounded-ratio}, for all $x \in R_{\Delta} = \mathcal{B}(x^{\star}, \Delta)$ (with $\Delta = \beta^{-5/12}$, which is still $o(1)$ as $\beta \to \infty$), we have
\[
K_{\textup{min}} = 1 - \mathcal{O}(\beta^{-1/4}), \qquad K_{\textup{max}} = 1 + \mathcal{O}(\beta^{-1/4}).
\]

Substituting these bounds into the Dirichlet coefficient from \cref{thm:weak-poincare-perturbation} gives
{
\[
\frac{(1+\epsilon)^2 K_{\textup{max}}}{\gamma K_{\textup{min}}} = \frac{(1+\beta^{-1/4})^2(1 + \mathcal{O}(\beta^{-1/4}))}{\beta\omega(1 - \mathcal{O}(\beta^{-1/4}))} = \frac{1 + \mathcal{O}(\beta^{-1/4})}{\omega\beta}.
\]
}

Since $\beta\omega\|x - x^{\star}\|^2 \leq \|Z\|^2$ under $m$, where $x = x^{\star} + (\beta H)^{-1/2} Z$ with $Z \sim \mathcal{N}(0, I_d)$, we have
{
\[
\Pr_m(\|x - x^{\star}\| > \Delta/2) \leq \Pr\!\big(\|Z\|^2 > \beta\omega\Delta^2/4\big) = \Pr\!\big(\|Z\|^2 > \omega\beta^{1/6}/4\big).
\]
}
By \cref{lem:wpi-gaussian-tail}, for $\beta$ large enough that $\omega\beta^{1/6} \geq 16d$, this gives $\Pr_m(x \notin R_{\Delta/2}) \leq \exp(-\omega\beta^{1/6}/32)$.

Since $f(x) \geq \frac{a}{2}\|x - x^{\star}\|^2$ globally, the Gibbs measure $\pi_\beta$ is dominated by a Gaussian in the tails. In particular, the normalizing constants of $\pi_\beta$ and $m$ satisfy $Z_\pi = Z_m(1 + \mathcal{O}(\beta^{-1/4})) \geq Z_m/2$ for $\beta$ sufficiently large by the proof of \cref{lemma:bounded-ratio}, and dividing~\eqref{eq:C2-outside-ball} by $Z_\pi$ gives $\Pr_{\pi_\beta}(\|x - x^{\star}\| > \Delta) \leq 2(\det H/a^d)^{1/2}\exp(-a\beta^{1/6}/8)$ once $a\beta^{1/6} \geq 4d$.

The coefficient of $\|h\|_\infty^2$ in \cref{thm:weak-poincare-perturbation} becomes
{
\[
4\Pr_{\pi_\beta}(x \notin R_{\Delta}) + \Big(1 + \frac{1}{\epsilon}\Big)\Big(1 + \frac{16(1+\epsilon)}{\gamma\Delta^2}\Big) K_{\textup{max}} \Pr_m(x \notin R_{\Delta/2}).
\]
}
Since $\gamma\Delta^2 = \beta\omega \cdot \beta^{-5/6} = \omega\beta^{1/6}$ and $\epsilon = \beta^{-1/4}$, the prefactor satisfies
{
\[
\big(1 + \beta^{1/4}\big)\Big(1 + \frac{16(1+\beta^{-1/4})}{\omega\beta^{1/6}}\Big) K_{\textup{max}} \leq 3\beta^{1/4}
\]
}
for $\beta$ sufficiently large. Therefore, the full $\|h\|_\infty^2$ coefficient is at most $8(\det H/a^d)^{1/2}\exp(-a\beta^{1/6}/8) + 3\beta^{1/4}\exp(-\omega\beta^{1/6}/32)$. Absorbing the factor $\beta^{1/4}$ into the exponential shows that this coefficient is bounded by $C\exp(-c\beta^{1/6})$ with $c = \min\{a/8, \omega/64\}$, concluding the proof.
\end{proof}

{
We bound the initial density ratio and the Gibbs tail probabilities, and apply \cref{thm:weak-poincare-perturbation} to obtain the weak Poincar\'{e} inequality required by \cref{thm:annealed-ula-convergence}.

\begin{lemma}[Density ratios, tail probabilities, and a weak Poincar\'{e} inequality]
\label{lem:ula-block-estimates}
Let $g\in C^2(\mathbb{R}^{\ell})$, where $1\leq\ell\leq b$, have an $L$-Lipschitz gradient and a minimizer $a$ with $g(a)=0$, $\|a\|\leq B_0$, and $g(y)\geq m\|y-a\|^2/2$. For $\beta>0$, define the following distributions and constants:
\begin{align}
Z_\beta&:=\int_{\mathbb R^\ell}e^{-\beta g(y)}\,dy,\nonumber\\
\pi_\beta(y)&:=Z_\beta^{-1}e^{-\beta g(y)},
&\gamma_0&:=\mathcal N(0,(2L)^{-1}I_\ell).\nonumber\\
\kappa&:=(L/m)^{b/2}, & A&:=(2L/m)^{b/2},\nonumber\\
K_0&:=Ae^{LB_0^2}, & U&:=1+\log K_0.
\label{eq:ula-block-constants}
\end{align}
These choices give the bounds
\begin{align}
\frac m2\|y-a\|^2&\leq g(y)\leq\frac L2\|y-a\|^2,
\label{eq:ula-quadratic-bounds}\\
\left\|\frac{d\gamma_0}{d\pi_1}\right\|_\infty&\leq K_0,
\label{eq:ula-initial-ratio}\\
\frac{Z_{\beta/2}}{Z_\beta}&\leq A,\qquad
\pi_\beta(g>u)\leq Ae^{-\beta u/2}\quad(\beta,u>0).
\label{eq:ula-energy-tail}
\end{align}
For $s\in(0,1]$, define
\begin{equation}
\begin{aligned}
S&:=\max\{32,16b,32\log(8\kappa/s)\},\\
\delta&:=\frac{m}{4\kappa}e^{-(L/m-1)S/2}.
\end{aligned}
\label{eq:ula-wpi-parameters}
\end{equation}
The generator Dirichlet form $\mathcal E_{\pi_\beta}(v,v)=\beta^{-1}\int\|\nabla v\|^2d\pi_\beta$ satisfies
\begin{equation}
\mathcal E_{\pi_\beta}(v,v)\geq
\delta\operatorname{Var}_{\pi_\beta}(v)-\delta s\|v\|_\infty^2.
\label{eq:ula-uniform-wpi}
\end{equation}
For fixed $b,m,L$, the inverse of $\delta$ is polynomial in $1/s$.
\end{lemma}
\begin{proof}
Taylor's theorem gives \eqref{eq:ula-quadratic-bounds}. Gaussian integration therefore bounds $Z_\beta$ between $(2\pi/(\beta L))^{\ell/2}$ and $(2\pi/(\beta m))^{\ell/2}$. Completing the square gives
\begin{equation*}
g(y)-L\|y\|^2\leq L\|a\|^2-\frac L2\|y+a\|^2,
\end{equation*}
which proves \eqref{eq:ula-initial-ratio}. The ratio of the Gaussian normalizer bounds is at most $A$. Applying Markov's inequality to $e^{\beta g/2}$ under $\pi_\beta$ proves \eqref{eq:ula-energy-tail}.

For \eqref{eq:ula-uniform-wpi}, apply \cref{thm:weak-poincare-perturbation} with comparison measure $m_\beta=\mathcal N(a,(\beta m)^{-1}I_\ell)$, comparison parameter $\epsilon=1$, and radius $r_\beta=\sqrt{S/(\beta m)}$. Equation \eqref{eq:ula-quadratic-bounds} gives the required density-ratio bounds:
\begin{equation}
\begin{aligned}
\frac{d\pi_\beta}{dm_\beta}&\leq\kappa,
&K_{\mathrm{min}}&\geq e^{-(L/m-1)S/2},\\
\gamma r_\beta^2&=S,
\end{aligned}
\label{eq:ula-comparison-parameters}
\end{equation}
where $\gamma=\beta m$. Applying \cref{lem:wpi-gaussian-tail} to the comparison Gaussian at radii $r_\beta$ and $r_\beta/2$ bounds the defect in that comparison by
\begin{equation}
4\kappa e^{-S/8}+2(1+32/S)\kappa e^{-S/32}
\leq8\kappa e^{-S/32}\leq s.
\label{eq:ula-wpi-defect}
\end{equation}
Substituting \eqref{eq:ula-comparison-parameters} and \eqref{eq:ula-wpi-defect} into \cref{thm:weak-poincare-perturbation}, and using the normalization \eqref{eq:wpi-variance-form}, proves \eqref{eq:ula-uniform-wpi}. The dependence on $s$ follows directly from \eqref{eq:ula-wpi-parameters}.
\end{proof}
}

At low temperature, quadratic growth of a block potential gives a weak Poincar\'e inequality with an exponentially small defect.

{
\begin{corollary}[Weak Poincar\'e Bound with Exponentially Small Defect]
\label{cor:wpi-exp-defect}
Let $\mathcal{X}' \subseteq \mathbb{R}^{d'}$ be a convex domain, and let $g \in C^2(\mathcal{X}')$ have a unique global minimizer $x^{\star}$, which lies in the interior of $\mathcal{X}'$, and satisfy $g(x) - g(x^{\star}) \geq \frac{m}{2}\|x - x^{\star}\|^2$ for some $m > 0$ and all $x \in \mathcal{X}'$. Then there exist constants $\beta_0, \xi, c, C > 0$, depending only on $g$ and $\mathcal{X}'$, such that for all $\beta \geq \beta_0$, the Gibbs measure $\pi_\beta \propto \exp(-\beta g)$ on $\mathcal{X}'$ satisfies
\[
\textup{Var}_{\pi_\beta}[h] \leq \frac{4}{\xi\beta}\int \|\nabla h\|^2 \, d\pi_\beta + Ce^{-c\beta}\|h\|_\infty^2
\]
for all bounded smooth $h$.
\end{corollary}
\begin{proof}
Adding a constant to $g$ does not change $\pi_\beta$, so we may assume $g(x^{\star}) = 0$. Choose $r_0 > 0$ and $\xi > 0$ such that the closed ball $B := \{x : \|x - x^{\star}\| \leq r_0\}$ lies in $\mathcal{X}'$ and $\nabla^2 g \succeq \xi I$ on $B$, and let $L_0 := \max_{x \in B}\|\nabla^2 g(x)\|$. We compare $\pi_\beta$ with the Gibbs measure of the function
\[
\tilde{g}(x) := \sup_{y \in B}\Big\{g(y) + \langle \nabla g(y), x - y \rangle + \frac{\xi}{2}\|x - y\|^2\Big\}, \qquad x \in \mathbb{R}^{d'}.
\]
As a supremum of $\xi$-strongly convex functions, $\tilde{g}$ is $\xi$-strongly convex. Since $\nabla^2 g \succeq \xi I$ on the convex set $B$, we have $\tilde{g} = g$ on $B$, and the choice $y = x^{\star}$, at which $g$ and $\nabla g$ vanish, gives $\tilde{g}(x) \geq \frac{\xi}{2}\|x - x^{\star}\|^2$.

Let $m_\beta \propto \exp(-\beta\tilde{g})$ on $\mathcal{X}'$, and let $Z_m$ and $Z_\pi$ be the normalizing constants of $m_\beta$ and $\pi_\beta$. Since $\beta\tilde{g}$ is $\beta\xi$-strongly convex and $\mathcal{X}'$ is convex, $m_\beta$ satisfies the Poincar\'e inequality with spectral gap $\gamma = \beta\xi$ by the Bakry--\'Emery criterion~\citep{Bakry2013AnalysisAG}, applied to smooth approximations of $\tilde{g}$ if $\tilde{g}$ is not smooth. The density ratio $d\pi_\beta/dm_\beta = \frac{Z_m}{Z_\pi}e^{-\beta(g - \tilde{g})}$ equals the constant $Z_m/Z_\pi$ on $B$.

We apply \cref{thm:weak-poincare-perturbation} to $\pi_\beta$ with the reference measure $m_\beta$, with $\epsilon = 1$ and with $\Delta = r_0$, so that $R_\Delta = B$ and $K_{\textup{max}} = K_{\textup{min}} = Z_m/Z_\pi$. The coefficient of $\int \|\nabla h\|^2 \, d\pi_\beta$ is then $4/(\beta\xi)$, and the coefficient of $\|h\|_\infty^2$ is
\begin{align*}
&4\Pr_{\pi_\beta}(\|x - x^{\star}\| > r_0) \\
&\qquad + 2\frac{Z_m}{Z_\pi}\Big(1 + \frac{32}{\beta\xi r_0^2}\Big)\Pr_{m_\beta}(\|x - x^{\star}\| > r_0/2).
\end{align*}

It remains to bound the two probabilities and the ratio $Z_m/Z_\pi$. Let $\phi = g$ and $a = m$, or $\phi = \tilde{g}$ and $a = \xi$, and let $P_\phi \propto e^{-\beta\phi}$ on $\mathcal{X}'$, so that $P_g = \pi_\beta$ and $P_{\tilde{g}} = m_\beta$. Then $\phi(x) \geq \frac{a}{2}\|x - x^{\star}\|^2$ on $\mathcal{X}'$, and $\phi(x) \leq \frac{L_0}{2}\|x - x^{\star}\|^2$ on $B$ by Taylor's theorem, since $\phi = g$ on $B$. For $s \in \{r_0, r_0/2\}$, \cref{lem:wpi-gaussian-tail} bounds the numerator of $P_\phi(\|x - x^{\star}\| > s)$ by the Gaussian integral of $e^{-\beta a\|x - x^{\star}\|^2/2}$ over $\{\|x - x^{\star}\| > s\}$, and its denominator is bounded below by the integral of $e^{-\beta L_0\|x - x^{\star}\|^2/2}$ over $B$, which is at least half of the corresponding Gaussian integral over $\mathbb{R}^{d'}$ for $\beta$ large. This gives, for $\beta$ large enough,
\[
P_\phi(\|x - x^{\star}\| > s) \leq 2\Big(\frac{L_0}{a}\Big)^{d'/2}\exp\Big(-\frac{\beta a s^2}{8}\Big).
\]
For $\phi = \tilde{g}$ and $s = r_0$, the right-hand side is at most $1/2$ for $\beta$ large, so $m_\beta(B) \geq 1/2$. Since $\tilde{g} = g$ on $B$, it follows that $Z_\pi \geq \int_B e^{-\beta g} = Z_m\, m_\beta(B) \geq Z_m/2$, that is, $Z_m/Z_\pi \leq 2$.

Inserting these bounds, with $(\phi, s) = (g, r_0)$ in the first term and $(\phi, s) = (\tilde{g}, r_0/2)$ in the second, bounds the coefficient of $\|h\|_\infty^2$ by $Ce^{-c\beta}$ for all $\beta \geq \beta_0$, where $c = \min\{m r_0^2/8, \xi r_0^2/32\}$, $\beta_0$ is large enough for the bounds above, and $C$ depends only on $d'$, $m$, $\xi$, $L_0$, $r_0$ and $\beta_0$.
\end{proof}
}

\section{Analysis of the Unadjusted Langevin Algorithm}

\label{sec:ula-proofs}

We consider an annealed version of Langevin dynamics of the form 
\begin{align}
    \label{eq:annealed-sde}
    dX_t = -\nabla f(X_t) d t + \sqrt{\eta(\theta_t)} d B_t,
\end{align}
where $\theta : [0,T] \to \mathbb{R}_{>0}$ is a learning rate schedule.
We denote the law of this process at time $T$ by $\mu_T$. Fix a decreasing sequence of learning rates ${\eta}_1 > {\eta}_2 > \cdots > {\eta}_M > 0$. In phase $k = 1, \ldots, M$, we run with constant learning rate ${\eta}_k$
for a continuous-time duration $\tau_k$ (i.e. a piecewise-constant schedule). The duration $\tau_k$ of each phase is chosen so that the Langevin SDE at constant
learning rate ${\eta}_k$,
\begin{align}
    \label{eq:phase-sde}
    dX_t = -\nabla f(X_t)\dd t + \sqrt{{\eta}_k}\,dB_t,
\end{align}
has approximately reached its stationary distribution $\pi_{{\eta}_k}
\propto \exp(-2f/{\eta}_k)$. Formally, the learning rate schedule is piecewise constant: if $T_k = \sum_{m=1}^{k} \tau_m$, then $\forall t \in [T_{k-1}, T_{k}], \eta(\theta_t) = \eta_k$.

We will also require the following weak condition on $f$.
\begin{definition}[Confining Function]
\label{defn:confining_func}
    A function $\phi$ is said to be confining  when $\lim_{\|x\|\to \infty} \phi(x) = \infty$ and for all $s>0$,
\begin{equation*}
  \int_{\mathbb{R}^d}e^{-\frac{2\phi(x)}{s}}\mathrm{d}x<+\infty.  
\end{equation*}
\end{definition}

It is known that for confining functions that are Morse, the spectral gap of Langevin dynamics can be estimated as follows.
\begin{restatable}[Theorem 2 \citep{shi2023learning}]{theorem}{langsgdgap}
\label{thm:lang_sgd_gap}
Suppose $f$ is {a confining function} (see Definition \ref{defn:confining_func}), \begin{align*}
    \forall \beta > 0, \lim_{\lVert x \rVert \rightarrow \infty}\left(\beta\norm{\nabla f(x)}^2 - \Delta f(x)\right) = \infty,
\end{align*} $f$ has at least two local minima, and $f$ is a Morse function. Then there exists a $\beta_0$, depending on $f$, such that $\forall \beta > \beta_0$, the {reciprocal of the }relaxation time $\delta^{(C)}(\beta)$ of Langevin dynamics {with potential} $f$ satisfies
\begin{align*}
    \delta^{(C)}(\beta) =(\alpha + o(1))e^{-\beta H_f},
\end{align*}
where $H_f$ is a quantity known as the Morse saddle barrier that depends only on the geometry of the function $f$, and $\alpha >0$ is a constant depending on $f$.
\end{restatable}

We start with the following result bounding the continuous-time convergence of Langevin dynamics in \Renyi{} divergence under WPI.

We use the generator-normalized WPI convention \eqref{eq:wpi-convention} from \cref{sec:weak-poincare-proofs}. When $\mathcal{X}$ is a bounded convex set, as in \cref{thm:lang_contin_convergence}, the diffusion \eqref{eq:phase-sde} is reflected at the boundary of $\mathcal{X}$.

We state \cref{lem:renyi_convergence_ld_wpi} and the convergence bounds that follow it in terms of the $q$-th \Renyi{} divergence of a probability measure $\nu$ with respect to a probability measure $\pi$, defined as $\mathcal{R}_q(\nu || \pi) = \frac{1}{q-1}\log\left[\int \left(\frac{d\nu}{d\pi}\right)^q d\pi\right]$.

\begin{lemma}
\label{lem:renyi_convergence_ld_wpi}
Let $f(x) : \mathcal{X} \rightarrow \mathbb{R}$. Consider running  Langevin dynamics generated by $f$ for time $T$ producing distribution $\mu_T$. Let $\pi$ be the  stationary distribution of the Langevin dynamics. Suppose that $K := \sup_{x \in \mathcal{X}} \frac{d\mu_{0}}{d\pi}(x) < \infty$ and let $q \geq 2$. If $\pi$ satisfies a $(\delta, \epsilon)$-WPI with $\epsilon K^q \leq \delta/4$, then $u(t) := \mathcal{R}_q(\mu_t || \pi)$ satisfies
{
\begin{align}
u(T) &\leq \max\left\{u(0) - \frac{\delta T}{q}, 1\right\},
\label{eq:renyi-wpi-decrease}\\
\intertext{and, for $T \geq T_0 := q(u(0) - 1)_+/\delta$,}
u(T) &\leq \exp\left(-\frac{2\delta (T - T_0)}{q}\right) + \frac{2\epsilon K^q}{\delta}.
\label{eq:renyi-wpi-convergence}
\end{align}
}

\end{lemma}
\begin{proof}
The following is a modified version of the proof of Theorem~2.2.25 in \citet{sinhobook}, originally from \citep{vempala2019rapid}, to account for a weak Poincar\'e inequality.

Let $\rho_0 := \frac{d\mu_0}{d\pi}$ and $P_t$ be the transition operator of Langevin dynamics. Then from reversibility of Langevin dynamics $\rho_t = P_t\rho_0$. %
Taking the derivative gives

{
\begin{align}
\label{eq:renyi-derivative}
\frac{d}{dt}\mathcal{R}_q(\mu_{t} || \pi) = -\frac{4}{q}\frac{\mathcal{E}(\rho_t^{q/2}, \rho_t^{q/2})}{\int \rho_t^q d\pi}.
\end{align}
}

The $(\delta, \epsilon)$-weak Poincar\'e inequality implies that
{
\begin{align*}
\mathcal{E}_{\pi}(\rho_t^{q/2}, \rho_t^{q/2}) &\geq  \delta\textup{Var}_{\pi}\left(\rho_t^{q/2}\right) - \textup{Error}_{t} \\
\textup{Error}_t &= \epsilon \lVert \rho_t \rVert_{\infty}^{q}.\\
\intertext{Note that}
\lVert \rho_t \rVert_{\infty} &\leq K,
\end{align*}
}  
where the inequality holds because $\rho_t = P_t\rho_0$ and the Markov operator $P_t$ is a contraction in the supremum norm. Hence $\textup{Error}_t \leq \epsilon K^q$.

Combining these bounds gives
{
\begin{align*}
\delta^{-1}\mathcal{E}_{\pi}(\rho_t^{q/2}, \rho_t^{q/2}) &\geq  \textup{Var}_{\pi}\left(\rho_t^{q/2}\right) - \delta^{-1}\textup{Error}_t\\
&= \int (\rho_t)^q d\pi - \left(\int (\rho_t)^{q/2} d\pi\right)^2 - \delta^{-1}\textup{Error}_t\\
&=\left(\int (\rho_t)^q d\pi\right)\left(1 - \frac{\exp((q-2)\mathcal{R}_{q/2}(\rho_t\pi || \pi))}{\exp((q-1)\mathcal{R}_{q}(\rho_t\pi || \pi))}\right) - \delta^{-1}\textup{Error}_{ t}\\
&\geq \left(\int \rho_t^q d\pi\right)\left(1 - {\exp(-\mathcal{R}_{q}(\rho_t\pi || \pi))}\right) - \delta^{-1}\textup{Error}_{t},
\end{align*}
}
where the last inequality follows from the monotonicity of \Renyi{} divergence in $q$.
Plugging this into the derivative of \Renyi{} divergence in \eqref{eq:renyi-derivative}, we get
{
\begin{align*}
    -\frac{d}{dt}\mathcal{R}_{q}(\mu_{t} || \pi) &\geq \frac{\delta 4}{q}(1 - {\exp(-\mathcal{R}_{q}(\rho_t\pi || \pi))}) - \frac{4}{q\int \rho_t^q d\pi}\textup{Error}_{t}\\
    &=\frac{\delta 4}{q}(1 - {\exp(-\mathcal{R}_{q}(\mu_t || \pi))}) - \frac{4}{q}\textup{Error}_{t}\exp(-(q-1)\mathcal{R}_{q}(\mu_{t} || \pi)).
\end{align*}
}

Since $F(t) := \int \rho_t^q d\pi = \exp((q-1)u(t)) \geq 1$ and $\textup{Error}_t \leq \epsilon K^q$, it follows that
{
\begin{align*}
u'(t) \leq -\frac{4\delta}{q}\left(1 - \exp(-u(t))\right) + \frac{4\epsilon K^q}{q}.
\end{align*}
}
If $u(t) \geq 1$, then $1 - \exp(-u(t)) \geq 1/2$, and the condition $\epsilon K^q \leq \delta/4$ gives $u'(t) \leq -\delta/q$. Hence $u$ decreases at rate at least $\delta/q$ while $u \geq 1$, and once $u \leq 1$ it stays at most $1$. This proves \eqref{eq:renyi-wpi-decrease} and shows that $u(t) \leq 1$ for all $t \geq T_0$. If $0 \leq u(t) \leq 1$, then $1 - \exp(-u(t)) \geq u(t)/2$, so
{
\begin{align*}
u'(t) \leq -\frac{2\delta}{q}u(t) + \frac{4\epsilon K^q}{q}.
\end{align*}
}
Applying Gr\"onwall's inequality on $[T_0, T]$ with $u(T_0) \leq 1$ gives \eqref{eq:renyi-wpi-convergence}.

\end{proof}

\cref{lem:renyi_convergence_ld_wpi} requires a bounded initial density ratio $K$, while the proof of \cref{thm:lang_contin_convergence} controls only a \Renyi{} divergence $\mathcal{R}_q$ of finite order $q$, which gives no bound on the density ratio. The following elementary lemma restricts a measure to a set of large mass on which its density ratio is bounded, and bounds the \Renyi{} divergence of the normalized restriction.

{
\begin{lemma}
\label{lem:truncation}
Let $\nu$ and $\pi$ be probability measures on $\mathcal{X}$, let $s > 1$, $U \geq 0$ and $\varrho \in (0, 1/2]$, and suppose that $\mathcal{R}_s(\nu || \pi) \leq U$. Let $S := \{x \in \mathcal{X} : \frac{d\nu}{d\pi}(x) \leq \varrho^{-1/(s-1)}e^{U}\}$ and $\tilde{\nu} := \nu(\cdot \cap S)/\nu(S)$. Then
\begin{align*}
    \nu(S^c) \leq \varrho, \qquad \sup_{x \in \mathcal{X}}\frac{d\tilde{\nu}}{d\pi}(x) \leq \frac{\varrho^{-1/(s-1)}e^{U}}{1-\varrho}, \nonumber\\
    \mathcal{R}_s(\tilde{\nu} || \pi) \leq U + \frac{s}{s-1}\log\frac{1}{1-\varrho}. 
\end{align*}
\end{lemma}
\begin{proof}
Write $g := \frac{d\nu}{d\pi}$. By the definition of $\mathcal{R}_s$, the expectation of $g^{s-1}$ under $\nu$ is $\int g^{s} d\pi = e^{(s-1)\mathcal{R}_s(\nu || \pi)} \leq e^{(s-1)U}$. Markov's inequality under $\nu$ applied to $g^{s-1}$ therefore gives $\nu(S^c) = \nu(g^{s-1} > \varrho^{-1}e^{(s-1)U}) \leq \varrho$, and hence $\nu(S) \geq 1 - \varrho$.

By the normalization of $\tilde{\nu}$ and this lower bound, its density with respect to $\pi$ is $\frac{d\tilde{\nu}}{d\pi} = g\mathbf{1}_S/\nu(S) \leq \varrho^{-1/(s-1)}e^{U}/(1-\varrho)$, and $\int \left(\frac{d\tilde{\nu}}{d\pi}\right)^s d\pi \leq \nu(S)^{-s}\int g^s d\pi \leq (1-\varrho)^{-s}e^{(s-1)U}$. Taking the logarithm and dividing by $s-1$ gives the bound on $\mathcal{R}_s(\tilde{\nu} || \pi)$.
\end{proof}
}

\subsection{Convergence of Annealed Langevin Diffusion}

We use \cref{lem:renyi_convergence_ld_wpi,lem:truncation} to bound the continuous-time convergence of annealed Langevin dynamics.

\begin{theorem}[Convergence of annealed Langevin diffusion]
\label{thm:lang_contin_convergence}
    Suppose the sequence of learning rates $\{ \eta_{r} \}_{r=1}^{M}$ is such that the corresponding Langevin dynamics starting from $\mu_0$ with constant rate $\eta_{r}$ has stationary distribution $\pi^{(r)} \propto \exp(-2f/\eta_r)$  and satisfies a $(\delta_{r}, \epsilon_r)$ WPI. Let $q \geq 2$ and $p \in (0, 1/4]$, suppose that $\bar{R} := \max\{1, \mathcal{R}_q(\mu_0 || \pi^{(1)})\} < \infty$, and let $\tau_r = \frac{2q}{\delta_r}$.
    Then if $\Delta_r \operatorname{osc}(f) \leq 1$ and $\epsilon_r \leq \frac{p^2\delta_r}{4M^2}e^{-q(\bar{R}+2)}$ for each $r = 2, \ldots, M$, where $\Delta_r := 2/\eta_r - 2/\eta_{r-1}$ and $\operatorname{osc}(f) := \sup_{x \in \mathcal{X}} f(x) - \inf_{x \in \mathcal{X}} f(x)$, then
  {
\begin{align*}
            \mu_T = \alpha\nu + (1-\alpha)\bar{\nu}, \qquad \alpha \geq 1 - p, \qquad \mathcal{R}_q(\nu || \pi^{(M)}) \leq \bar{R},
  \end{align*}
}
    for some $\alpha \in [0, 1]$ and probability measures $\nu$ and $\bar{\nu}$ on $\mathcal{X}$, where $T = \sum_{r=2}^{M}\tau_r$ and $\mu_T$ is the output of Langevin dynamics after running for time $T$.
\end{theorem}
\begin{proof}
  The proof of this result is inspired by the proof of Theorem 4.13 \citep{huang2025weakpoincare}. 

Let $\pi^{(r)} \propto \exp(-2f/\eta_r)$, corresponding to the stationary distribution for the $r$-th chain. As mentioned earlier $\eta(t)$ is such that for all $t \in [T_{r-1}, T_{r})$, $\eta(t) = \eta_r$. Also let $T_1 = 0$ and $T_r = T_{r-1} + \tau_r$ for $r = 2, \ldots, M$, so that $T_M = T$. The distribution $\mu_t$ corresponds to the output of the annealed Langevin chain at time $t \in [0, T]$. Recall that $\Delta_r = 2/\eta_r - 2/\eta_{r-1}$.

We start by bounding the following quantities for consecutive Gibbs measures. 
Note that for two consecutive Gibbs distributions,
{
\begin{align*}
    \sup_{x \in \mathcal{X}}\frac{d\pi^{(r-1)}}{d\pi^{(r)}}(x) = \sup_{x \in \mathcal{X}}\frac{e^{\Delta_r f(x)}}{\mathbb{E}_{\pi^{(r)}}[e^{\Delta_r f}]} \leq e^{\Delta_r \operatorname{osc}(f)}.
\end{align*}
}
Hence $D_r := \mathcal{R}_{\infty}(\pi^{(r-1)} || \pi^{(r)}) \leq \Delta_r \operatorname{osc}(f) \leq 1$, where $\mathcal{R}_{\infty}(\pi_1 || \pi_2) := \log \sup_{x \in \mathcal{X}}\frac{d\pi_1}{d\pi_2}(x)$. We also have the following inequalities for every probability measure $\mu$ on $\mathcal{X}$ and every $s > 1$:
{
\begin{align}
    \mathcal{R}_{s}(\mu || \pi^{(r)}) &\leq \mathcal{R}_{s}(\mu || \pi^{(r-1)}) + D_r, \nonumber\\
    \sup_{x \in \mathcal{X}}\frac{d\mu}{d\pi^{(r)}}(x) &\leq e^{D_r}\sup_{x \in \mathcal{X}}\frac{d\mu}{d\pi^{(r-1)}}(x). \label{eqn:infinity_overlap}
\end{align}
}

We now show by induction on $r = 1, \ldots, M$ that each chain is acting as a warm start for the next, in the sense that the law at time $T_r$ decomposes as
{
\begin{align}
    \mu_{T_r} = \alpha_r\nu_r + (1-\alpha_r)\bar{\nu}_r, \qquad \alpha_r \geq 1 - (r-1)\rho, \nonumber\\
    \mathcal{R}_q(\nu_r || \pi^{(r)}) \leq \bar{R}, \label{eq:good-bad-decomposition}
\end{align}
}
for some $\alpha_r \in [0, 1]$ and probability measures $\nu_r$ and $\bar{\nu}_r$ on $\mathcal{X}$, where $\rho := p/M$. For $r = 1$, this holds with $\nu_1 = \mu_0$ and $\alpha_1 = 1$ by the definition of $\bar{R}$.

Now let $2 \leq r \leq M$ and suppose that \eqref{eq:good-bad-decomposition} holds for $r-1$. Applying \cref{lem:truncation} to $\nu_{r-1}$ and $\pi^{(r-1)}$ with $s = q$, $U = \bar{R}$ and $\varrho = \rho \leq 1/4$ gives a set $S$ with $\nu_{r-1}(S^c) \leq \rho$ such that $\tilde{\nu} := \nu_{r-1}(\cdot \cap S)/\nu_{r-1}(S)$ satisfies $\sup_{x \in \mathcal{X}}\frac{d\tilde{\nu}}{d\pi^{(r-1)}}(x) \leq K := \rho^{-1/(q-1)}e^{\bar{R}}/(1-\rho)$ and $\mathcal{R}_q(\tilde{\nu} || \pi^{(r-1)}) \leq \bar{R} + 1$, where the last bound uses $q/(q-1) \leq 2$ and $\rho \leq 1/4$.

Changing the reference measure to the stationary distribution $\pi^{(r)}$ of the next chain, using \eqref{eqn:infinity_overlap} and $D_r \leq 1$, gives $\sup_{x \in \mathcal{X}}\frac{d\tilde{\nu}}{d\pi^{(r)}}(x) \leq eK$ and $\mathcal{R}_q(\tilde{\nu} || \pi^{(r)}) \leq \bar{R} + 2$.

We apply \cref{lem:renyi_convergence_ld_wpi} to the Langevin dynamics with constant learning rate $\eta_r$ started from $\tilde{\nu}$. Its hypothesis holds, because the density ratio of $\tilde{\nu}$ with respect to $\pi^{(r)}$ is at most $eK$, and $\rho^{-q/(q-1)} \leq \rho^{-2}$, $\rho \leq 1/4$ and the assumption on $\epsilon_r$ give $\epsilon_r(eK)^q \leq \frac{\delta_r}{4}\left(\frac{4}{3e}\right)^q \leq \frac{\delta_r}{4}$. Write $\mu P^{(r)}_t$ for the law at time $t$ of this dynamics started from a probability measure $\mu$, and let $\nu_r := \tilde{\nu}P^{(r)}_{\tau_r}$. Since $\mathcal{R}_q(\tilde{\nu} || \pi^{(r)}) \leq \bar{R} + 2$, $\delta_r\tau_r/q = 2$ and $\bar{R} \geq 1$, \eqref{eq:renyi-wpi-decrease} gives
{
\begin{align*}
    \mathcal{R}_q(\nu_r || \pi^{(r)}) \leq \max\left\{\bar{R} + 2 - \frac{\delta_r\tau_r}{q}, 1\right\} = \bar{R}.
\end{align*}
}

The law $\mu P^{(r)}_{\tau_r}$ depends linearly on $\mu$, so $\mu_{T_r} = \mu_{T_{r-1}}P^{(r)}_{\tau_r}$ is the sum of $\alpha_{r-1}\nu_{r-1}(S)\,\nu_r$ and the image under $P^{(r)}_{\tau_r}$ of the remaining mass $\alpha_{r-1}\nu_{r-1}(\cdot \cap S^c) + (1-\alpha_{r-1})\bar{\nu}_{r-1}$. Hence \eqref{eq:good-bad-decomposition} holds for $r$ with $\alpha_r := \alpha_{r-1}\nu_{r-1}(S)$, which is at least $1 - (r-1)\rho$ because $\alpha_{r-1} \geq 1 - (r-2)\rho$ and $\nu_{r-1}(S) \geq 1 - \rho$, and with $\bar{\nu}_r$ the normalized image of the remaining mass.

With the induction complete, taking $r = M$ and using $T_M = T$ gives the claim with $\nu = \nu_M$, $\bar{\nu} = \bar{\nu}_M$ and $\alpha = \alpha_M \geq 1 - (M-1)\rho \geq 1 - p$.

\end{proof}

\subsection{Discretization of Langevin Diffusion}

To discretize the annealed Langevin dynamics in time, consider a sequence of step sizes $h_1, \dots, h_M$, one for each phase of the piecewise-constant learning rate $\eta(\theta_t)$. We consider a schedule of step sizes $h(t)$ such that $\forall t \in [T_{k-1}, T_{k}], h(t) = h_k$. Then,  the first-order
discretization of~\eqref{eq:annealed-sde}, evaluated at the discrete times
$t_n = \sum_{m=0}^{n-1} h(t_m)$, takes the form
{
\begin{align}
    \label{eq:annealed-sgd}
    {X}_{t_{n+1}} = {X}_{t_n} - h(t_n) \nabla f(X_{t_n}) + \sqrt{h(t_n)\,\eta(\theta_{t_n})} \, Z_n,
    \qquad Z_n \sim \mathcal{N}(0, I_d),
\end{align}
}
which is referred to as the unadjusted Langevin algorithm (ULA). Since the total time is  $T = \sum_{k=1}^M \tau_k$, for  $N_k = \tau_k/h_k$, the total 
number of steps of the discrete process is $N = \sum_{k=1}^M N_k$. If we set $h_k = \eta_k$, then ULA corresponds to  annealed SGD (Definition \ref{def:annealed-sgd-intro}). %

We write $\hat{\mu}_t$ for the law of the discretized process at time $t$. To bound the discretization error, we compare the law $\mathrm{W}_{T}$ of the path, up to time $T$, of the continuous SDE
\begin{align*}
    dX_t = -\nabla f(X_t) d t + \sqrt{\eta(\theta_t)} d B_t,
\end{align*}
and the law $P_{T}$ of the path of the SDE that corresponds to the discretized dynamics (ULA)
\begin{align*}
    dX_t = -\nabla f(X_{t^-}) d t + \sqrt{\eta(\theta_t)} d \tilde{B}_t,
\end{align*}
where $t^- = t_n$ for $t \in [t_n, t_{n+1})$. Note that $\theta$ does not need to be evaluated at $t^-$ since the learning rate schedule is piecewise constant.

To keep things general, for the result below, we will let $h_k$ and $\eta_k$ be independent. Later on we will discuss when we can set $h_k = \eta_k$. Again, $\pi^{(k)} \propto \exp(-2f/\eta_k)$.

\begin{theorem}[Annealed Unadjusted Langevin \Renyi{} discretization bound]
    \label{thm:ula_discretization}
    Let $f : \mathbb{R}^d \rightarrow \mathbb{R}$ be such that $\exp(-2f/\eta_k)$ is integrable for $k = 1, \dots, M$. Assume that $\nabla f$ is $L$-Lipschitz. Let $\hat{\mu}_T$ be the law of the discretized process at time $T$ and $\mu_T$ the law of the corresponding Langevin diffusion at time $T$, where both processes start from the same distribution $\mu_0$. Let $T_{m} = \sum_{r=1}^{m} \tau_r$, and assume that each $N_k = \tau_k/h_k$ is an integer. There are universal constants $B, C > 0$ such that for all $q \geq 2$, if $h_k \leq \min((C T q^2 L^2)^{-1}, L^{-1})$ for $k = 1, \dots, M$, then
    {
\begin{align}
 \mathcal{R}_q(\hat{\mu}_T || \mu_T)\nonumber &\leq \frac{1}{q-1}\max_{k=1, \dots, M}\bigg[B T q^2 L^3 h_k^2\left(d + \mathcal{R}_q(\mu_{T_{k-1}} || \pi^{(k)}) + 1\right)\nonumber\\&+ B T q^2L^2
    h_k d\bigg] \label{eq:N-step-bound}
\end{align}
}
\end{theorem}

\begin{proof}

We follow the method given by~\citep{sinhobook}. We start with Girsanov's theorem to control the \Renyi{} divergence in a single step.
\begin{align*}
    &\exp\bigl((q-1)\, R_q(\mathrm{P}_{T} \| \mathrm{W}_{T})\bigr)= \mathbb{E}_{\mathrm{W}_{T}}\!\left[\left(\frac{d\mathrm{P}_{T}}{d\mathrm{W}_T}\right)^{\!q}\right]\\
    &= \mathbb{E}_{\mathrm{W}_T}\exp\!\left(
        q \int_0^T \frac{\langle \nabla f(X_t) - \nabla f(X_{t^-}),\, dB_t \rangle}{\sqrt{\eta(\theta_t)}}
        - \frac{q}{2} \int_0^T \frac{\|\nabla f(X_t) - \nabla f(X_{t^-})\|^2}{\eta(\theta_t)} d t
    \right).
\end{align*}
Following the standard argument, we apply Cauchy--Schwarz and the supermartingale
property of the exponential martingale. For any $\lambda > 0$,
{
\begin{align*}
&\left(\mathbb{E}_{\mathrm{W}_T}\!\left[\left(\frac{d\mathrm{P}_T}{d\mathrm{W}_T}\right)^{\!q}\right]\right)^{\!2}\\
    &\leq
    \mathbb{E}_{\mathrm{W}_T}\exp\!\left(
        2q \int_0^T \frac{\langle \nabla f(X_t) - \nabla f(X_{t^-}),\, dB_t \rangle}{\sqrt{\eta(\theta_t)}}
        - \lambda \int_0^T \frac{\|\nabla f(X_t) - \nabla f(X_{t^-})\|^2}{\eta(\theta_t)} \dd t
    \right) \\
    &\quad \times \mathbb{E}_{\mathrm{W}_T}\exp\!\left(
        \bigl(\lambda - q\bigr) \int_0^T \frac{\|\nabla f(X_t) - \nabla f(X_{t^-})\|^2}{\eta(\theta_t)} \dd t
    \right).
\end{align*}
}
Setting $\lambda = 2q^2$ makes the first factor on the right-hand side of the Cauchy--Schwarz
inequality equal to the expectation of a non-negative stochastic exponential, and since this
stochastic exponential is a supermartingale, its expectation is at most $1$. Noting that the left-hand
side is at least $1$ (by Jensen), we may drop the square and, since $\lambda - q \leq 2q^2$, obtain
{
\begin{align*}
    \mathbb{E}_{\mathrm{W}_T}\!\left[\left(\frac{d\mathrm{P}_T}{d\mathrm{W}_T}\right)^{\!q}\right]
    &\leq
    \mathbb{E}_{\mathrm{W}_T}\exp\!\left(
        2q^2 \int_0^T \frac{\|\nabla f(X_t) - \nabla f(X_{t^-})\|^2}{\eta(\theta_t)} d t
    \right)\\
    &\leq \mathbb{E}_{\mathrm{W}_T}\exp\!\left(
        2q^2 L^2 \int_0^T \frac{\| X_t - X_{t^-}\|^2}{\eta(\theta_t)} d t
    \right)\\
    &\leq \frac{1}{T}\sum_{n=0}^{N-1}\int_{t_n}^{t_{n+1}}\mathbb{E}_{\mathrm{W}_T}\exp\!\left(
        2q^2 L^2 T \frac{\| X_t - X_{t^{-}}\|^2}{\eta(\theta_{t^-})}
    \right)dt,
\end{align*}
}
where the last step uses Jensen's inequality.

It remains to control $\|X_t - X_{t_{n}}\|$ on a single discretization interval. Fix the
$n$-th interval $[t_n, t_{n+1})$ with $t_{n+1} - t_n = h_n$, and write $x = X_{t_n}$.
Under $\mathrm{W}_T$, the process satisfies
\begin{align*}
    X_t - x &= -\int_{t_n}^{t} \nabla f(X_s) \dd s + \sqrt{\eta(\theta_{t_n})}(B_t - B_{t_n}).\\
\intertext{By the triangle inequality and the $L$-Lipschitz property of $\nabla f$,}
    \|X_t - x\| &\leq \int_{t_n}^{t} \bigl(\|\nabla f(x)\| + L\|X_s - x\|\bigr) \dd s
    + \sqrt{\eta(\theta_{t_n})}\,\|B_t - B_{t_n}\|.\\
\intertext{Gr\"onwall's inequality gives, for all $t \in [t_n, t_{n+1}]$,}
    \|X_t - x\| &\leq e^{L h(t_n)}\!\left(h(t_n) \|\nabla f(x)\|
    + \sqrt{\eta(\theta_{t_n})}\sup_{s \in [t_n, t_{n+1}]} \|B_s - B_{t_n}\|\right).
\end{align*}

We now set $\lambda := 32 q^2 T L^2$ and $S_n := \sup_{s \in [t_n, t_{n+1}]} \|B_s - B_{t_n}\|$.
Thus, since $L h(t_n) \leq 1$,
{
\begin{align}
    \label{eq:one-step-bound}
    \frac{2Tq^2L^2}{\eta(\theta_{t_n})} \|X_t - X_{t_n}\|^2
    \leq \lambda \frac{h(t_n)^2}{\eta(\theta_{t_n})} \|\nabla f(X_{t_n})\|^2
    + \lambda S_n^2.
\end{align}
}

By the reflection principle, a standard one-dimensional Brownian motion $w$ satisfies, for $\lambda h(t_n) \leq \frac{1}{4}$, $\E\exp(\lambda \sup_{s \in [0, h(t_n)]} w_s^2) \leq 2(1 - 2\lambda h(t_n))^{-1/2} - 1 \leq \exp(4\lambda h(t_n))$. Applying the one-dimensional bound to each of the $d$ independent coordinates of $B_s - B_{t_n}$, and bounding the supremum of a sum by the sum of the suprema, gives
{
\begin{align}
    \label{eq:bm-mgf}
    \log\E\exp\!\left(\lambda S_n^2\right)
    \leq 4\lambda d\, h(t_n).
\end{align}
}

Since $S_n$ is independent of $X_{t_n}$ under $\mathrm{W}_T$, and since $\lambda h(t_n) \leq 32/C \leq \frac{1}{4}$ for $C \geq 128$ by $h(t_n) \leq (C T q^2 L^2)^{-1}$, combining~\eqref{eq:one-step-bound} and~\eqref{eq:bm-mgf} gives
{
\begin{align*}
&\mathbb{E}_{\mathrm{W}_T}\left[\mathbb{E}_{\mathrm{W}_T}\left[\exp\left(\frac{2Tq^2L^2}{\eta(\theta_{t_n})} \|X_t - X_{t_n}\|^2 \right)~|~ X_{t_n} = x\right]\right]\\
&
    \leq \mathbb{E}_{\mathrm{W}_T}\left[\exp\left(\lambda\frac{h(t_n)^2}{\eta(\theta_{t_n})} \|\nabla f(X_{t_n})\|^2\right)\right]\cdot\mathbb{E}_{\mathrm{W}_T}\left[\exp\left(\lambda S_n^2\right)\right]\\
    &\leq \mathbb{E}_{\mathrm{W}_T}\left[\exp\left(\lambda\frac{h(t_n)^2}{\eta(\theta_{t_n})} \|\nabla f(X_{t_n})\|^2\right)\right]\cdot \exp\left(4\lambda d\, h(t_n)\right).
\end{align*}
}

So, with $Y_n := \exp(\lambda h(t_n)^2 \|\nabla f(X_{t_n})\|^2/\eta(\theta_{t_n}))$,
{
\begin{align}
    \mathcal{R}_q(\mathrm{P}_T \| \mathrm{W}_T) &\leq \frac{1}{q-1}\max_{n=0, \dots, N-1}\bigg[\log\mathbb{E}_{\mathrm{W}_T}\left[Y_n\right] \nonumber\\ &+ \log\mathbb{E}_{\mathrm{W}_T}\left[\exp\left(\lambda S_n^2\right)\right]\bigg] \label{eq:ula-path-max}
\end{align}
}

We now proceed to bound the first term inside the maximum in \eqref{eq:ula-path-max}. Let $\pi_{t_n}$ denote the stationary distribution for the $n$-th segment, i.e. $t \in [t_n, t_{n+1})$, and let $k$ be the phase that contains $[t_n, t_{n+1})$, so that $\pi_{t_n} = \pi^{(k)}$ and $\eta(\theta_{t_n}) = \eta_k$. From \citep[Lemma 6.2.7]{sinhobook}, applied to the potential $2f/\eta_k$, for any $u \geq 0$:
{
\begin{align*}
\pi^{(k)}\{ \lVert \nabla f(X)\rVert \geq C' \sqrt{\eta_k L}(\sqrt{d} + u)\} \leq \exp(-u^2).
\end{align*}
}

Then by H\"older's inequality (see \citep[Lemma 6.2.8]{sinhobook}), with $R_n := \mathcal{R}_q(\mu_{t_n} || \pi^{(k)})$, every measurable set $A \subseteq \mathbb{R}^d$ satisfies
{
\begin{align*}
    \mu_{t_n}(A) \leq \exp\left(\left(1 - \frac{1}{q}\right)R_n\right)\pi^{(k)}(A)^{1 - 1/q}.
\end{align*}
}
Thus, for any $0 < \delta < \frac{1}{2}$, taking $A$ to be the set where $\lVert \nabla f\rVert \geq C'\sqrt{\eta_k L}(\sqrt{d} + u)$ with $u^2 = R_n + 2\log(1/\delta)$ shows that with at least $1 - \delta$ probability under $\mu_{t_n}$,
{
\begin{align}
    \label{eq:grad-high-prob}
    \lVert \nabla f (X_{t_n})\rVert^2 \leq 2C'^2\eta_k L\left(d + R_n + 2\log(1/\delta)\right).
\end{align}
}
Thus we can apply the following lemma.
{
\begin{lemma}
    \label{lem:moment-from-tail}
    Let $Y \geq 0$ be a random variable, $v > 0$, and $0 \leq \xi \leq \frac{1}{2}$. If $Y^2 \leq v\delta^{-\xi}$ with probability at least $1 - \delta$ for every $0 < \delta < \frac{1}{2}$, then $\mathbb{E}[Y] \leq e^{\xi}\sqrt{v}$.
\end{lemma}
\begin{proof}[Proof of Lemma~\ref{lem:moment-from-tail}]
The assumption gives $\mathbb{P}(Y > \sqrt{v}\,\delta^{-\xi/2}) \leq \delta$ for every $0 < \delta < \frac{1}{2}$. If $\xi = 0$, then $Y \leq \sqrt{v}$ almost surely, so $\mathbb{E}[Y] \leq \sqrt{v}$. If $\xi > 0$, let $t_0 = 2^{\xi/2}\sqrt{v}$. For $t > t_0$, the choice $\delta = (\sqrt{v}/t)^{2/\xi} < \frac{1}{2}$ yields $\mathbb{P}(Y > t) \leq (\sqrt{v}/t)^{2/\xi}$, and integrating the tail gives
\begin{align*}
    \mathbb{E}[Y] \leq t_0 + \int_{t_0}^{\infty} \left(\frac{\sqrt{v}}{t}\right)^{2/\xi} \dd t = 2^{\xi/2}\left(1 + \frac{\xi}{2(2 - \xi)}\right)\sqrt{v} \leq e^{\xi}\sqrt{v}.
\end{align*}
\end{proof}
}
By~\eqref{eq:grad-high-prob}, for every $0 < \delta < \frac{1}{2}$, with probability at least $1 - \delta$,
{
\begin{align*}
Y_n^2 \leq v\delta^{-\xi}, \qquad v := \exp\left(4C'^2\lambda L h(t_n)^2(d + R_n)\right), \qquad \xi := 8C'^2\lambda L h(t_n)^2.
\end{align*}
}

If $h(t_n) \leq \min((C T q^2 L^2)^{-1}, L^{-1})$ with $C \geq 512 C'^2$, then $\xi \leq \frac{1}{2}$, so the conditions of Lemma~\ref{lem:moment-from-tail} are satisfied and
{
\begin{align*}
&\log\mathbb{E}_{\mathrm{W}_T}\left[Y_n\right] \leq \xi + \frac{1}{2}\log v \leq 8C'^2\lambda L h(t_n)^2\left(d + R_n + 1\right).
\end{align*}
}

From data processing we have that
\begin{align*}
    \mathcal{R}_q(\hat{\mu}_T || \mu_T) \leq \mathcal{R}_q(P_T || W_T).
\end{align*}

Lastly, note that for $\sum_{r=1}^{k-1}N_r \leq  n < \sum_{r=1}^{k}N_r,$ all $\pi_{t_n}$ are the same, i.e. the chain remains time homogeneous. Thus we can take the max over $k = 1, \dots, M$ and, by data processing, upper bound $R_n$ by $\mathcal{R}_q(\mu_{T_{k-1}} || \pi^{(k)})$, where $T_{k-1} = \sum_{r=1}^{k-1}\tau_r$. Note that $\pi^{(k)}$ denotes the stationary distribution for Langevin dynamics with learning rate $\eta_k$. Thus the result follows by inserting the bound on $\log\mathbb{E}_{\mathrm{W}_T}[Y_n]$ and~\eqref{eq:bm-mgf} into the bound \eqref{eq:ula-path-max} on $\mathcal{R}_q(\mathrm{P}_T \| \mathrm{W}_T)$, with $\lambda = 32 q^2 T L^2$ and $h(t_n) = h_k$.

\end{proof}

\subsection{Convergence of Annealed ULA}

{
The following theorem combines diffusion convergence and discretization over an annealing schedule with independently chosen noise levels and step sizes. Setting the step size equal to the noise level in a phase gives the SGD update in that phase.

\begin{theorem}[Convergence of annealed ULA]
\label{thm:annealed-ula-convergence}
\label{lem:ula-one-phase}
Let $f:\mathcal{X}\to\mathbb{R}$, where $\mathcal{X}\subseteq\mathbb{R}^d$, have an $L$-Lipschitz gradient. Consider $M$ ULA phases \eqref{eq:annealed-sgd} with positive noise levels $\eta_r$ and step sizes $h_r$, chosen independently. Suppose that $\pi^{(r)}\propto e^{-2f/\eta_r}$ is integrable and satisfies a $(\delta_r,\epsilon_r)$-WPI in the normalization \eqref{eq:wpi-convention}. Assume that \cref{lem:renyi_convergence_ld_wpi,thm:ula_discretization} apply to the diffusion and ULA kernels on $\mathcal{X}$ in every phase.

Let $q\geq2$, $\vartheta_r\in(0,1]$, $K_r\geq1$ and $\rho_r\in[0,1)$ for $1\leq r\leq M$. Suppose that
\begin{equation}
\epsilon_r K_r^{2q-1}\leq\frac{\delta_r\vartheta_r}{8}.
\label{eq:annealed-ula-defect}
\end{equation}
With $B,C$ the constants of \cref{thm:ula_discretization}, choose
\begin{align}
\tau_r&:=\frac{2q-1}{\delta_r}\left(\log K_r+1+\log\frac{8}{\vartheta_r}\right),\nonumber\\
h_r&\leq\min\left\{\frac1L,\frac{1}{4C\tau_r q^2L^2},
\frac{\vartheta_r}{4B\tau_r q^2L^2(2d+\log K_r+1)}\right\},\qquad \tau_r/h_r\in\mathbb{N}.
\label{eq:ula-phase-steps}
\end{align}
Write $Q_r$ for the ULA transition kernel with duration $\tau_r$, noise level $\eta_r$ and step size $h_r$.

Suppose that the initial law has a decomposition $\mu_0=\alpha_0\nu_0+(1-\alpha_0)\bar\nu_0$, where $\alpha_0\geq1-p_0$ and $p_0\in[0,1)$. Starting from $\nu_0$, define recursively
\begin{equation*}
\sigma_r:=\frac{\nu_{r-1}(\cdot\cap G_r)}{\nu_{r-1}(G_r)},
\qquad \nu_r:=\sigma_r Q_r,
\end{equation*}
where the measurable sets $G_r$ satisfy
\begin{equation}
\nu_{r-1}(G_r)\geq1-\rho_r,
\qquad \left\|\frac{d\sigma_r}{d\pi^{(r)}}\right\|_\infty\leq K_r.
\label{eq:annealed-ula-density}
\end{equation}

Then the law $\mu_r:=\mu_0Q_1\cdots Q_r$ of the actual ULA iterates satisfies, for every $1\leq r\leq M$,
\begin{align}
\mu_r&=\alpha_r\nu_r+(1-\alpha_r)\bar\nu_r,\nonumber\\
\alpha_r&\geq1-p_0-\sum_{j=1}^r\rho_j,
\qquad \mathcal R_q(\nu_r\|\pi^{(r)})\leq\vartheta_r,
\label{eq:ula-phase-conclusion}
\end{align}
for probability measures $\bar\nu_r$. In particular,
\begin{equation}
\textup{TV}(\mu_M,\pi^{(M)})
\leq p_0+\sum_{r=1}^M\rho_r+\sqrt{\vartheta_M/2}.
\label{eq:annealed-ula-tv}
\end{equation}
The sets $G_r$ restrict measures only in the analysis; the algorithm performs the original ULA updates and uses $\sum_{r=1}^M\tau_r/h_r$ gradient queries.
\end{theorem}
\begin{proof}
Fix a phase $r$, and write $P^{(r)}_{\tau_r}$ for the corresponding Langevin transition kernel. The density bound \eqref{eq:annealed-ula-density} implies $\mathcal R_s(\sigma_r\|\pi^{(r)})\leq\log K_r$ for every $s>1$. At order $2q-1$, the transient time in \cref{lem:renyi_convergence_ld_wpi} is at most $(2q-1)(\log K_r-1)_+/\delta_r$.

Equations \eqref{eq:renyi-wpi-convergence}, \eqref{eq:annealed-ula-defect} and \eqref{eq:ula-phase-steps} therefore give
\begin{equation}
\mathcal R_{2q-1}(\sigma_rP^{(r)}_{\tau_r}\|\pi^{(r)})
\leq(\vartheta_r/8)^2+\frac{2\epsilon_rK_r^{2q-1}}{\delta_r}
\leq\frac{\vartheta_r}{2}.
\label{eq:ula-phase-continuous}
\end{equation}

Apply \cref{thm:ula_discretization} to this phase at order $2q$. Since $h_r\leq1/L$, \eqref{eq:N-step-bound} and \eqref{eq:ula-phase-steps} yield
\begin{align}
\mathcal R_{2q}(\nu_r\|\sigma_rP^{(r)}_{\tau_r})
&\leq\frac{4B\tau_rq^2}{2q-1}
\left[L^3h_r^2(d+\log K_r+1)+L^2h_rd\right]\nonumber\\
&\leq\frac{4B\tau_rq^2L^2h_r}{2q-1}(2d+\log K_r+1)
\leq\frac{\vartheta_r}{2q-1}.
\label{eq:ula-phase-discrete}
\end{align}

The weak triangle inequality for \Renyi{} divergence combines \eqref{eq:ula-phase-continuous} and \eqref{eq:ula-phase-discrete}:
\begin{align*}
\mathcal R_q(\nu_r\|\pi^{(r)})
&\leq\frac{q-1/2}{q-1}\mathcal R_{2q}(\nu_r\|\sigma_rP^{(r)}_{\tau_r})
+\mathcal R_{2q-1}(\sigma_rP^{(r)}_{\tau_r}\|\pi^{(r)})\nonumber\\
&\leq\frac{\vartheta_r}{2(q-1)}+\frac{\vartheta_r}{2}
\leq\vartheta_r.
\end{align*}

We now propagate the initial decomposition. If $\mu_{r-1}\geq\alpha_{r-1}\nu_{r-1}$, then restriction to $G_r$ and application of the Markov kernel $Q_r$ give
\begin{equation}
\mu_r\geq\alpha_{r-1}\nu_{r-1}(G_r)\sigma_rQ_r
=\alpha_r\nu_r,
\qquad \alpha_r:=\alpha_{r-1}\nu_{r-1}(G_r).
\label{eq:annealed-ula-mass}
\end{equation}
The remaining measure is nonnegative, so normalization gives $\bar\nu_r$. Induction yields $\alpha_r\geq\alpha_0\prod_{j=1}^r(1-\rho_j)\geq1-p_0-\sum_{j=1}^r\rho_j$, proving \eqref{eq:ula-phase-conclusion}.

Finally, Pinsker's inequality and monotonicity of \Renyi{} divergence give $\textup{TV}(\nu_M,\pi^{(M)})\leq\sqrt{\vartheta_M/2}$. Convexity of total variation in its first argument gives \eqref{eq:annealed-ula-tv}. Each ULA step uses one gradient query, which gives the stated query count.
\end{proof}
}

\subsection{Sampling for Block-Separable Functions}

{
For the block analysis, let $g:\mathbb R^\ell\to\mathbb R$ have an $L$-Lipschitz gradient, a minimizer $a$ with $g(a)=0$ and $\|a\|\leq B_0$, and $g(y)\geq m\|y-a\|^2/2$, where $\ell\leq b$. Write $\pi_\beta\propto e^{-\beta g}$ and use $A=(2L/m)^{b/2}$, $K_0=Ae^{LB_0^2}$ and $U=1+\log K_0$, as in \cref{lem:ula-block-estimates}.

The restriction in the following lemma is used only in the analysis. ULA itself does not test the objective against a cutoff or discard iterates.

\begin{lemma}[Density bounds relative to the Gibbs distribution of the next phase]
\label{lem:ula-warm-restriction}
With these block constants, suppose $\beta'\geq\beta>0$, $\mathcal R_2(\nu\|\pi_\beta)\leq U$, and $\rho\in(0,1/32]$. Define
\begin{align}
H&:=\frac4\beta\left(\frac U2+\frac{\log A}{2}+\log\frac2\rho\right),\nonumber\\
G&:=\left\{y:\frac{d\nu}{d\pi_\beta}(y)\leq\frac{2e^U}{\rho},\quad g(y)\leq H\right\}.
\label{eq:ula-kept-set}
\end{align}
If $(\beta'-\beta)H\leq1$, the normalized restriction $\sigma=\nu(\cdot\cap G)/\nu(G)$ satisfies
\begin{align}
\nu(G^c)&\leq\rho,
\label{eq:ula-kept-mass}\\
\left\|\frac{d\sigma}{d\pi_{\beta'}}\right\|_\infty&\leq K:=\frac{4e^{U+1}}\rho,\qquad
\mathcal R_2(\sigma\|\pi_{\beta'})\leq U+2.
\label{eq:ula-warm-restriction}
\end{align}
\end{lemma}
\begin{proof}
Apply \cref{lem:truncation} with order $2$ and discarded mass $\rho/2$. It bounds the mass above the density cutoff in \eqref{eq:ula-kept-set} by $\rho/2$. Cauchy--Schwarz and \eqref{eq:ula-energy-tail} give the other half:
\begin{equation*}
\nu(g>H)\leq e^{U/2}\pi_\beta(g>H)^{1/2}
\leq e^{U/2}A^{1/2}e^{-\beta H/4}=\rho/2.
\end{equation*}
Thus \eqref{eq:ula-kept-mass} holds. Normalization bounds the density of $\sigma$ relative to $\pi_\beta$ by $4e^U/\rho$ and its order-two divergence by $U+2\log(1/(1-\rho))\leq U+1$.

Since $g\geq0$, we have $Z_{\beta'}\leq Z_\beta$. On $G$,
\begin{equation}
\frac{d\pi_\beta}{d\pi_{\beta'}}
=\frac{Z_{\beta'}}{Z_\beta}e^{(\beta'-\beta)g}\leq e.
\label{eq:ula-local-overlap}
\end{equation}
Changing reference measures using \eqref{eq:ula-local-overlap} increases the logarithmic density bound and the order-two divergence by at most $1$, proving \eqref{eq:ula-warm-restriction}.
\end{proof}

The next lemma shows that, for an arithmetic schedule of learning rates, consecutive inverse temperatures satisfy the condition of \cref{lem:ula-warm-restriction} in every phase.

\begin{lemma}[Number of phases for arithmetically decreasing learning rates]
\label{lem:ula-arithmetic-schedule}
For $\beta_*\geq2$ and $\zeta\in(0,1/4]$, define
\begin{align}
M&:=\left\lceil1024(\beta_*+1)^2(d/\zeta)^2(U+\log A+1)^2\right\rceil,\nonumber\\
\rho&:=\frac{\zeta}{8dM},\qquad
\eta_r:=2-(r-1)\frac{2-2/\beta_*}{M-1},\qquad \beta_r:=2/\eta_r.
\label{eq:ula-phase-schedule}
\end{align}
Set $\beta_0=\beta_1$. With $H_r$ given by \eqref{eq:ula-kept-set} at $\beta=\beta_{r-1}$, every phase satisfies $(\beta_r-\beta_{r-1})H_r\leq1$.
\end{lemma}
\begin{proof}
The assertion is immediate for $r=1$. Put $D=\beta_*+1$, $a=d/\zeta$, and $J=U+\log A+1$. Then
\begin{align}
512D^2a^2J^2&\leq M-1,\qquad M\leq2048D^2a^2J^2,\nonumber\\
\log(16aM)&\leq11+2D+3a+2J\leq8DaJ.
\label{eq:ula-schedule-log}
\end{align}
Here $D\geq3$, $a\geq4$, and $J\geq1$. For $r\geq2$, the relative inverse-temperature increment is at most $(\beta_*-1)/(M-1)$. Consequently, \eqref{eq:ula-schedule-log} gives
\begin{align*}
(\beta_r-\beta_{r-1})H_r
&\leq\frac{4D}{M-1}\left(\frac J2+\log(16aM)\right)\nonumber\\
&\leq\frac{34}{512aJ}<1.
\end{align*}
\end{proof}

We can now apply these estimates to ULA without repeating the diffusion convergence or discretization arguments.

\begin{theorem}[Annealed Unadjusted Langevin Algorithm for Block-Separable Functions]
\label{thm:ula-sampling-formal}
Let $f : \mathbb{R}^d \to \mathbb{R}$ be block-separable (\cref{def:separable-func}), $f(x) = g(Rx)$ with an orthogonal matrix $R$ and $g(y) = \sum_{i=1}^{k} g_i(\hat{y}_i)$, where $k \leq d$ and the blocks $\hat{y}_i \in \mathbb{R}^{d_i}$ have dimensions $d_i \leq b$.

Suppose that each $g_i \in C^2(\mathbb{R}^{d_i})$ has an $L$-Lipschitz gradient and a unique global minimizer $a_i$ with $\|a_i\| \leq B_0$, and that
\[
    g_i(y) - g_i(a_i) \geq \frac{m}{2}\|y - a_i\|^2
\]
for all $y \in \mathbb{R}^{d_i}$, where $b$, $m$, $L$ and $B_0$ do not depend on $d$.

Let $\beta_* \geq 2$, $\zeta \in (0, 1/4]$ and $\pi_{\beta_*}(x) := Z^{-1}\exp(-\beta_* f(x))$, where $Z := \int_{\mathbb{R}^d}\exp(-\beta_* f(x))\,dx$.

Then there exists an unadjusted Langevin algorithm \eqref{eq:annealed-sgd} started from $\mathcal{N}(0, (2L)^{-1}I_d)$, with learning rates $\{\eta_r\}_{r=1}^{M}$ that decrease arithmetically from $\eta_1 = 2$ to $\eta_M = 2/\beta_*$ and sufficiently small step sizes $\{h_r\}_{r=1}^{M}$ chosen independently of the learning rates.

The algorithm uses a number of gradient queries polynomial in $d$, $\beta_*$ and $1/\zeta$. Its output distribution $\mu$ satisfies
\[
    \textup{TV}(\mu, \pi_{\beta_*}) \leq \zeta.
\]
\end{theorem}
\begin{proof}
Normalize each block so that $g_i(a_i)=0$, and take the constants in \eqref{eq:ula-block-constants}. The orthogonal invariance of the Gaussian initialization and noise makes the rotated ULA iterates independent across blocks. Write $\lambda_{i,r}$ for the law of block $i$ after phase $r$ and $\pi_{i,\beta}$ for its Gibbs distribution.

Choose the arithmetic schedule and $\rho$ from \eqref{eq:ula-phase-schedule}, and set
\begin{equation}
\vartheta:=\frac{\zeta^2}{64d^2},\qquad K:=\frac{4e^{U+1}}\rho,\qquad
s:=\frac{\vartheta}{8K^3}.
\label{eq:ula-phase-budget}
\end{equation}
Use the WPI constant $\delta$ in \eqref{eq:ula-wpi-parameters} for this $s$. For every phase, choose its duration and discretization step by
\begin{align}
\tau&:=\frac3\delta\left(\log K+1+\log\frac8\vartheta\right),
\label{eq:ula-phase-duration}\\
\bar h&:=\min\left\{\frac1L,\frac{1}{16C\tau L^2},
\frac{\vartheta}{16B\tau L^2(2b+\log K+1)}\right\},\nonumber\\
h_r&:=\frac{\tau}{\lceil\tau/\bar h\rceil}.
\label{eq:ula-implementable-step}
\end{align}
These choices use only the dimension, accuracy, final temperature and the stated regularity bounds, not the block decomposition or the minimizers.

Apply \cref{thm:annealed-ula-convergence} to each block with $q=2$, $p_0=0$, $\vartheta_r=\vartheta$, $K_r=K$, $\rho_r=\rho$, $\delta_r=\delta$ and $\epsilon_r=\delta s$. The choices \eqref{eq:ula-phase-duration}--\eqref{eq:ula-implementable-step} satisfy \eqref{eq:ula-phase-steps}, and \eqref{eq:ula-phase-budget} gives $\epsilon_rK_r^3=\delta\vartheta/8$.

The density conditions \eqref{eq:annealed-ula-density} are verified inductively. Initially, \eqref{eq:ula-initial-ratio} gives $\mathcal R_2(\nu_{i,0}\|\pi_{i,\beta_0})\leq\log K_0\leq U$. Given the same bound with $U$ after a phase, \cref{lem:ula-arithmetic-schedule,lem:ula-warm-restriction} provide the set for the next phase, with discarded mass at most $\rho$ and density at most $K$. The \Renyi{} divergence in the phase conclusion \eqref{eq:ula-phase-conclusion} is at most $\vartheta\leq U$, so this construction continues through all phases.

The convergence theorem for annealed ULA therefore gives
\begin{align}
\lambda_{i,r}&=\alpha_{i,r}\nu_{i,r}+(1-\alpha_{i,r})\bar\nu_{i,r},\nonumber\\
\alpha_{i,r}&\geq1-r\rho,\qquad
\mathcal R_2(\nu_{i,r}\|\pi_{i,\beta_r})\leq\vartheta\leq U.
\label{eq:ula-induction}
\end{align}

At the last phase, Pinsker's inequality and \eqref{eq:ula-induction} yield
\begin{align*}
\textup{TV}(\lambda_{i,M},\pi_{i,\beta_*})
&\leq M\rho+\sqrt{\vartheta/2}<\frac{\zeta}{2d},\nonumber\\
\textup{TV}(\mu,\pi_{\beta_*})
&\leq\sum_{i=1}^k\textup{TV}(\lambda_{i,M},\pi_{i,\beta_*})
\leq\frac\zeta2\leq\zeta.
\end{align*}
Here we used independence across blocks and invariance of total variation under rotation.

Finally, $M$, $\rho^{-1}$, $K$ and $s^{-1}$ are polynomial in $d$, $\beta_*$ and $\zeta^{-1}$. The same holds for $\delta^{-1}$ by \cref{lem:ula-block-estimates}, hence for $\tau$ and $\bar h^{-1}$ by \eqref{eq:ula-phase-duration}--\eqref{eq:ula-implementable-step}. Each iteration uses one gradient query, and the total number is $M\lceil\tau/\bar h\rceil$, which is polynomial in those parameters.
\end{proof}
}

\subsection{Global Optimization}

{
\cref{thm:ula-sampling-formal} implies a global optimization guarantee for the same algorithm on $\mathbb{R}^d$, because for $\beta_*$ of order $d/\epsilon^2$ the Gibbs distribution $\pi_{\beta_*}$ assigns most of its mass to the points $x$ with $f(x) - f^{\star} < \epsilon$.

\begin{theorem}[Global Optimization of Block-Separable Functions via the Unadjusted Langevin Algorithm]
\label{thm:ula-optimization-formal}
Let $f$ satisfy the hypotheses of \cref{thm:ula-sampling-formal}, let $f^{\star} := \min_{x \in \mathbb{R}^d} f(x)$, and let $\epsilon \in (0, 1]$. Then there exist $\beta_* = \Theta(d/\epsilon^2)$ and an unadjusted Langevin algorithm \eqref{eq:annealed-sgd} started from $\mathcal{N}(0, (2L)^{-1}I_d)$, with learning rates that decrease arithmetically from $\eta_1 = 2$ to $\eta_M = 2/\beta_*$ and sufficiently small step sizes chosen independently of the learning rates, that uses $\textup{poly}(d/\epsilon)$ gradient queries and outputs a point $x$ with $f(x) - f^{\star} < \epsilon$ with probability at least $2/3$.
\end{theorem}
\begin{proof}
Adding constants to the functions $g_i$ changes neither the algorithm nor the event $f(x) - f^{\star} < \epsilon$, so, as in the proof of \cref{thm:ula-sampling-formal}, we assume $g_i(a_i) = 0$ for every $i \in [k]$. Then $g_i \geq 0$ and $f^{\star} = 0$.

We first bound the mean of $f$ under the Gibbs distribution $\pi_\beta \propto e^{-\beta f}$ on $\mathbb{R}^d$. For $\beta > 0$ and $u \geq 0$, the function $u \mapsto ue^{-\beta u/2}$ is at most $2/(e\beta)$, so $ue^{-\beta u} \leq \frac{2}{e\beta}e^{-\beta u/2}$. We apply this inequality with $u = g_i(y)$ to the block Gibbs measures $\pi_{i,\beta}(y) = Z_{i,\beta}^{-1}e^{-\beta g_i(y)}$, where $Z_{i,\beta} = \int_{\mathbb{R}^{d_i}}e^{-\beta g_i(y)}\,dy$, and use the bound $Z_{i,\beta/2}/Z_{i,\beta} \leq A$ of \eqref{eq:ula-energy-tail}, where $A = (2L/m)^{b/2}$. Since the rotation $x \mapsto Rx$ sends $\pi_\beta$ to the product of the measures $\pi_{i,\beta}$ and $k \leq d$, we obtain, for every $\beta > 0$,
\begin{align}
\mathbb{E}_{\pi_{i,\beta}}[g_i] &= \frac{1}{Z_{i,\beta}}\int_{\mathbb{R}^{d_i}}g_i(y)e^{-\beta g_i(y)}\,dy \leq \frac{2}{e\beta}\,\frac{Z_{i,\beta/2}}{Z_{i,\beta}} \leq \frac{2A}{e\beta}, \nonumber\\
\mathbb{E}_{\pi_\beta}[f - f^{\star}] &= \sum_{i=1}^{k}\mathbb{E}_{\pi_{i,\beta}}[g_i] \leq \frac{2Ak}{e\beta} \leq \frac{2Ad}{e\beta}. \label{eq:ula-gibbs-mean}
\end{align}

We choose $\beta_* := \max\{2, 24Ad/(e\epsilon^2)\}$, which is $\Theta(d/\epsilon^2)$ because $A$ does not depend on $d$. Markov's inequality, \eqref{eq:ula-gibbs-mean} with $\beta = \beta_*$ and $\epsilon \leq 1$ give
\begin{align}
\pi_{\beta_*}(f - f^{\star} \geq \epsilon) &\leq \frac{\mathbb{E}_{\pi_{\beta_*}}[f - f^{\star}]}{\epsilon} \leq \frac{2Ad}{e\beta_*\epsilon} \leq \frac{\epsilon}{12} \leq \frac{1}{12}, \nonumber\\
\pi_{\beta_*}(f - f^{\star} < \epsilon) &\geq \frac{11}{12}. \label{eq:ula-gibbs-success}
\end{align}

We now apply \cref{thm:ula-sampling-formal} with this $\beta_* \geq 2$ and $\zeta := 1/12$. The output distribution $\mu$ of the resulting algorithm satisfies $\textup{TV}(\mu, \pi_{\beta_*}) \leq 1/12$. By the definition of the total variation distance and \eqref{eq:ula-gibbs-success}, the output $X \sim \mu$ therefore satisfies
\begin{align*}
\mathbb{P}_{\mu}[f(X) - f^{\star} < \epsilon] &\geq \pi_{\beta_*}(f - f^{\star} < \epsilon) - \textup{TV}(\mu, \pi_{\beta_*}) \nonumber\\
&\geq \frac{11}{12} - \frac{1}{12} = \frac{5}{6} > \frac{2}{3}. 
\end{align*}
\cref{thm:ula-sampling-formal} also bounds the number of gradient queries by a polynomial in $d$, $\beta_*$ and $1/\zeta = 12$. Since $d \leq d/\epsilon$ and $\beta_* \leq (2 + 24A/e)(d/\epsilon)^2$, this number is polynomial in $d/\epsilon$.
\end{proof}
}

\section{Analysis of Stochastic Gradient Descent}
\label{sec:annealing_related_proofs}

We analyze stochastic gradient descent with decaying and fixed learning rates.

The stationary law of discrete SGD need not equal the Gibbs law of Langevin diffusion. For a quadratic block $g(y)=\lambda\|y\|^2/2$ on $\mathbb R^\ell$, with $\lambda>0$ and $0<\eta\lambda<2$, SGD has stationary covariance $\eta[2\lambda(1-\eta\lambda/2)]^{-1}I_\ell$, whereas the Gibbs covariance is $\eta(2\lambda)^{-1}I_\ell$. A spectral-transfer condition must therefore concern the actual SGD kernel and account for its stationary bias.

\subsection{Sampling with Annealed Stochastic Gradient Descent}

{
The following discrete analogue of the diffusion convergence bound uses the stationary reverse kernel. It does not require the SGD kernel to be reversible.

\begin{lemma}[R\'{e}nyi convergence from a discrete weak Poincar\'e inequality]
\label{lem:sgd-renyi-contraction}
Let $P$ be a Markov kernel with invariant probability measure $\rho$, and let $P^*$ be its adjoint in $L^2(\rho)$. Suppose that, for every bounded function $v$,
\begin{equation}
\frac{\|v\|_{L^2(\rho)}^2-\|P^*v\|_{L^2(\rho)}^2}{\eta}
\geq\frac{\delta}{\kappa}\operatorname{Var}_\rho(v)
-\kappa\delta s\|v\|_\infty^2,
\label{eq:sgd-discrete-wpi}
\end{equation}
where $\eta,\delta>0$, $\kappa\geq1$, $s\geq0$ and $a:=\eta\delta/\kappa\leq1$. Let $r\geq2$ and let the initial density satisfy $\|d\mu_0/d\rho\|_\infty\leq K<\infty$. Write $u_n:=\mathcal R_r(\mu_0P^n\|\rho)$ and $b:=\kappa^2sK^r$. If $b\leq1/4$, then
\begin{equation}
u_n\leq\exp\left[-\frac{a(n-n_0)}{2(r-1)}\right]+2b
\quad(n\geq n_0),
\qquad n_0:=\left\lceil\frac{4(r-1)(u_0-1)_+}{a}\right\rceil.
\label{eq:sgd-discrete-convergence}
\end{equation}
\end{lemma}
\begin{proof}
Stationarity makes $P^*$ a Markov kernel. The densities $f_n:=d(\mu_0P^n)/d\rho$ therefore satisfy $f_{n+1}=P^*f_n$ and $0\leq f_n\leq K$. Let $F_n:=\int f_n^r\,d\rho=e^{(r-1)u_n}$. Jensen's inequality and \eqref{eq:sgd-discrete-wpi} applied to $f_n^{r/2}$ give
\begin{align*}
F_n-F_{n+1}
&\geq\|f_n^{r/2}\|_{L^2(\rho)}^2-\|P^*(f_n^{r/2})\|_{L^2(\rho)}^2\\
&\geq a\left[\operatorname{Var}_\rho(f_n^{r/2})-b\right].
\end{align*}
Monotonicity of R\'{e}nyi divergence gives $\operatorname{Var}_\rho(f_n^{r/2})\geq F_n(1-e^{-u_n})$. Since $F_n\geq1$ and $\log(1+x)\leq x$, it follows that
\[
u_{n+1}\leq u_n-\frac{a}{r-1}(1-e^{-u_n})+\frac{ab}{r-1}.
\]

When $u_n>1$, the decrement is at least $a/[4(r-1)]$, so $u_{n_0}\leq1$. Data processing makes $u_n$ non-increasing. For $0\leq u_n\leq1$, the bound $1-e^{-u_n}\geq u_n/2$ yields
\[
u_{n+1}\leq\left(1-\frac{a}{2(r-1)}\right)u_n+\frac{ab}{r-1}.
\]
Iterating this scalar recursion proves \eqref{eq:sgd-discrete-convergence}.
\end{proof}
}

{
For a block potential $g_i:\mathbb R^{d_i}\to\mathbb R$, let $P_{i,\eta}$ be the actual SGD kernel
\[
Y'=Y-\eta\nabla g_i(Y)+\eta Z,
\qquad Z\sim\mathcal N(0,I_{d_i}),
\]
and let $\pi_i^\eta\propto e^{-2g_i/\eta}$ be the corresponding Langevin Gibbs law. Its Dirichlet form is $\mathcal E_{\pi_i^\eta}(v,v)=\frac\eta2\int\|\nabla v\|^2\,d\pi_i^\eta$. We call $(\delta,\epsilon)$ a continuous WPI pair if $\mathcal E_{\pi_i^\eta}(v,v)\geq\delta\operatorname{Var}_{\pi_i^\eta}(v)-\epsilon\|v\|_\infty^2$. The continuous pairs used below are supplied by \cref{lem:ula-block-estimates,cor:wpi-exp-defect}.

\begin{assumption}[Spectral transfer from Langevin diffusion to SGD]
\label{ass:sgd-spectral-transfer}
\label{ass:sgd-transfer}
\label{ass:block-sgd-initialization}
Fix $q\geq2$. There are constants $\eta_{\mathrm{tr}}>0$, $C_\kappa\geq1$ and $\nu_\kappa\geq0$, independent of $d$, and $1\leq\kappa_d\leq C_\kappa[\log(ed)]^{\nu_\kappa}$ such that, for every block $g_i$ and $0<\eta\leq\eta_{\mathrm{tr}}$, its SGD kernel $P_{i,\eta}$ has an invariant law $\rho_i^\eta$. Let $P_{i,\eta}^*$ be its adjoint in $L^2(\rho_i^\eta)$ and $\pi_i^\eta\propto e^{-2g_i/\eta}$ its Gibbs law. For every bounded $v$ and every continuous WPI pair $(\delta,\epsilon)$ used in the analysis, we require
\begin{align*}
\frac{\|v\|_{L^2(\rho_i^\eta)}^2-\|P_{i,\eta}^*v\|_{L^2(\rho_i^\eta)}^2}{\eta}
&\geq\frac\delta{\kappa_d}\operatorname{Var}_{\rho_i^\eta}(v)
-\kappa_d\epsilon\|v\|_\infty^2,\\
\max\{\mathcal R_t(\rho_i^\eta\|\pi_i^\eta),
\mathcal R_t(\pi_i^\eta\|\rho_i^\eta)\}
&\leq1,\qquad 2\leq t\leq8q.
\end{align*}
Here $\mathcal R_t$ is the order-$t$ R\'{e}nyi divergence. The constants may depend on $q$ and the uniform block regularity bounds.
\end{assumption}

The spectral inequality concerns one-step dissipation of the discrete kernel, not the distribution produced by an annealing run. Constant-factor transfer corresponds to bounded $\kappa_d$; polylogarithmic distortion changes the phase lengths but not the polynomial-versus-exponential runtime conclusion. The constant comparisons control initialization, changes between consecutive invariant laws, and the density moments needed to prove the final stationary bias.

\begin{lemma}[One-step Gibbs discrepancy for SGD]
\label{lem:sgd-one-step-gibbs}
Let $g_i$ have an $L$-Lipschitz gradient and an integrable Gibbs law $\pi_i^\eta\propto e^{-2g_i/\eta}$ on $\mathbb R^{d_i}$. Fix $t\geq2$, and let $B,C$ be the constants of \cref{thm:ula_discretization}. If
\[
0<\eta\leq\min\{L^{-1},(\sqrt C\,tL)^{-1}\},
\]
then one actual SGD update satisfies
\begin{equation}
\mathcal R_t(\pi_i^\eta P_{i,\eta}\|\pi_i^\eta)
\leq\frac{3Bt^2}{t-1}\,d_iL^2\eta^2.
\label{eq:sgd-one-step-gibbs}
\end{equation}
\end{lemma}
\begin{proof}
Start the Langevin diffusion from $\pi_i^\eta$, so its law remains $\pi_i^\eta$. Apply \cref{thm:ula_discretization} with one phase, duration $T=\eta$, step size $h=\eta$, and order $t$. Its step-size condition follows from $C\eta^2t^2L^2\leq1$. Equation \eqref{eq:N-step-bound} gives
\[
\mathcal R_t(\pi_i^\eta P_{i,\eta}\|\pi_i^\eta)
\leq\frac{B\eta t^2}{t-1}
\left[L^3\eta^2(d_i+1)+L^2\eta d_i\right].
\]
Using $L\eta\leq1$ and $d_i+1\leq2d_i$ proves \eqref{eq:sgd-one-step-gibbs}.
\end{proof}

The next result derives the stationary bias from the SGD update and the Gibbs weak Poincar\'e inequality. Constant comparisons with the Gibbs law control the density moments used in the proof.

\begin{lemma}[Stationary bias from constant Gibbs comparisons]
\label{lem:sgd-stationary-bias}
Let $g_i\in C^2(\mathbb R^{d_i})$, where $1\leq d_i\leq b$, have an $L$-Lipschitz gradient and a minimizer $a_i$ with $g_i(a_i)=0$ and $g_i(y)\geq m\|y-a_i\|^2/2$. Fix $q\geq2$. Suppose that, for $0<\eta\leq\eta_*$, the SGD kernel $P_{i,\eta}$ has an invariant law $\rho_i^\eta$ satisfying
\[
\mathcal R_s(\rho_i^\eta\|\pi_i^\eta)\leq1,
\qquad 2\leq s\leq8q,
\]
and that the Gibbs law satisfies
\begin{equation}
\mathcal E_{\pi_i^\eta}(w,w)
\geq\delta_*\operatorname{Var}_{\pi_i^\eta}(w)
-\delta_* C_*e^{-c_*/\eta}\|w\|_\infty^2
\label{eq:sgd-bias-wpi}
\end{equation}
for bounded smooth $w$, with fixed positive $\delta_*,C_*,c_*$. There are constants $\eta_{\mathrm{stat}}>0$ and $B_{\mathrm{stat}}\geq1$, depending only on $q,b,m,L,\eta_*,\delta_*,C_*,c_*$, such that
\begin{equation}
\mathcal R_t(\rho_i^\eta\|\pi_i^\eta)\leq B_{\mathrm{stat}}\eta,
\qquad 2\leq t\leq2q-1,\quad 0<\eta\leq\eta_{\mathrm{stat}}.
\label{eq:sgd-proved-stationary-bias}
\end{equation}
\end{lemma}
\begin{proof}
Write $\pi=\pi_i^\eta$ and $\rho=\rho_i^\eta$, and interpolate one Euler step by
\[
X_h=X_0-h\nabla g_i(X_0)+\sqrt\eta W_h,
\qquad X_0\sim\rho,\quad 0\leq h\leq\eta,
\]
where $W_h$ is standard Brownian motion independent of $X_0$. Let $Q_h$ be this transition kernel, $\mu_h=\rho Q_h$, and $f_h=d\mu_h/d\pi$. Stationarity gives $\mu_0=\mu_\eta=\rho$.

The proof of \cref{lem:sgd-one-step-gibbs}, with duration and step size $h\leq\eta$, gives $\mathcal R_s(\pi Q_h\|\pi)\leq C_s h^2$ for $2\leq s\leq8q$. Choose $\eta_{\mathrm{stat}}\leq\min\{\eta_*,1\}$ small enough that these bounds are at most $1/2$ and all their step-size conditions hold. The weak triangle inequality and data processing then give
\begin{align}
\mathcal R_s(\mu_h\|\pi)
&\leq\frac{s-1/2}{s-1}\mathcal R_{2s}(\rho\|\pi)
+\mathcal R_{2s-1}(\pi Q_h\|\pi)\leq2,
\qquad 2\leq s\leq4q,\nonumber\\
\int f_h^s\,d\pi&\leq e^{2(s-1)},\qquad 2\leq s\leq4q.
\label{eq:sgd-bias-density-moments}
\end{align}

The quadratic bounds and Gaussian integration give $\mathbb E_\pi\|\nabla g_i\|^8\leq C_g\eta^4$. H\"older's inequality and $\mathcal R_2(\rho\|\pi)\leq1$ therefore imply $\mathbb E_\rho\|\nabla g_i\|^4\leq C_g'\eta^2$. Consequently,
\begin{equation}
\mathbb E\|X_h-X_0\|^4
\leq8h^4 C_g'\eta^2+8\eta^2h^2d_i(d_i+2)
\leq C_\Delta\eta^2h^2.
\label{eq:sgd-bias-displacement}
\end{equation}
All constants here are uniform over $d_i\leq b$.

Fix $r\in[2,2q-1]$, and put $F(h)=\int f_h^r\,d\pi$ and $u(h)=\mathcal R_r(\mu_h\|\pi)$. The drift error of the interpolated law is
\[
e_h(x):=\nabla g_i(x)-\mathbb E[\nabla g_i(X_0)\mid X_h=x].
\]
Conditional Jensen's inequality, Cauchy--Schwarz, \eqref{eq:sgd-bias-density-moments}, and \eqref{eq:sgd-bias-displacement} give
\begin{align}
D(h)&:=\frac1{F(h)}\int f_h^r\|e_h\|^2\,d\pi\nonumber\\
&\leq\frac{L^2}{F(h)}
\left(\int f_h^{2r-1}\,d\pi\right)^{1/2}
\left(\mathbb E\|X_h-X_0\|^4\right)^{1/2}
\leq C_e\eta h,
\label{eq:sgd-bias-drift-error}
\end{align}
where $F(h)\geq1$. The Fokker--Planck equation, integration by parts, and Young's inequality yield, with $v_h=f_h^{r/2}$,
\begin{align}
u'(h)&=-\frac{\eta r}{2F(h)}\int f_h^{r-2}\|\nabla f_h\|^2\,d\pi
+\frac r{F(h)}\int f_h^{r-1}\langle\nabla f_h,e_h\rangle\,d\pi\nonumber\\
&\leq-\frac2{rF(h)}\mathcal E_\pi(v_h,v_h)+\frac r\eta D(h).
\label{eq:sgd-bias-differential}
\end{align}
Gaussian transition densities and the moment bounds justify these operations by smooth truncation.

To apply the WPI to the unbounded $v_h$, write $s_\eta=C_*e^{-c_*/\eta}$. For a nonnegative $v\in L^4(\pi)$, apply \eqref{eq:sgd-bias-wpi} to smooth approximations of $v\wedge K$. The bounds
\begin{align}
\operatorname{Var}_\pi(v)
&\leq2\operatorname{Var}_\pi(v\wedge K)
+2\|v-v\wedge K\|_{L^2(\pi)}^2\nonumber\\
&\leq\frac2{\delta_*}\mathcal E_\pi(v,v)
+2s_\eta K^2+\frac2{K^2}\int v^4\,d\pi\nonumber\\
\intertext{and the choice $K^4=(\int v^4\,d\pi)/s_\eta$ give}
\operatorname{Var}_\pi(v)
&\leq\frac2{\delta_*}\mathcal E_\pi(v,v)
+4\sqrt{s_\eta}\left(\int v^4\,d\pi\right)^{1/2}.
\label{eq:sgd-bias-extended-wpi}
\end{align}
For $v=v_h$, R\'{e}nyi monotonicity gives $\operatorname{Var}_\pi(v_h)/F(h)\geq1-e^{-u(h)}$. Equation \eqref{eq:sgd-bias-density-moments} bounds $\int v_h^4\,d\pi$ uniformly and gives $u(h)\leq2$, so $1-e^{-u(h)}\geq u(h)/4$. Combining \eqref{eq:sgd-bias-drift-error}, \eqref{eq:sgd-bias-differential}, and \eqref{eq:sgd-bias-extended-wpi} yields uniform constants $A_0,A_1>0$ for which
\begin{equation}
u'(h)\leq-\lambda u(h)+A_0e^{-c_*/(2\eta)}+A_1h,
\qquad\lambda:=\frac{\delta_*}{8q}.
\label{eq:sgd-bias-scalar}
\end{equation}
Integrating \eqref{eq:sgd-bias-scalar} over $[0,\eta]$ and using $u(\eta)=u(0)$ gives
\[
u(0)\leq\frac{A_0}{\lambda}e^{-c_*/(2\eta)}+\frac{A_1}{\lambda}\eta.
\]
Since $e^{-c_*/(2\eta)}\leq2\eta/(ec_*)$, taking $B_{\mathrm{stat}}=\max\{1,\lambda^{-1}(A_1+2A_0/(ec_*))\}$ proves \eqref{eq:sgd-proved-stationary-bias}.
\end{proof}
}

{
We first state the annealing argument for one block. The neighboring Gibbs comparison and the continuous WPI pairs will be verified explicitly in the application to block-separable functions below.

\begin{theorem}[Convergence of annealed SGD under spectral transfer]
\label{thm:sgd-spectral-annealing}
\label{thm:formal_thm_gen_nonconvex_sampling}
Let $P_r=P_{\eta_r}$ be the SGD kernels at noise levels $\eta_1>\cdots>\eta_M>0$, with Gibbs laws $\pi_r\propto e^{-2g/\eta_r}$ and invariant laws $\rho_r$. Suppose that \cref{ass:sgd-spectral-transfer} holds, that $\eta_1\leq\eta_{\mathrm{tr}}$, and that
\begin{equation}
\begin{aligned}
\mathcal R_{8q-6}(\pi_{r-1}\|\pi_r)&\leq1\quad(2\leq r\leq M),\\
\left\|\frac{d\mu_0}{d\pi_1}\right\|_\infty&\leq K_{\mathrm{init}}<\infty.
\end{aligned}
\label{eq:sgd-spectral-initialization}
\end{equation}
Let $p\in(0,1/48]$ and $\theta\in(0,1]$. Suppose that $g$ and its low-temperature Gibbs family satisfy the smoothness, quadratic-growth, and WPI hypotheses of \cref{lem:sgd-stationary-bias} on $0<\eta\leq\eta_*\leq\eta_{\mathrm{tr}}$. Let $\eta_{\mathrm{stat}},B_{\mathrm{stat}}$ be the constants proved there, and choose $\eta_M\leq\min\{\eta_{\mathrm{stat}},\theta/(2B_{\mathrm{stat}})\}$. Define
\begin{align}
c_t&:=\frac{t-1/2}{t-1},
&H&:=\max\{6,c_q\log K_{\mathrm{init}}+1\},\nonumber\\
\varrho&:=p/M,
&K&:=\frac{e^H\varrho^{-1/(q-1)}}{1-\varrho},\nonumber\\
v&:=\frac{\theta}{2c_q}.
\label{eq:sgd-spectral-phase-parameters}
\end{align}
For each phase, choose a continuous WPI pair $(\delta_r,\epsilon_r)$ satisfying
\begin{align}
\kappa_d^2\frac{\epsilon_r}{\delta_r}K^{2q}&\leq\frac v4,
\qquad a_r:=\frac{\eta_r\delta_r}{\kappa_d}\leq1,
\label{eq:sgd-spectral-defect-budget}\\
\intertext{and run}
N_r&:=\left\lceil\frac{4(2q-1)(\log K-1)_+}{a_r}\right\rceil
+\left\lceil\frac{2(2q-1)}{a_r}\log\frac2v\right\rceil
\label{eq:sgd-spectral-phase-length}
\end{align}
SGD updates. Then the actual output law satisfies
\begin{align}
\mu&=\alpha\nu+(1-\alpha)\bar\nu,\nonumber\\
\alpha&\geq1-p\geq1-3p,
\qquad \mathcal R_q(\nu\|\pi_M)\leq\theta.
\label{eq:sgd-shared-conclusion}
\end{align}
The algorithm performs $\sum_rN_r$ gradient queries.
\end{theorem}
\begin{proof}
The weak triangle inequality gives
\[
\mathcal R_q(\mu_0\|\rho_1)
\leq c_q\mathcal R_{2q}(\mu_0\|\pi_1)
+\mathcal R_{2q-1}(\pi_1\|\rho_1)
\leq c_q\log K_{\mathrm{init}}+1\leq H.
\]
For consecutive invariant laws, apply the same inequality twice, with $t=2q-1$:
\begin{align*}
\mathcal R_t(\rho_{r-1}\|\rho_r)
&\leq c_t\mathcal R_{2t}(\rho_{r-1}\|\pi_{r-1})
+c_{2t-1}\mathcal R_{4t-2}(\pi_{r-1}\|\pi_r)
+\mathcal R_{4t-3}(\pi_r\|\rho_r)\\
&\leq c_t+c_{2t-1}+1<4.
\end{align*}
All stationary comparison orders are at most $8q$, and the order of the neighboring Gibbs comparison is $4t-2=8q-6$.

We now construct the retained component phase by phase, starting with $\nu_0:=\mu_0$. Given $\mathcal R_{2q}(\nu_{r-1}\|\rho_{r-1})\leq v\leq1$, the weak triangle inequality and the bound for consecutive invariant laws give $\mathcal R_q(\nu_{r-1}\|\rho_r)\leq c_qv+4<6\leq H$. In the first phase use the initial bound just proved. Applying \cref{lem:truncation} at order $q$ with discarded mass $\varrho$ gives a normalized restriction $\sigma_r$ with density at most $K$ relative to $\rho_r$.

The spectral-transfer condition gives \eqref{eq:sgd-discrete-wpi} with $s=\epsilon_r/\delta_r$. Apply \cref{lem:sgd-renyi-contraction} at order $2q$, and let $n_{0,r}$ be its transient step count. Its defect parameter is at most $v/4$ by \eqref{eq:sgd-spectral-defect-budget}, and the phase length \eqref{eq:sgd-spectral-phase-length} gives
\[
\mathcal R_{2q}(\nu_r\|\rho_r)
\leq e^{-a_r(N_r-n_{0,r})/[2(2q-1)]}+v/2\leq v,
\qquad \nu_r:=\sigma_rP_r^{N_r}.
\]
Thus the same bound holds after every phase. The restrictions are used only in the proof. The mass-propagation argument \eqref{eq:annealed-ula-mass}, applied to the Markov kernels $P_r^{N_r}$, gives $\alpha\geq(1-\varrho)^M\geq1-p$ for the retained component of the actual SGD law.

At the last phase, the weak triangle inequality and the proved stationary-bias estimate \eqref{eq:sgd-proved-stationary-bias} give
\[
\mathcal R_q(\nu_M\|\pi_M)
\leq c_q\mathcal R_{2q}(\nu_M\|\rho_M)
+\mathcal R_{2q-1}(\rho_M\|\pi_M)
\leq c_qv+B_{\mathrm{stat}}\eta_M\leq\theta.
\]
This proves \eqref{eq:sgd-shared-conclusion}; each SGD step uses one gradient query.
\end{proof}

\begin{restatable}[Nonconvex Sampling in Total Variation]{theorem}{annealedTVSampling}
\label{thm:annealed-tvd}
Under the hypotheses and phase choices of \cref{thm:sgd-spectral-annealing}, the actual annealed SGD output satisfies
\[
\textup{TV}(\mu,\pi_M)\leq p+\sqrt{\theta/2}\leq3p+\sqrt{\theta/2}.
\]
\end{restatable}
\begin{proof}
Apply the derivation of the total-variation bound \eqref{eq:annealed-ula-tv} to the output decomposition \eqref{eq:sgd-shared-conclusion}, whose exceptional mass is at most $p$.
\end{proof}
}

\subsection{Sampling for Block-Separable Functions}

We apply the weak Poincar\'e estimates of \cref{sec:weak-poincare-proofs} to the block Gibbs measures and combine the resulting block distributions.

{
Let $f:\mathbb R^d\to\mathbb R$ have the form $f(x)=g(Rx)$, where $R$ is orthogonal and $g(y)=\sum_{i=1}^k g_i(\hat y_i)$, with blocks $\hat y_i\in\mathbb R^{d_i}$ satisfying $d_i\leq b$ and $\sum_i d_i=d$. Each $g_i\in C^2(\mathbb R^{d_i})$ has an $L$-Lipschitz gradient and a unique minimizer $a_i$, and satisfies
\[
g_i(y)-g_i(a_i)\geq\frac m2\|y-a_i\|^2
\qquad\text{for all }y\in\mathbb R^{d_i}.
\]
{
Assume that there are common constants $r_0>0$ and $0<\xi\leq L$ such that
\begin{equation}
\nabla^2g_i(y)\succeq\xi I_{d_i}
\qquad\text{whenever }\|y-a_i\|\leq r_0,
\quad 1\leq i\leq k.
\label{eq:sgd-uniform-local-convexity}
\end{equation}
A sufficient condition is the existence of a common radius on which $\|\nabla^2g_i(y)-\nabla^2g_i(a_i)\|_{\mathrm{op}}\leq m/2$ since quadratic growth gives $\nabla^2g_i(a_i)\succeq mI_{d_i}$, and then \eqref{eq:sgd-uniform-local-convexity} holds with $\xi=m/2$.
}

For $\beta>0$, define the block Gibbs laws by
\[
\pi_{i,\beta}(dy):=Z_{i,\beta}^{-1}e^{-\beta g_i(y)}\,dy,
\qquad Z_{i,\beta}:=\int_{\mathbb R^{d_i}}e^{-\beta g_i(y)}\,dy.
\]
{
Apply \cref{cor:wpi-exp-defect} to every block with the same radius $r_0$ and curvature $\xi$ from \eqref{eq:sgd-uniform-local-convexity}. In its proof, $d_i\leq b$ and the local Hessian bound satisfies $L_0\leq L$. Thus the Gaussian tail bounds and the comparison of the normalizing constants are uniform over all blocks. For example, take
\[
c:=\min\{mr_0^2/8,\xi r_0^2/32\},\qquad
C:=8(L/m)^{b/2}+16(L/\xi)^{b/2}.
\]
Once $\beta_0$ is sufficiently large in terms of $b,m,L,r_0,\xi$, the two tail contributions in that proof are bounded by $Ce^{-c\beta}$ for every $\beta\geq\beta_0$. Therefore the constants $c,C,\beta_0$ depend only on these uniform regularity constants, and every block satisfies
}
\[
\operatorname{Var}_{\pi_{i,\beta}}(v)
\leq\frac{4}{\xi\beta}\int\|\nabla v\|^2\,d\pi_{i,\beta}
+Ce^{-c\beta}\|v\|_\infty^2,
\qquad\beta\geq\beta_0,
\]
valid for every block $i$ and every bounded smooth function $v:\mathbb R^{d_i}\to\mathbb R$.

Fix $q\geq2$, $p\in(0,1/48]$ and $\theta\in(0,1]$. Let $\eta_{\mathrm{stat}},B_{\mathrm{stat}}$ be the constants of \cref{lem:sgd-stationary-bias}, applied with
\begin{align*}
\delta_*&:=\xi/4, & C_*&:=C,\\
c_*&:=2c, & \eta_*&:=\min\{\eta_{\mathrm{tr}},2/\beta_0\}.\\
\intertext{Under \cref{ass:sgd-spectral-transfer}, choose}
\eta_1&:=\min\{1,\eta_{\mathrm{tr}},1/(8qL)\},
&\beta_1&:=2/\eta_1.
\end{align*}
{
Initialize at a specified deterministic point $x_0\in\mathbb R^d$. Write $y_{i,0}:=(Rx_0)_i$ and assume that the initial block energies satisfy
\[
G_0:=\max_i[g_i(y_{i,0})-g_i(a_i)]<\infty.
\]
After one SGD update at the constant learning rate $\eta_1$, the law of each block of the iterate satisfies the density bound needed by the annealing proof.

\begin{lemma}[One-step initialization from a point]
\label{lem:sgd-point-initialization}
Under the block regularity conditions above, let
\[
X_1=x_0-\eta_1\nabla f(x_0)+\eta_1 Z,
\qquad Z\sim\mathcal N(0,I_d).
\]
The blocks of $RX_1$ are independent. If $\mu_{i,\mathrm{init}}$ denotes the law of block $i$, then
\begin{equation}
\left\|\frac{d\mu_{i,\mathrm{init}}}{d\pi_{i,\beta_1}}\right\|_\infty
\leq K_{\mathrm{init}}
:=(2m\eta_1)^{-b/2}
\exp\left(\frac{2LG_0(1+L\eta_1)^2}{m\eta_1(1-2L\eta_1)}\right).
\label{eq:sgd-point-initialization}
\end{equation}
For fixed $G_0$ and regularity constants, $K_{\mathrm{init}}$ is independent of $d$ and the target accuracy.
\end{lemma}
\begin{proof}
Subtract $g_i(a_i)$ from each block potential. This changes neither the update nor the Gibbs law. Orthogonal invariance of $Z$ gives independent Gaussian blocks with means
\begin{align*}
\bar y_i&:=y_{i,0}-\eta_1\nabla g_i(y_{i,0})\\
\intertext{and covariances $\eta_1^2 I_{d_i}$. The quadratic bounds give}
g_i(y)&\leq\frac L2\|y-a_i\|^2,\\
Z_{i,\beta_1}&\leq(\pi\eta_1/m)^{d_i/2},\\
\|\bar y_i-a_i\|^2&\leq\frac{2G_0}{m}(1+L\eta_1)^2.
\end{align*}
Consequently, for $z:=y-a_i$ and $w_i:=\bar y_i-a_i$,
\[
\frac{d\mu_{i,\mathrm{init}}}{d\pi_{i,\beta_1}}(y)
\leq(2m\eta_1)^{-d_i/2}
\exp\left(-\frac{\|z-w_i\|^2}{2\eta_1^2}+\frac L{\eta_1}\|z\|^2\right).
\]
The exponent is maximized at $z=w_i/(1-2L\eta_1)$, with maximum $L\|w_i\|^2/[\eta_1(1-2L\eta_1)]$. Since $m\leq L$, $2L\eta_1\leq1/(4q)<1$ and $d_i\leq b$, this proves \eqref{eq:sgd-point-initialization}.
\end{proof}

The algorithm performs this update in the original coordinates and does not need $R$, the block decomposition, or the minimizers. The bound on $G_0$ is a hypothesis on the given starting point. The annealing analysis starts from the Gaussian distributions of the blocks after the first SGD update. The first SGD update uses the same learning rate as the first annealing phase, so it adds one gradient query but no additional temperature level.
}

\begin{theorem}[Annealed SGD for Block-Separable Functions]
\label{thm:annealed_sgd_block_sep}
Let $f(x)=\sum_{i=1}^k g_i((Rx)_i)$ satisfy the block regularity conditions above, and suppose that \cref{ass:sgd-spectral-transfer} holds. Let $\epsilon\in(0,1]$ and choose
{
\[
\beta_M\geq\max\{C_{\mathrm{opt}}d/\epsilon^2,4/\eta_1,\beta_0,
2/\eta_{\mathrm{stat}},4B_{\mathrm{stat}}d/\theta\},
\]
}
where $C_{\mathrm{opt}}$ is the constant of \cref{lem:sufficient_beta_for_opt}.

There is an arithmetic learning-rate schedule from $\eta_1$ to $2/\beta_M$ that uses the SGD update in every phase and whose output satisfies
\begin{equation}
\mu=\alpha\nu+(1-\alpha)\bar\nu,
\qquad\alpha\geq1-p\geq1-3p,
\qquad\mathcal R_q(\nu\|\pi_{\beta_M})\leq\theta,
\label{eq:sgd-block-output}
\end{equation}
where $\pi_{\beta_M}\propto e^{-\beta_M f}$.
{
For fixed $q$ and regularity constants, the number of gradient queries, including the gradient query for the initial SGD update, is polynomial in $d$, $\beta_M$, $K_{\mathrm{init}}$, $1/p$ and $1/\theta$. If $\beta_M=\Theta(d/\epsilon^2)$ and $G_0$ is bounded independently of $d$ and $\epsilon$, this is $\textup{poly}(d/\epsilon)$. A constant or polylogarithmic distortion $\kappa_d$ preserves this polynomial bound.
}
\end{theorem}
\begin{proof}
{
Rotate to $R=I$ and subtract each $g_i(a_i)$. By \cref{lem:sgd-point-initialization}, the first SGD update from $x_0$ produces independent block laws whose densities relative to $\pi_{i,\beta_1}$ are at most $K_{\mathrm{init}}$. We apply the phasewise annealing argument to these laws and count the first update separately.
}

Let $A:=(2L/m)^{b/2}$ and $D:=\max\{1,16q,4A/e\}$. Choose
\begin{align*}
M&:=3+\lceil D(\beta_M+1)^2\rceil,\nonumber\\
\eta_r&:=\eta_1-(r-1)\frac{\eta_1-2/\beta_M}{M-1},
\qquad\beta_r:=2/\eta_r.
\end{align*}
For consecutive temperatures, put $\beta=\beta_{r-1}$, $\Delta=\beta_r-\beta$ and $t=8q-6$. The schedule gives $\Delta/\beta\leq(\beta_M-1)/(M-1)\leq1/D$, so $(t-1)\Delta/\beta\leq1/2$.

The block Gibbs mean bound \eqref{eq:ula-gibbs-mean} gives $\mathbb E_{\pi_{i,u}}g_i\leq2A/(eu)$. With $Z_i(u):=\int e^{-ug_i}$, it follows that
\begin{align*}
\mathcal R_t(\pi_{i,\beta}\|\pi_{i,\beta+\Delta})
&=\log\frac{Z_i(\beta+\Delta)}{Z_i(\beta)}
+\frac1{t-1}\log\frac{Z_i(\beta-(t-1)\Delta)}{Z_i(\beta)}\\
&\leq\frac{2A}{e(t-1)}\log\frac{\beta}{\beta-(t-1)\Delta}
\leq\frac{4A\Delta}{e\beta}\leq1.
\end{align*}
Thus the neighboring Gibbs condition \eqref{eq:sgd-spectral-initialization} is proved for every block.

{
Use the common $H,K$ from \eqref{eq:sgd-spectral-phase-parameters}, with $\varrho=p/(dM)$, and let $v:=\theta/(2c_qd)$. These parameters are defined in terms of the ambient dimension, not the unknown number of blocks. Set
}
\begin{equation*}
s_0:=\frac{v}{4\kappa_d^2K^{2q}}.
\end{equation*}
Let $\delta_0(s_0)$ be the uniform continuous WPI rate in \eqref{eq:ula-wpi-parameters}, and put
\[
\beta_c:=\max\{\beta_0,c^{-1}\log(C/s_0)\}.
\]
Choose the proved continuous WPI pairs
\begin{equation*}
(\delta_r,\epsilon_r):=
\begin{cases}
(\delta_0(s_0),\delta_0(s_0)s_0),&\beta_r<\beta_c,\\
(\xi/4,(\xi/4)Ce^{-c\beta_r}),&\beta_r\geq\beta_c.
\end{cases}
\end{equation*}
The first comes from \cref{lem:ula-block-estimates}, the second from \cref{cor:wpi-exp-defect}. Both have $\epsilon_r/\delta_r\leq s_0$, giving \eqref{eq:sgd-spectral-defect-budget}. Also $\delta_r\leq L/4$, so $a_r=\eta_r\delta_r/\kappa_d\leq1$. Run the common phase lengths \eqref{eq:sgd-spectral-phase-length} in all blocks.

{
Apply \cref{thm:sgd-spectral-annealing} to each block with $p/d$ and $\theta/d$. The constant comparisons in \cref{ass:sgd-spectral-transfer} give the bounds involving the invariant laws that are needed at initialization and between consecutive phases. At the final temperature, \cref{lem:sgd-stationary-bias} gives a stationary bias of at most $B_{\mathrm{stat}}\eta_M\leq\theta/(2d)$ by the choice of $\beta_M$. Each block output has a component of mass at least $1-p/d$ and R\'enyi divergence at most $\theta/d$ from its final Gibbs law. Let $\nu$ be the product of the retained block laws, rotated back to the original coordinates. Independence of the blocks, Bernoulli's inequality and the tensorization of R\'enyi divergence give
\[
\alpha\geq(1-p/d)^k\geq1-kp/d\geq1-p,
\qquad
\mathcal R_q(\nu\|\pi_{\beta_M})\leq k\theta/d\leq\theta,
\]
proving \eqref{eq:sgd-block-output} since $k\leq d$.
}

{
Finally, $M$ is polynomial in $\beta_M$ and $K$ is polynomial in $K_{\mathrm{init}}$ and $dM/p$ for fixed $q$. The uniform rate satisfies $\delta_0(s_0)^{-1}\leq C_\delta s_0^{-\gamma}$, with $\gamma:=16(L/m-1)$. Thus the phase lengths have polynomial cost, and the overhead from the distortion in the spectral-transfer assumption is at most a constant times $\kappa_d^{1+2\gamma}$. Since $B_{\mathrm{stat}}$ and $\eta_{\mathrm{stat}}$ depend only on the uniform regularity constants and not on $d$, the condition on the final temperature is compatible with $\beta_M=\Theta(d/\epsilon^2)$ for fixed $\theta$ and $\epsilon\leq1$. If $\kappa_d\leq C_\kappa[\log(ed)]^{\nu_\kappa}$, this overhead is polylogarithmic in $d$. Every schedule parameter is determined by $d$, the target accuracies, and the given regularity constants, without knowledge of $k$, $R$, or the decomposition. Each vector update queries $\nabla f$ once and updates all blocks.
}
Including the initial SGD update, the total number of gradient queries is $1+\sum_{r=1}^M N_r$; equivalently, the first phase contains $N_1+1$ updates.
\end{proof}
}

\subsection{Global Optimization}
\label{sec:optimization_block_sep}

In this subsection, we show that the Gibbs measure of a block-separable function concentrates on approximate global minimizers at low temperature, and we use this concentration to prove a global optimization guarantee for annealed SGD on these functions.

{
\begin{lemma}
\label{lem:sufficient_beta_for_opt}
Let $f:\mathbb R^d\to\mathbb R$ satisfy the smoothness and quadratic-growth hypotheses of \cref{thm:annealed_sgd_block_sep}, and let $f^{\star}:=\min_x f(x)$. Set
\begin{align}
A&:=\left(\frac{2L}{m}\right)^{b/2},\nonumber\\
C_{\mathrm{opt}}&:=\max\left\{2,\frac{12A}{e}\right\}.\nonumber\\
\intertext{For every $\epsilon\in(0,1]$ and $\beta\geq C_{\mathrm{opt}}d/\epsilon^2$,}
\mathbb P_{\pi_{\beta,f}}\bigl[f(X)-f^{\star}<\epsilon\bigr]&\geq\frac56.
\label{eq:sgd-gibbs-concentration}
\end{align}
\end{lemma}
\begin{proof}
Translating each block minimizer to the origin and applying the Gibbs mean estimate \eqref{eq:ula-gibbs-mean} gives
\begin{align*}
\mathbb E_{\pi_{\beta,f}}[f-f^{\star}]&\leq\frac{2Ad}{e\beta}.\\
\intertext{Markov's inequality and the choice of $C_{\mathrm{opt}}$ therefore yield}
\pi_{\beta,f}(f-f^{\star}\geq\epsilon)
&\leq\frac{2Ad}{e\beta\epsilon}
\leq\frac{\epsilon}{6}\leq\frac16,
\end{align*}
which proves \eqref{eq:sgd-gibbs-concentration}.
\end{proof}
}

\begin{theorem}[Global Optimization of Block Separable Functions Via Annealing]
\label{thm:global_opt_sep}
{
Let $f:\mathbb R^d\to\mathbb R$ satisfy the block regularity hypotheses of \cref{thm:annealed_sgd_block_sep}, with minimizer $x^{\star}$. Let annealed SGD start from a specified deterministic point $x_0$ with $G_0=\max_i[g_i((Rx_0)_i)-g_i(a_i)]$ bounded independently of $d$ and $\epsilon$, and let it include the initial SGD update of \cref{lem:sgd-point-initialization}. Suppose that \cref{ass:sgd-spectral-transfer} holds with $q=2$. Set $p=1/48$, $\theta=1/72$, and choose $\beta_M=\Theta(d/\epsilon^2)$ satisfying the temperature bound in \cref{thm:annealed_sgd_block_sep}. Then annealed SGD uses $\textup{poly}(d/\epsilon)$ gradient queries and outputs $x$ with $f(x)-f(x^{\star})<\epsilon$ with probability at least $2/3$.
}
\end{theorem}
\begin{proof}
    Theorem \ref{thm:annealed_sgd_block_sep}, with $q=2$, $p=\frac{1}{48}$ and $\theta=\frac{1}{72}$, implies that the output distribution $\mu$ of the algorithm satisfies $\mu = \alpha\nu + (1-\alpha)\bar{\nu}$ with $\alpha \geq 1 - 3p = \frac{15}{16}$ and $\mathcal{R}_2(\nu || \pi_{\beta_M}) \leq \frac{1}{72}$, where $\nu$ and $\bar{\nu}$ are probability measures on $\mathbb R^d$ and $\pi_{\beta_M} \propto \exp(-\beta_M f)$. By Pinsker's inequality and the monotonicity of \Renyi{} divergences in the order, the total variation distance satisfies $\textup{TV}(\nu, \pi_{\beta_M}) \leq \sqrt{\mathcal{R}_2(\nu || \pi_{\beta_M})/2} \leq \frac{1}{12}$.

    Since $\beta_M \geq \beta^{\star} = C_{\mathrm{opt}}d/\epsilon^2$, Lemma \ref{lem:sufficient_beta_for_opt} shows that sampling from $\pi_{\beta_M}$ obtains an $\epsilon$-approximate global optimizer with probability $\geq \frac{5}{6}$. Consequently, $\mathbb{P}_{\mu}[f(X) - f(x^{\star}) < \epsilon] \geq \alpha\,\mathbb{P}_{\nu}[f(X) - f(x^{\star}) < \epsilon] \geq \frac{15}{16}\left(\frac{5}{6} - \frac{1}{12}\right) = \frac{45}{64} > \frac{2}{3}$. Since $\beta_M = \Theta(d/\epsilon^2)$, \cref{thm:annealed_sgd_block_sep} also shows that the number of gradient queries is polynomial in $d/\epsilon$. Thus the result follows.
\end{proof}

\section{Lower Bound for SGD with Fixed Learning Rate}
\label{app:lower-bound}
The recent results of Azizian et al.~\citep{azizian2025the} give tight bounds on the expected time for fixed-step-size SGD to reach the global minimum, under a condition that bounds the norm of the gradient noise at every point $x$ by a fixed multiple of $1 + \|x\|$~\citep[Assumption~2(b)]{azizian2025the}. Similar bounds can be obtained from the analysis of Shi et al.~\citep{shi2023learning}; however, these are specialized to mixing times and therefore do not yield complete convergence lower bounds, which are instead governed by first hitting times.

The main lower bound is \cref{thm:lower-bound-main}. The noise in fixed-step-size SGD (\cref{def:fixed-sgd-intro}) is Gaussian and therefore unbounded, so we prove \cref{thm:lower-bound-main} by bounding the first hitting time of the $\epsilon$-optimal set directly. For SGD on a separable objective $f(x) = \sum_{i=1}^d \phi(x_i)$ started at $(a, \ldots, a)$, where $a$ is a spurious local minimizer of $\phi$, the first hitting time is at least the time at which a fixed coordinate leaves a neighborhood of $a$. \Cref{lem:exit-time} shows that the expected exit time from this neighborhood is at least exponential in $1/\eta$. For step sizes $\eta = \Ocal(\epsilon/d)$, as in \cref{def:fixed-sgd-intro}, the resulting lower bound on the expected first hitting time is exponential in $d$.

{
\begin{lemma}
\label{lem:exit-time}
Let $\phi:\mathbb{R}\to\mathbb{R}$ be $C^2$ with an $L$-Lipschitz derivative, and let $a$ be a local minimizer of $\phi$ with $\phi''(a) > 0$. Fix $h > 0$, let $D$ be the connected component of $\{y \in \mathbb{R} : \phi(y) < \phi(a) + h\}$ that contains $a$, and assume that $D$ is bounded and that $a$ is the only critical point of $\phi$ in the closure $\bar{D}$. Consider the SGD recursion
\begin{align*}
    y_{n+1} = y_n - \eta\,\phi'(y_n) + \eta\,\xi_n, \qquad y_0 = a,
\end{align*}
where $\xi_0, \xi_1, \ldots$ are i.i.d.\ $\mathcal{N}(0,1)$, and let $\sigma_D := \min\{n \geq 0 : y_n \notin D\}$. Then there exists $\eta_0 > 0$, depending only on $\phi$, $a$ and $h$, such that for all $\eta \in (0, \eta_0]$,
\begin{align*}
    \mathbb{E}[\sigma_D] \geq \exp(h/\eta).
\end{align*}
\end{lemma}

\begin{proof}
The proof compares the growth of the Lyapunov function $W(y) := \exp\big((\phi(y) - \phi(a))/\eta\big)$ before the exit time with its value at the exit time. Note that $W(y_0) = W(a) = 1$, and assume $L\eta \leq 1/2$. We first bound the expected one-step change of $W$ from a point $y \in D$. Let $g := \phi'(y)$ and $y' := y - \eta g + \eta\xi$ with $\xi \sim \mathcal{N}(0,1)$. Since $\phi'$ is $L$-Lipschitz,
\begin{align*}
    \phi(y') - \phi(y) \leq \eta\,(-g^2 + g\xi) + \frac{L\eta^2}{2}(\xi - g)^2,
\end{align*}
so that $W(y')/W(y) \leq \exp\big(-g^2 + g\xi + \frac{L\eta}{2}(\xi - g)^2\big)$. For $L\eta < 1$, a Gaussian integral gives
\begin{align}
    \mathbb{E}\Big[\exp\Big(g\xi + \frac{L\eta}{2}(\xi - g)^2\Big)\Big] &= (1 - L\eta)^{-1/2}\exp(g^2/2),\nonumber\\
\intertext{and taking expectations yields}
    \mathbb{E}[W(y')] &\leq (1 - L\eta)^{-1/2}\exp(-g^2/2)\,W(y)\nonumber\\
    &\leq \exp(L\eta - g^2/2)\,W(y),
\label{eq:lyapunov-drift}
\end{align}
where the second inequality uses $L\eta \leq 1/2$.

We now turn \eqref{eq:lyapunov-drift} into a bound on the one-step change of $W$ that is uniform over $y \in D$. If $g^2 \geq 2L\eta$, then \eqref{eq:lyapunov-drift} gives $\mathbb{E}[W(y')] \leq W(y)$.

Otherwise $g^2 < 2L\eta$, and $W(y)$ is bounded because $y$ is close to $a$. Indeed, since $a$ is the only critical point of $\phi$ in $\bar{D}$ and $\phi''(a) > 0$, the quotient $\phi'(y)/(y - a)$, extended by $\phi''(a)$ at $y = a$, is continuous and positive on the compact interval $\bar{D}$, so
\begin{align*}
    c_D &:= \min_{y \in \bar{D}} \frac{\phi'(y)}{y - a} > 0.\\
\intertext{Hence}
    |y - a| &\leq |g|/c_D \leq \sqrt{2L\eta}/c_D.
\end{align*}
Since $\phi'(a) = 0$ and $\phi'$ is $L$-Lipschitz, $\phi(y) - \phi(a) \leq L(y - a)^2/2$, and therefore $W(y) \leq \exp(L^2/c_D^2)$. Together with \eqref{eq:lyapunov-drift} and $e^{L\eta} - 1 \leq 2L\eta$, this gives
\begin{align*}
    \mathbb{E}[W(y')] - W(y) &\leq 2L\eta\exp(L^2/c_D^2).\\
\intertext{In both cases,}
    \mathbb{E}[W(y')] - W(y) &\leq C_1\eta \quad \text{for all } y \in D,\\
    C_1 &:= 2L\exp(L^2/c_D^2).
\end{align*}

For $n < \sigma_D$ we have $y_n \in D$, and $\xi_n$ is independent of $y_0, \ldots, y_n$. Applying the bound $\mathbb{E}[W(y')] - W(y) \leq C_1\eta$ with $y = y_n$ and $\xi = \xi_n$ therefore shows that the stopped process
\begin{align*}
    W(y_{n \wedge \sigma_D}) - C_1\eta\,(n \wedge \sigma_D)
\end{align*}
is a supermartingale. As $W(y_0) = 1$, it follows that
\begin{equation}\label{eq:lyapunov-stopped}
    \mathbb{E}\big[W(y_{n \wedge \sigma_D})\big] \leq 1 + C_1\eta\,n \qquad \text{for all } n \geq 0.
\end{equation}

We next show that $W(y_{\sigma_D}) \geq \exp(h/\eta)$ unless a large noise value occurs before the exit. Write $D = (u_-, u_+)$. The endpoints satisfy $\phi(u_\pm) = \phi(a) + h$, and $\phi'(u_\pm) \neq 0$ because $a$ is the only critical point of $\phi$ in $\bar{D}$. As $\phi < \phi(a) + h$ on $D$, this gives $\phi'(u_+) > 0$ and $\phi'(u_-) < 0$, so there is $r > 0$ such that $\phi \geq \phi(a) + h$ on $[u_+, u_+ + r]$ and on $[u_- - r, u_-]$.

To bound the length of a single step, let $G := \max_{y \in \bar{D}} |\phi'(y)|$ and assume in addition that $\eta \leq r/(2G)$. If $y_k \in D$ and $|\xi_k| < r/(2\eta)$, then
\begin{align*}
    |y_{k+1} - y_k| \leq \eta G + \eta|\xi_k| < r,
\end{align*}
so $y_{k+1} \notin D$ implies $\phi(y_{k+1}) \geq \phi(a) + h$, that is, $W(y_{k+1}) \geq \exp(h/\eta)$.

To control the noise at the last step before the exit, let $E_n$ be the event that $|\xi_k| < r/(2\eta)$ for all $k < n$. On $E_n \cap \{\sigma_D \leq n\}$, the last step before the exit starts at $y_{\sigma_D - 1} \in D$ and satisfies $|\xi_{\sigma_D - 1}| < r/(2\eta)$, so
\begin{align*}
    W(y_{n \wedge \sigma_D}) = W(y_{\sigma_D}) \geq \exp(h/\eta).
\end{align*}

Combining the bound $W(y_{n \wedge \sigma_D}) \geq \exp(h/\eta)$ on $E_n \cap \{\sigma_D \leq n\}$ with \eqref{eq:lyapunov-stopped} controls the probability that $\sigma_D \leq n$. Markov's inequality, a union bound over $k < n$, and the Gaussian tail bound $\mathbb{P}(|\xi_k| \geq t) \leq 2e^{-t^2/2}$ give
\begin{align*}
    \mathbb{P}(\sigma_D \leq n)
    &\leq \mathbb{P}\big(W(y_{n \wedge \sigma_D}) \geq e^{h/\eta}\big) + \mathbb{P}(E_n^c) \\
    &\leq (1 + C_1\eta\,n)\,e^{-h/\eta} + 2n\,e^{-r^2/(8\eta^2)}.
\end{align*}
Take $n = \lceil 2\exp(h/\eta) \rceil$. Both terms on the right-hand side then tend to $0$ as $\eta \to 0$, so there is $\eta_0 \in (0, \min\{1/(2L),\, r/(2G)\}]$, depending only on $\phi$, $a$ and $h$, such that $\mathbb{P}(\sigma_D \leq n) \leq 1/2$ for all $\eta \in (0, \eta_0]$. For these $\eta$,
\begin{align*}
    \mathbb{E}[\sigma_D] \geq n\,\mathbb{P}(\sigma_D > n) \geq \frac{n}{2} \geq \exp(h/\eta).
\end{align*}
\end{proof}
}

We now apply \cref{lem:exit-time} to one coordinate of fixed-step-size SGD on a separable objective.

\begin{theorem}[Lower bound for fixed-step-size SGD]
\label{thm:lower-bound-main}
Let $f(x) = \sum_{i=1}^d \phi(x_i)$ where $\phi:\mathbb{R} \to \mathbb{R}$ is a coercive $C^2$ Morse function with an $L$-Lipschitz derivative and exactly three critical points: a unique global minimizer $x^*$ satisfying $\phi(x^*) = 0$ and $\phi''(x^*) = \omega > 0$, a spurious local minimizer $a$ with $\phi(a) = c > 0$, and a local maximizer $s$ between $a$ and $x^*$ with barrier height $\Delta := \phi(s) - \phi(a) > 0$. Fix $\epsilon \in (0, c)$ and $C > 0$. Consider fixed-step-size SGD (\cref{def:fixed-sgd-intro}) with any learning rate $\eta \in (0, C\epsilon/d]$, initialized at $x_0 = (a, \ldots, a)$. Then for $\tau = \inf\{t \geq 0 : f(x_t) - f^* < \epsilon\}$, there exists $d_0$, depending only on $\phi$, $\epsilon$ and $C$, such that for all $d \geq d_0$,
{
\begin{align*}
    \mathbb{E}[\tau] \geq \exp\!\left(\frac{\Delta\,d}{2C\epsilon}\right) = \exp(\Omega(d)).
\end{align*}
}
\end{theorem}

\begin{remark}
The step-size regime $\eta \leq C\epsilon/d$ of \cref{thm:lower-bound-main} is the regime $\eta = \Ocal(\epsilon/d)$ of \cref{def:fixed-sgd-intro}. The regime $\eta = \Ocal(\epsilon/d)$ can be motivated by the quadratic approximation of $\phi$ near the global minimizer. Within the basin of $x^*$, the per-coordinate SGD dynamics linearize as $z_{t+1} = (1 - \eta\omega)\,z_t + \eta\,\xi_{t,i}$, where $z_t = x_{t,i} - x^*$ and $\xi_{t,i}$ is the $i$-th coordinate of $\xi_t$. This recursion is exact for the quadratic $\phi(x) = \omega(x - x^*)^2/2$. In this quadratic model, $z_t$ is an autoregressive process of order 1, which is stable for $\eta < 2/\omega$ and has stationary variance
{
\begin{align*}
    \sigma_{\mathrm{stat}}^2 = \frac{\eta}{2\omega(1 - \eta\omega/2)}.
\end{align*}
}
The expected suboptimality of a point $X$ drawn from the stationary measure of the quadratic model is therefore
{
\begin{align*}
    \mathbb{E}[f(X) - f^*] = \sum_{i=1}^d \frac{\omega}{2}\,\sigma_{\mathrm{stat}}^2
    = \frac{d\,\eta}{4(1 - \eta\omega/2)},
\end{align*}
}
which is of order $\epsilon$ when $\eta$ is of order $\epsilon/d$.
\end{remark}

\begin{proof}[Proof of \cref{thm:lower-bound-main}]
By separability, the coordinates of $x_t$ evolve independently, and the first coordinate $y_t := x_{t,1}$ follows the recursion
{
\begin{align*}
    y_{t+1} = y_t - \eta\,\phi'(y_t) + \eta\,\xi_{t,1}, \qquad y_0 = a,
\end{align*}
}
where the first coordinates $\xi_{0,1}, \xi_{1,1}, \ldots$ of the noise vectors are i.i.d.\ $\mathcal{N}(0,1)$. This is the recursion of \cref{lem:exit-time}. Since $f^* = d\,\phi(x^*) = 0$ and all terms of $f(x) - f^* = \sum_{i=1}^d \big(\phi(x_i) - \phi(x^*)\big)$ are nonnegative, every $\epsilon$-optimal point $x$ satisfies $\phi(x_i) - \phi(x^*) < \epsilon$ for all $i$. In particular, $f(x_t) - f^* < \epsilon$ implies $\phi(y_t) < \epsilon$.

We next find a neighborhood of $a$ on which $\phi \geq c$. Let $D$ be the connected component of $\{y \in \mathbb{R} : \phi(y) < \phi(a) + \Delta/2\}$ that contains $a$, and let $\bar{D}$ be its closure. Because $\phi$ is coercive, $D$ is a bounded interval. On $\bar{D}$ we have $\phi \leq \phi(a) + \Delta/2 < \phi(s)$, so $\bar{D}$ does not contain $s$. As $s$ lies between $a$ and $x^*$, the interval $\bar{D}$ does not contain $x^*$ either. Hence $a$ is the only critical point of $\phi$ in $\bar{D}$, and since $\phi = \phi(a) + \Delta/2$ at the endpoints of $D$, the minimum of $\phi$ over $\bar{D}$ is $\phi(a) = c$. Thus $\phi \geq c > \epsilon$ on $D$, so $x_t$ is not $\epsilon$-optimal for $t < \sigma_D := \min\{t \geq 0 : y_t \notin D\}$, and $\tau \geq \sigma_D$.

Since $\phi$ is a Morse function, $\phi''(a) > 0$. Hence \cref{lem:exit-time} applies to $D$ with $h = \Delta/2$ and gives $\eta_0 > 0$, depending only on $\phi$ and not on $d$, such that $\mathbb{E}[\sigma_D] \geq \exp(\Delta/(2\eta))$ for all $\eta \in (0, \eta_0]$. Since $\eta \leq C\epsilon/d \to 0$ as $d \to \infty$, the condition $\eta \leq \eta_0$ holds for all sufficiently large $d$, namely for $d \geq d_0 := \lceil C\epsilon/\eta_0 \rceil$. For such $d$,
{
\begin{align*}
    \mathbb{E}[\tau] \geq \mathbb{E}[\sigma_D] \geq \exp\!\left(\frac{\Delta}{2\eta}\right)
    \geq \exp\!\left(\frac{\Delta\,d}{2C\epsilon}\right) = \exp(\Omega(d)).
\end{align*}
}
\end{proof}

{
\begin{corollary}[Separation from the same starting point]
\label{cor:sgd-common-initialization}
Let $\phi$ satisfy the hypotheses of \cref{thm:lower-bound-main}, and suppose additionally that
\[
\phi(y)\geq\frac m2|y-x^*|^2\qquad(y\in\mathbb R)
\]
for a fixed $m>0$ and that \cref{ass:sgd-spectral-transfer} holds for its SGD kernels with $q=2$ and constants independent of $d$. For an arbitrary orthogonal matrix $R$, define
\[
f_d(x):=\sum_{i=1}^d\phi((Rx)_i),\qquad
x_0:=R^\top(a,\ldots,a).
\]
Fix $0<\epsilon\leq1$ with $\epsilon<c$, and fix $C_{\mathrm{lr}}>0$. Both algorithms receive the same point $x_0$ as input. For all sufficiently large $d$:
\begin{enumerate}[leftmargin=*]
\item Fixed-step-size SGD with any $0<\eta\leq C_{\mathrm{lr}}\epsilon/d$ has an expected first hitting time of the $\epsilon$-optimal set of at least
\[
\exp\left(\frac{\Delta d}{2C_{\mathrm{lr}}\epsilon}\right).
\]
\item Annealed SGD, including its first update at the constant learning rate $\eta_1$, uses $\mathrm{poly}(d/\epsilon)$ gradient queries and outputs an $\epsilon$-optimal point with probability at least $2/3$.
\end{enumerate}
\end{corollary}
\begin{proof}
The rotation $x\mapsto Rx$ leaves SGD invariant, since $RZ$ is standard Gaussian whenever $Z$ is standard Gaussian. It also preserves objective values and maps the starting point $x_0$ to $(a,\ldots,a)$. Thus \cref{thm:lower-bound-main} gives the lower bound for fixed-step-size SGD.

For the upper bound, continuity of $\phi''$ and $\phi''(x^*)>0$ give a radius $r_0>0$ such that $\phi''(y)\geq\xi:=\phi''(x^*)/2$ whenever $|y-x^*|\leq r_0$. Since every block is the same fixed potential $\phi$, this radius and curvature are independent of $d$ and satisfy \eqref{eq:sgd-uniform-local-convexity}. The repeated blocks $g_i=\phi$ therefore have dimension one, fixed regularity constants, and initial energies
\[
g_i(a)-g_i(x^*)=c,\qquad G_0=c.
\]
The density bound after the first SGD update is therefore independent of $d$ and $\epsilon$ by \cref{lem:sgd-point-initialization}. Applying \cref{thm:global_opt_sep} gives the claimed success probability and polynomial query count from the same point $x_0$.
\end{proof}
}

\section{Additional Details for Numerical Results}
\label{app:numeric-details}

\paragraph{Common setup.}
The experiments in \cref{fig:dimension,tab:high-dim-sgd-arithmetic} use SGD with independent standard Gaussian noise:
\[
    x_{k+1}=x_k-\eta_k\nabla f(x_k)+\eta_k Z_k,
    \qquad Z_k\sim\mathcal N(0,I_d).
\]
All runs start from $X_0=c_f\mathbf{1}_d+0.05\,\xi$, where $\xi\sim\mathcal{N}(0,I_d)$, $c_f\simeq0.9601496$ for both double wells and $c_f\simeq0.9510504$ for Rastrigin.
After each update, we check whether $\|x_k-x^\star\|_\infty<\delta$, where $\delta=0.5$ for the double-well functions and $\delta=0.3$ for Rastrigin. We stop the run as soon as this condition is satisfied and record the number of updates performed.

\paragraph{Arithmetic annealing.}
An annealing schedule specifies $M$ learning-rate levels:
\[
    \eta_j=\eta_0+\frac{j}{M-1}(\eta_{\mathrm{end}}-\eta_0),
    \qquad j=0,\ldots,M-1.
\]
Normally each level is held constant for $N_p$ updates before decreasing to the next.

\paragraph{Low-dimensional experiments (Figure~\ref{fig:dimension}).}
For the low-dimensional grid $d\in\{1,5,10,15,20,25,30,35\}$, we use ten trials per method and dimension.
The annealed settings are:
\begin{itemize}[leftmargin=7mm]
    \item \textbf{Separable double well:} $\eta_0=0.18/(1+0.08\log d)$, $\eta_{\mathrm{end}}=0.145/[1+(d-1)/40]$, $M=\max(3d,15)$, and $N_p=10^5$.
    \item \textbf{Rastrigin:} $\eta_0=0.10$, $\eta_{\mathrm{end}}=0.064/[1+(d-1)/200]$, $M=\max(3d,15)$, and $N_p=20$.
    \item \textbf{Coupled double well:} $\eta_0=0.155$, $\eta_{\mathrm{end}}=0.150/[1+(d-1)/140]$, $M=\max(3d,20)$, and $N_p=100$.
\end{itemize}
The Rastrigin schedule allows $\max(2000,1000d)$ additional updates at the last annealing level with $\eta_j\geq0.065$ and $20\lceil\max(1000,50d)/20\rceil$ additional updates at the final rate, if the stopping criterion has not yet been met.
Coupled double well allows $\max(50000,5000d)$ additional updates at the last level with $\eta_j\geq0.15$ and $2\times10^6$ additional terminal-rate updates. These allowances are used only while the stopping criterion remains unmet.

Fixed-step-size SGD keeps the corresponding terminal rate throughout the run. These terminal-rate rules were selected in pilot experiments and then locked; fixed-step-size SGD was not independently optimized for the fastest constant rate at each dimension.
The cap for fixed-step-size SGD runs is $10^9$ updates except for the separable double well at $d=30,35$ ($2\times10^{10}$ and $5\times10^{10}$) and the coupled double well at $d=35$ ($10^{10}$).

\paragraph{High-dimensional experiments (Table~\ref{tab:high-dim-sgd-arithmetic}).}
For the high-dimensional grid $d\in\{100,200,300,400\}$, we run ten annealed trials per function and dimension using the same learning-rate formulas, phase lengths, and additional holds as above, with $M=3d$.

\paragraph{Reporting and run outcomes.}
The table reports medians $\pm$ one empirical standard deviation computed with divisor 10 over all ten recorded update counts and run times, with update counts truncated at the limits recorded with the data. Timing starts after the initial state is generated and excludes JIT warmup and queueing.
Medians and interquartile ranges use all ten planned trials when identified; groups containing nonfinite outcomes are excluded from fitting.

\end{document}